
\mag1200
\input amstex

\expandafter\ifx\csname beta.def\endcsname\relax \else\endinput\fi
\expandafter\edef\csname beta.def\endcsname{%
 \catcode`\noexpand\@=\the\catcode`\@\space}

\let\atbefore @

\catcode`\@=11

\overfullrule\z@

\def\PaperA4{\hsize 6.25truein \vsize 9.63truein}

\def\PaperUS{\hsize 6.6truein \vsize 9truein} 

\def\foliorm{\ifMag\eightrm\else\ninerm\fi}

\let\@ft@\expandafter \let\@tb@f@\atbefore

\newif\ifMag
\def\Magset{\ifnum\mag>\@m\Magtrue\fi}
\Magset

\newif\ifUS

\newdimen\p@@ \p@@\p@
\def\m@ths@r{\ifnum\mathsurround=\z@\z@\else\maths@r\fi}
\def\maths@r{1.6\p@@} \def\mathsurzero{\def\maths@r{\z@}}

\mathsurround\maths@r
\font\Brm=cmr12 \font\Bbf=cmbx12 \font\Bit=cmti12 \font\ssf=cmss10
\font\Bsl=cmsl10 scaled 1200 \font\Bmmi=cmmi10 scaled 1200
\font\BBf=cmbx12 scaled 1200 \font\BMmi=cmmi10 scaled 1440

\def\atletter{\edef\atrestore{\catcode`\noexpand\@=\the\catcode`\@\space}
 \catcode`\@=11}

\newread\@ux \newwrite\@@x \newwrite\@@cd
\let\@np@@\input
\def\@np@t#1{\openin\@ux#1\relax\ifeof\@ux\else\closein\@ux\relax\@np@@ #1\fi}
\def\input#1 {\openin\@ux#1\relax\ifeof\@ux\wrs@x{No file #1}\else
 \closein\@ux\relax\@np@@ #1\fi}
\def\Input#1 {\relax} 

\def\wr@@x#1{} \def\wrs@x{\immediate\write\sixt@@n}

\def\readldf{\@np@t{\jobname.ldf}}
\def\writeldf{\def\wr@@x{\immediate\write\@@x}\def\wr@x@{\write\@@x}
 \def\cl@selbl{\wr@@x{\string\def\string\nextpage{\the\pageno}}%
 \wr@@x{\string\endinput}\immediate\closeout\@@x}
 \immediate\openout\@@x\jobname.ldf}
\let\cl@selbl\relax

\def\nextpage{1}

\def\tod@y{\ifcase\month\or
 January\or February\or March\or April\or May\or June\or July\or
 August\or September\or October\or November\or December\fi\space\,
\number\day,\space\,\number\year}

\newcount\c@time
\def\h@@r{hh}\def\m@n@te{mm}
\def\wh@tt@me{\c@time\time\divide\c@time 60\edef\h@@r{\number\c@time}%
 \multiply\c@time -60\advance\c@time\time\edef
 \m@n@te{\ifnum\c@time<10 0\fi\number\c@time}}
\def\t@me{\h@@r\/{\rm:}\m@n@te} \let\whattime\wh@tt@me
\def\today{\tod@y\wr@@x{\string\todaydef{\tod@y}}}
\def\nowtime{\t@me{\let\/\ic@\wr@@x{\string\nowtimedef{\t@me}}}}
\def\todaydef#1{} \def\nowtimedef#1{}

\def\em#1{{\it #1\/}} \def\emph#1{{\sl #1\/}}

\def\fitem#1{\par\setbox\z@\hbox{#1}\hangindent\wd\z@
 \hglue-2\parindent\kern\wd\z@\indent\llap{#1}\ignore}

\def\itemflat#1{\par\setbox\z@\hbox{\rm #1\enspace}\hang\ifnum\wd\z@>\parindent
 \noindent\unhbox\z@\ignore\else\textindent{\rm#1}\fi}

\newcount\itemlet
\def\newbi{\itemlet 96} \newbi
\def\bitem{\gad\itemlet \par\hangindent1.5\parindent
 \hglue-.5\parindent\textindent{\rm\rlap{\char\the\itemlet}\hp{b})}}
\def\atem{\newbi\bitem}

\newcount\itemrm
\def\newitrm{\itemrm 0}
\def\iitem{\gad\itemrm \par\hangindent1.5\parindent
 \hglue-.5\parindent\textindent{\rm\hp{v}\llap{\romannumeral\the\itemrm})}}
\def\Item{\newitrm\iitem}

\newcount\itemar

\def\iitema{\gad\itemrm \par\hangindent1.5\parindent
 \hglue-.5\parindent\textindent{\rm\hp{0}\llap{\the\itemrm}.}}

\def\center{\par\begingroup\leftskip\z@ plus \hsize \rightskip\leftskip
 \parindent\z@\parfillskip\z@skip \def\\{\unskip\break}}
\def\endcenter{\endgraf\endgroup}

\def\Abstract{\begingroup\narrower\nt{\bf Abstract.}\enspace\ignore}
\def\endAbs{\endgraf\endgroup}

\let\b@gr@@\begingroup \let\B@gr@@\begingroup
\def\b@gr@{\b@gr@@\let\b@gr@@\undefined}
\def\B@gr@{\B@gr@@\let\B@gr@@\undefined}

\def\@fn@xt#1#2#3{\let\@ch@r=#1\def\n@xt{\ifx\t@st@\@ch@r
 \def\n@@xt{#2}\else\def\n@@xt{#3}\fi\n@@xt}\futurelet\t@st@\n@xt}

\def\@fwd@@#1#2#3{\setbox\z@\hbox{#1}\ifdim\wd\z@>\z@#2\else#3\fi}
\def\s@twd@#1#2{\setbox\z@\hbox{#2}#1\wd\z@}

\def\r@st@re#1{\let#1\s@v@} \def\s@v@d@f{\let\s@v@}

\def\p@sk@p#1#2{\par\skip@#2\relax\ifdim\lastskip<\skip@\relax\removelastskip
 \ifnum#1=\z@\else\penalty#1\relax\fi\vskip\skip@
 \else\ifnum#1=\z@\else\penalty#1\relax\fi\fi}
\def\sk@@p#1{\par\skip@#1\relax\ifdim\lastskip<\skip@\relax\removelastskip
 \vskip\skip@\fi}

\newbox\p@b@ld
\def\poorbold#1{\setbox\p@b@ld\hbox{#1}\kern-.01em\copy\p@b@ld\kern-\wd\p@b@ld
 \kern.02em\copy\p@b@ld\kern-\wd\p@b@ld\kern-.012em\raise.02em\box\p@b@ld}

\ifx\plainfootnote\undefined \let\plainfootnote\footnote \fi

\let\s@v@\proclaim \let\proclaim\relax
\def\r@R@fs#1{\let#1\s@R@fs} \let\s@R@fs\Refs \let\Refs\relax
\def\r@endd@#1{\let#1\s@endd@} \let\s@endd@\enddocument
\let\bye\relax

\def\myR@fs{\@fn@xt[\m@R@f@\m@R@fs} \def\m@R@fs{\@fn@xt*\m@r@f@@\m@R@f@@}
\def\m@R@f@@{\m@R@f@[References]} \def\m@r@f@@*{\m@R@f@[]}

\def\Twelvepoint{\twelvepoint \let\Bbf\BBf \let\Bmmi\BMmi
\font\Brm=cmr12 scaled 1200 \font\Bit=cmti12 scaled 1200
\font\ssf=cmss10 scaled 1200 \font\Bsl=cmsl10 scaled 1440
\font\BBf=cmbx12 scaled 1440 \font\BMmi=cmmi10 scaled 1728}

\newdimen\b@gsize

\newdimen\r@f@nd \newbox\r@f@b@x \newbox\adjb@x
\newbox\p@nct@ \newbox\k@yb@x \newcount\rcount
\newbox\b@b@x \newbox\p@p@rb@x \newbox\j@@rb@x \newbox\y@@rb@x
\newbox\v@lb@x \newbox\is@b@x \newbox\p@g@b@x \newif\ifp@g@ \newif\ifp@g@s
\newbox\inb@@kb@x \newbox\b@@kb@x \newbox\p@blb@x \newbox\p@bl@db@x
\newbox\ed@b@x \newif\ifed@ \newif\ifed@s \newif\if@fl@b \newif\if@fn@m
\newbox\p@p@nf@b@x \newbox\inf@b@x \newbox\b@@nf@b@x
\newtoks\@dd@p@n \newtoks\@ddt@ks

\newif\ifp@gen@

\def\p@@nt{.\kern.3em} \let\point\p@@nt

\let\proheadfont\bf \let\probodyfont\sl \let\demofont\it

\headline={\hfil}
\footline={\ifp@gen@\ifnum\pageno=\z@\else\hfil\foliorm\folio\fi\else
 \ifnum\pageno=\z@\hfil\foliorm\folio\fi\fi\hfil\global\p@gen@true}
\parindent1pc

\font@\tensmc=cmcsc10
\font@\sevenex=cmex7
\font@\sevenit=cmti7
\font@\eightrm=cmr8
\font@\sixrm=cmr6
\font@\eighti=cmmi8 \skewchar\eighti='177
\font@\sixi=cmmi6 \skewchar\sixi='177
\font@\eightsy=cmsy8 \skewchar\eightsy='60
\font@\sixsy=cmsy6 \skewchar\sixsy='60
\font@\eightex=cmex8
\font@\eightbf=cmbx8
\font@\sixbf=cmbx6
\font@\eightit=cmti8
\font@\eightsl=cmsl8
\font@\eightsmc=cmcsc8
\font@\eighttt=cmtt8
\font@\ninerm=cmr9
\font@\ninei=cmmi9 \skewchar\ninei='177
\font@\ninesy=cmsy9 \skewchar\ninesy='60
\font@\nineex=cmex9
\font@\ninebf=cmbx9
\font@\nineit=cmti9
\font@\ninesl=cmsl9
\font@\ninesmc=cmcsc9
\font@\ninemsa=msam9
\font@\ninemsb=msbm9
\font@\nineeufm=eufm9
\font@\eightmsa=msam8
\font@\eightmsb=msbm8
\font@\eighteufm=eufm8
\font@\sixmsa=msam6
\font@\sixmsb=msbm6
\font@\sixeufm=eufm6

\loadmsam\loadmsbm\loadeufm
\input amssym.tex

\def\footnoterule{\kern-3\p@\hrule width5pc\kern 2.6\p@}
\def\m@k@foot#1{\insert\footins
 {\interlinepenalty\interfootnotelinepenalty
 \ifMag\eightpoint\else\ninepoint\fi
 \splittopskip\ht\strutbox\splitmaxdepth\dp\strutbox
 \floatingpenalty\@MM\leftskip\z@\rightskip\z@
 \spaceskip\z@\xspaceskip\z@
 \leavevmode\footstrut\ignore#1\unskip\lower\dp\strutbox
 \vbox to\dp\strutbox{}}}
\def\ftext#1{\m@k@foot{\vsk-.8>\nt #1}}
\def\pr@cl@@m#1{\p@sk@p{-100}\medskipamount
 \def\endproclaim{\endgroup\p@sk@p{55}\medskipamount}\begingroup
 \nt\ignore\proheadfont#1\unskip.\enspace\probodyfont\ignore}
\outer\def\proclaim{\pr@cl@@m} \s@v@d@f\proclaim \let\proclaim\relax
\def\demo#1{\sk@@p\medskipamount\nt{\ignore\demofont#1\unskip.}\enspace
 \ignore}
\def\enddemo{\sk@@p\medskipamount}

\def\cite#1{{\rm[#1]}} 
 \def\Refs#1#2{\relax}

\def\big@#1#2{{\hbox{$\left#2\vcenter to#1\b@gsize{}%
 \right.\nulldelimiterspace\z@\m@th$}}}
\def\big{\big@\@ne}
\def\Big{\big@{1.5}}
\def\bigg{\big@\tw@}
\def\Bigg{\big@{2.5}}
\normallineskiplimit\p@

\def\tenpoint{\p@@\p@ \normallineskiplimit\p@@
 \mathsurround\m@ths@r \normalbaselineskip12\p@@
 \abovedisplayskip12\p@@ plus3\p@@ minus9\p@@
 \belowdisplayskip\abovedisplayskip
 \abovedisplayshortskip\z@ plus3\p@@
 \belowdisplayshortskip7\p@@ plus3\p@@ minus4\p@@
 \textonlyfont@\rm\tenrm \textonlyfont@\it\tenit
 \textonlyfont@\sl\tensl \textonlyfont@\bf\tenbf
 \textonlyfont@\smc\tensmc \textonlyfont@\tt\tentt
 \ifsyntax@ \def\big##1{{\hbox{$\left##1\right.$}}}%
  \let\Big\big \let\bigg\big \let\Bigg\big
 \else
  \textfont\z@\tenrm \scriptfont\z@\sevenrm \scriptscriptfont\z@\fiverm
  \textfont\@ne\teni \scriptfont\@ne\seveni \scriptscriptfont\@ne\fivei
  \textfont\tw@\tensy \scriptfont\tw@\sevensy \scriptscriptfont\tw@\fivesy
  \textfont\thr@@\tenex \scriptfont\thr@@\sevenex
	\scriptscriptfont\thr@@\sevenex
  \textfont\itfam\tenit \scriptfont\itfam\sevenit
	\scriptscriptfont\itfam\sevenit
  \textfont\bffam\tenbf \scriptfont\bffam\sevenbf
	\scriptscriptfont\bffam\fivebf
  \textfont\msafam\tenmsa \scriptfont\msafam\sevenmsa
	\scriptscriptfont\msafam\fivemsa
  \textfont\msbfam\tenmsb \scriptfont\msbfam\sevenmsb
	\scriptscriptfont\msbfam\fivemsb
  \textfont\eufmfam\teneufm \scriptfont\eufmfam\seveneufm
	\scriptscriptfont\eufmfam\fiveeufm
  \setbox\strutbox\hbox{\vrule height8.5\p@@ depth3.5\p@@ width\z@}%
  \setbox\strutbox@\hbox{\lower.5\normallineskiplimit\vbox{%
	\kern-\normallineskiplimit\copy\strutbox}}%
   \setbox\z@\vbox{\hbox{$($}\kern\z@}\b@gsize1.2\ht\z@
  \fi
  \normalbaselines\rm\dotsspace@1.5mu\ex@.2326ex\jot3\ex@}

\def\eightpoint{\p@@.8\p@ \normallineskiplimit\p@@
 \mathsurround\m@ths@r \normalbaselineskip10\p@
 \abovedisplayskip10\p@ plus2.4\p@ minus7.2\p@
 \belowdisplayskip\abovedisplayskip
 \abovedisplayshortskip\z@ plus3\p@@
 \belowdisplayshortskip7\p@@ plus3\p@@ minus4\p@@
 \textonlyfont@\rm\eightrm \textonlyfont@\it\eightit
 \textonlyfont@\sl\eightsl \textonlyfont@\bf\eightbf
 \textonlyfont@\smc\eightsmc \textonlyfont@\tt\eighttt
 \ifsyntax@\def\big##1{{\hbox{$\left##1\right.$}}}%
  \let\Big\big \let\bigg\big \let\Bigg\big
 \else
  \textfont\z@\eightrm \scriptfont\z@\sixrm \scriptscriptfont\z@\fiverm
  \textfont\@ne\eighti \scriptfont\@ne\sixi \scriptscriptfont\@ne\fivei
  \textfont\tw@\eightsy \scriptfont\tw@\sixsy \scriptscriptfont\tw@\fivesy
  \textfont\thr@@\eightex \scriptfont\thr@@\sevenex
	\scriptscriptfont\thr@@\sevenex
  \textfont\itfam\eightit \scriptfont\itfam\sevenit
	\scriptscriptfont\itfam\sevenit
  \textfont\bffam\eightbf \scriptfont\bffam\sixbf
	\scriptscriptfont\bffam\fivebf
  \textfont\msafam\eightmsa \scriptfont\msafam\sixmsa
	\scriptscriptfont\msafam\fivemsa
  \textfont\msbfam\eightmsb \scriptfont\msbfam\sixmsb
	\scriptscriptfont\msbfam\fivemsb
  \textfont\eufmfam\eighteufm \scriptfont\eufmfam\sixeufm
	\scriptscriptfont\eufmfam\fiveeufm
 \setbox\strutbox\hbox{\vrule height7\p@ depth3\p@ width\z@}%
 \setbox\strutbox@\hbox{\raise.5\normallineskiplimit\vbox{%
   \kern-\normallineskiplimit\copy\strutbox}}%
 \setbox\z@\vbox{\hbox{$($}\kern\z@}\b@gsize1.2\ht\z@
 \fi
 \normalbaselines\eightrm\dotsspace@1.5mu\ex@.2326ex\jot3\ex@}

\def\ninepoint{\p@@.9\p@ \normallineskiplimit\p@@
 \mathsurround\m@ths@r \normalbaselineskip11\p@
 \abovedisplayskip11\p@ plus2.7\p@ minus8.1\p@
 \belowdisplayskip\abovedisplayskip
 \abovedisplayshortskip\z@ plus3\p@@
 \belowdisplayshortskip7\p@@ plus3\p@@ minus4\p@@
 \textonlyfont@\rm\ninerm \textonlyfont@\it\nineit
 \textonlyfont@\sl\ninesl \textonlyfont@\bf\ninebf
 \textonlyfont@\smc\ninesmc \textonlyfont@\tt\ninett
 \ifsyntax@ \def\big##1{{\hbox{$\left##1\right.$}}}%
  \let\Big\big \let\bigg\big \let\Bigg\big
 \else
  \textfont\z@\ninerm \scriptfont\z@\sevenrm \scriptscriptfont\z@\fiverm
  \textfont\@ne\ninei \scriptfont\@ne\seveni \scriptscriptfont\@ne\fivei
  \textfont\tw@\ninesy \scriptfont\tw@\sevensy \scriptscriptfont\tw@\fivesy
  \textfont\thr@@\nineex \scriptfont\thr@@\sevenex
	\scriptscriptfont\thr@@\sevenex
  \textfont\itfam\nineit \scriptfont\itfam\sevenit
	\scriptscriptfont\itfam\sevenit
  \textfont\bffam\ninebf \scriptfont\bffam\sevenbf
	\scriptscriptfont\bffam\fivebf
  \textfont\msafam\ninemsa \scriptfont\msafam\sevenmsa
	\scriptscriptfont\msafam\fivemsa
  \textfont\msbfam\ninemsb \scriptfont\msbfam\sevenmsb
	\scriptscriptfont\msbfam\fivemsb
  \textfont\eufmfam\nineeufm \scriptfont\eufmfam\seveneufm
	\scriptscriptfont\eufmfam\fiveeufm
  \setbox\strutbox\hbox{\vrule height8.5\p@@ depth3.5\p@@ width\z@}%
  \setbox\strutbox@\hbox{\lower.5\normallineskiplimit\vbox{%
	\kern-\normallineskiplimit\copy\strutbox}}%
   \setbox\z@\vbox{\hbox{$($}\kern\z@}\b@gsize1.2\ht\z@
  \fi
  \normalbaselines\rm\dotsspace@1.5mu\ex@.2326ex\jot3\ex@}

\font@\twelverm=cmr10 scaled 1200
\font@\twelveit=cmti10 scaled 1200
\font@\twelvesl=cmsl10 scaled 1200
\font@\twelvebf=cmbx10 scaled 1200
\font@\twelvesmc=cmcsc10 scaled 1200
\font@\twelvett=cmtt10 scaled 1200
\font@\twelvei=cmmi10 scaled 1200 \skewchar\twelvei='177
\font@\twelvesy=cmsy10 scaled 1200 \skewchar\twelvesy='60
\font@\twelveex=cmex10 scaled 1200
\font@\twelvemsa=msam10 scaled 1200
\font@\twelvemsb=msbm10 scaled 1200
\font@\twelveeufm=eufm10 scaled 1200

\def\twelvepoint{\p@@1.2\p@ \normallineskiplimit\p@@
 \mathsurround\m@ths@r \normalbaselineskip12\p@@
 \abovedisplayskip12\p@@ plus3\p@@ minus9\p@@
 \belowdisplayskip\abovedisplayskip
 \abovedisplayshortskip\z@ plus3\p@@
 \belowdisplayshortskip7\p@@ plus3\p@@ minus4\p@@
 \textonlyfont@\rm\twelverm \textonlyfont@\it\twelveit
 \textonlyfont@\sl\twelvesl \textonlyfont@\bf\twelvebf
 \textonlyfont@\smc\twelvesmc \textonlyfont@\tt\twelvett
 \ifsyntax@ \def\big##1{{\hbox{$\left##1\right.$}}}%
  \let\Big\big \let\bigg\big \let\Bigg\big
 \else
  \textfont\z@\twelverm \scriptfont\z@\eightrm \scriptscriptfont\z@\sixrm
  \textfont\@ne\twelvei \scriptfont\@ne\eighti \scriptscriptfont\@ne\sixi
  \textfont\tw@\twelvesy \scriptfont\tw@\eightsy \scriptscriptfont\tw@\sixsy
  \textfont\thr@@\twelveex \scriptfont\thr@@\eightex
	\scriptscriptfont\thr@@\sevenex
  \textfont\itfam\twelveit \scriptfont\itfam\eightit
	\scriptscriptfont\itfam\sevenit
  \textfont\bffam\twelvebf \scriptfont\bffam\eightbf
	\scriptscriptfont\bffam\sixbf
  \textfont\msafam\twelvemsa \scriptfont\msafam\eightmsa
	\scriptscriptfont\msafam\sixmsa
  \textfont\msbfam\twelvemsb \scriptfont\msbfam\eightmsb
	\scriptscriptfont\msbfam\sixmsb
  \textfont\eufmfam\twelveeufm \scriptfont\eufmfam\eighteufm
	\scriptscriptfont\eufmfam\sixeufm
  \setbox\strutbox\hbox{\vrule height8.5\p@@ depth3.5\p@@ width\z@}%
  \setbox\strutbox@\hbox{\lower.5\normallineskiplimit\vbox{%
	\kern-\normallineskiplimit\copy\strutbox}}%
  \setbox\z@\vbox{\hbox{$($}\kern\z@}\b@gsize1.2\ht\z@
  \fi
  \normalbaselines\rm\dotsspace@1.5mu\ex@.2326ex\jot3\ex@}

\font@\twelvetrm=cmr10 at 12truept
\font@\twelvetit=cmti10 at 12truept
\font@\twelvetsl=cmsl10 at 12truept
\font@\twelvetbf=cmbx10 at 12truept
\font@\twelvetsmc=cmcsc10 at 12truept
\font@\twelvettt=cmtt10 at 12truept
\font@\twelveti=cmmi10 at 12truept \skewchar\twelveti='177
\font@\twelvetsy=cmsy10 at 12truept \skewchar\twelvetsy='60
\font@\twelvetex=cmex10 at 12truept
\font@\twelvetmsa=msam10 at 12truept
\font@\twelvetmsb=msbm10 at 12truept
\font@\twelveteufm=eufm10 at 12truept

\def\twelvetruepoint{\p@@1.2truept \normallineskiplimit\p@@
 \mathsurround\m@ths@r \normalbaselineskip12\p@@
 \abovedisplayskip12\p@@ plus3\p@@ minus9\p@@
 \belowdisplayskip\abovedisplayskip
 \abovedisplayshortskip\z@ plus3\p@@
 \belowdisplayshortskip7\p@@ plus3\p@@ minus4\p@@
 \textonlyfont@\rm\twelvetrm \textonlyfont@\it\twelvetit
 \textonlyfont@\sl\twelvetsl \textonlyfont@\bf\twelvetbf
 \textonlyfont@\smc\twelvetsmc \textonlyfont@\tt\twelvettt
 \ifsyntax@ \def\big##1{{\hbox{$\left##1\right.$}}}%
  \let\Big\big \let\bigg\big \let\Bigg\big
 \else
  \textfont\z@\twelvetrm \scriptfont\z@\eightrm \scriptscriptfont\z@\sixrm
  \textfont\@ne\twelveti \scriptfont\@ne\eighti \scriptscriptfont\@ne\sixi
  \textfont\tw@\twelvetsy \scriptfont\tw@\eightsy \scriptscriptfont\tw@\sixsy
  \textfont\thr@@\twelvetex \scriptfont\thr@@\eightex
	\scriptscriptfont\thr@@\sevenex
  \textfont\itfam\twelvetit \scriptfont\itfam\eightit
	\scriptscriptfont\itfam\sevenit
  \textfont\bffam\twelvetbf \scriptfont\bffam\eightbf
	\scriptscriptfont\bffam\sixbf
  \textfont\msafam\twelvetmsa \scriptfont\msafam\eightmsa
	\scriptscriptfont\msafam\sixmsa
  \textfont\msbfam\twelvetmsb \scriptfont\msbfam\eightmsb
	\scriptscriptfont\msbfam\sixmsb
  \textfont\eufmfam\twelveteufm \scriptfont\eufmfam\eighteufm
	\scriptscriptfont\eufmfam\sixeufm
  \setbox\strutbox\hbox{\vrule height8.5\p@@ depth3.5\p@@ width\z@}%
  \setbox\strutbox@\hbox{\lower.5\normallineskiplimit\vbox{%
	\kern-\normallineskiplimit\copy\strutbox}}%
  \setbox\z@\vbox{\hbox{$($}\kern\z@}\b@gsize1.2\ht\z@
  \fi
  \normalbaselines\rm\dotsspace@1.5mu\ex@.2326ex\jot3\ex@}

\font@\elevenrm=cmr10 scaled 1095
\font@\elevenit=cmti10 scaled 1095
\font@\elevensl=cmsl10 scaled 1095
\font@\elevenbf=cmbx10 scaled 1095
\font@\elevensmc=cmcsc10 scaled 1095
\font@\eleventt=cmtt10 scaled 1095
\font@\eleveni=cmmi10 scaled 1095 \skewchar\eleveni='177
\font@\elevensy=cmsy10 scaled 1095 \skewchar\elevensy='60
\font@\elevenex=cmex10 scaled 1095
\font@\elevenmsa=msam10 scaled 1095
\font@\elevenmsb=msbm10 scaled 1095
\font@\eleveneufm=eufm10 scaled 1095

\def\elevenpoint{\p@@1.1\p@ \normallineskiplimit\p@@
 \mathsurround\m@ths@r \normalbaselineskip12\p@@
 \abovedisplayskip12\p@@ plus3\p@@ minus9\p@@
 \belowdisplayskip\abovedisplayskip
 \abovedisplayshortskip\z@ plus3\p@@
 \belowdisplayshortskip7\p@@ plus3\p@@ minus4\p@@
 \textonlyfont@\rm\elevenrm \textonlyfont@\it\elevenit
 \textonlyfont@\sl\elevensl \textonlyfont@\bf\elevenbf
 \textonlyfont@\smc\elevensmc \textonlyfont@\tt\eleventt
 \ifsyntax@ \def\big##1{{\hbox{$\left##1\right.$}}}%
  \let\Big\big \let\bigg\big \let\Bigg\big
 \else
  \textfont\z@\elevenrm \scriptfont\z@\eightrm \scriptscriptfont\z@\sixrm
  \textfont\@ne\eleveni \scriptfont\@ne\eighti \scriptscriptfont\@ne\sixi
  \textfont\tw@\elevensy \scriptfont\tw@\eightsy \scriptscriptfont\tw@\sixsy
  \textfont\thr@@\elevenex \scriptfont\thr@@\eightex
	\scriptscriptfont\thr@@\sevenex
  \textfont\itfam\elevenit \scriptfont\itfam\eightit
	\scriptscriptfont\itfam\sevenit
  \textfont\bffam\elevenbf \scriptfont\bffam\eightbf
	\scriptscriptfont\bffam\sixbf
  \textfont\msafam\elevenmsa \scriptfont\msafam\eightmsa
	\scriptscriptfont\msafam\sixmsa
  \textfont\msbfam\elevenmsb \scriptfont\msbfam\eightmsb
	\scriptscriptfont\msbfam\sixmsb
  \textfont\eufmfam\eleveneufm \scriptfont\eufmfam\eighteufm
	\scriptscriptfont\eufmfam\sixeufm
  \setbox\strutbox\hbox{\vrule height8.5\p@@ depth3.5\p@@ width\z@}%
  \setbox\strutbox@\hbox{\lower.5\normallineskiplimit\vbox{%
	\kern-\normallineskiplimit\copy\strutbox}}%
  \setbox\z@\vbox{\hbox{$($}\kern\z@}\b@gsize1.2\ht\z@
  \fi
  \normalbaselines\rm\dotsspace@1.5mu\ex@.2326ex\jot3\ex@}

\def\m@R@f@[#1]{\mathsurzero{
 \s@ct{}{#1}}\wr@@c{\string\Refcd{#1}{\the\pageno}}\B@gr@
 \frenchspacing\rcount\z@\refkey{\k@yf@nt[##1]}\refno{\k@yf@nt[##1]}%
 \widest{AZ}\keyright\let\Key\key\let\refin\relax}
\def\widest#1{\s@twd@\r@f@nd{\r@fk@y{\k@yf@nt#1}\enspace}}
\def\widestno#1{\s@twd@\r@f@nd{\r@fn@{\k@yf@nt#1}\enspace}}
\def\widestlabel#1{\s@twd@\r@f@nd{\k@yf@nt#1\enspace}}
\def\refkey{\def\r@fk@y##1} \def\refno{\def\r@fn@##1}
\def\keyright{\def\r@fit@m{\hang\textindent}}
\def\keyflat{\def\r@fit@m##1{\setbox\z@\hbox{##1\enspace}\hang\noindent
 \ifnum\wd\z@<\parindent\indent\hglue-\wd\z@\fi\unhbox\z@}}

\def\R@fb@x{\global\setbox\r@f@b@x} \def\K@yb@x{\global\setbox\k@yb@x}
\def\ref{\par\b@gr@\r@ff@nt\R@fb@x\box\voidb@x\K@yb@x\box\voidb@x
 \@fn@mfalse\@fl@bfalse\b@g@nr@f}
\def\c@nc@t#1{\setbox\z@\lastbox
 \setbox\adjb@x\hbox{\unhbox\adjb@x\unhbox\z@\unskip\unskip\unpenalty#1}}
\def\adjust#1{\relax\ifmmode\penalty-\@M\null\hfil$\clubpenalty\z@
 \widowpenalty\z@\interlinepenalty\z@\offinterlineskip\endgraf
 \setbox\z@\lastbox\unskip\unpenalty\c@nc@t{#1}\nt$\hfil\penalty-\@M
 \else\endgraf\c@nc@t{#1}\nt\fi}
\def\adjustnext#1{\P@nct\hbox{#1}\ignore}
\def\adjusted#1{\def\@djp@{#1}\ignore}
\def\addtoks#1{\global\@ddt@ks{#1}\ignore}
\def\addnext#1{\global\@dd@p@n{#1}\ignore}

\def\cl@s@{\adjust{\@djp@}\endgraf\setbox\z@\lastbox
 \global\setbox\@ne\hbox{\unhbox\adjb@x\ifvoid\z@\else\unhbox\z@\unskip\unskip
 \unpenalty\fi}\egroup\ifnum\c@rr@nt=\k@yb@x\global\fi
 \setbox\c@rr@nt\hbox{\unhbox\@ne\box\p@nct@}\P@nct\null
 \the\@ddt@ks\global\@ddt@ks{}}
\def\@p@n#1{\def\c@rr@nt{#1}\setbox\c@rr@nt\vbox\bgroup\let\@djp@\relax
 \hsize\maxdimen\nt\the\@dd@p@n\global\@dd@p@n{}}
\def\b@g@nr@f{\bgroup\@p@n\z@}
\def\key{\cl@s@\ifvoid\k@yb@x\@p@n\k@yb@x\k@yf@nt\else\@p@n\z@\fi}
\def\label{\cl@s@\ifvoid\k@yb@x\global\@fl@btrue\@p@n\k@yb@x\k@yf@nt\else
 \@p@n\z@\fi}
\def\no{\cl@s@\ifvoid\k@yb@x\gad\rcount\global\@fn@mtrue
 \K@yb@x\hbox{\k@yf@nt\the\rcount}\fi\@p@n\z@}
\def\labelno{\cl@s@\ifvoid\k@yb@x\gad\rcount\@fl@btrue
 \@p@n\k@yb@x\k@yf@nt\the\rcount\else\@p@n\z@\fi}
\def\by{\cl@s@\@p@n\b@b@x} \def\paper{\cl@s@\@p@n\p@p@rb@x\p@p@rf@nt\ignore}
\def\jour{\cl@s@\@p@n\j@@rb@x} \def\yr{\cl@s@\@p@n\y@@rb@x}
\def\vol{\cl@s@\@p@n\v@lb@x\v@lf@nt\ignore}
\def\issue{\cl@s@\@p@n\is@b@x\iss@f@nt\ignore}
\def\page{\cl@s@\ifp@g@s\@p@n\z@\else\p@g@true\@p@n\p@g@b@x\fi}
\def\pages{\cl@s@\ifp@g@\@p@n\z@\else\p@g@strue\@p@n\p@g@b@x\fi}
\def\inbook{\cl@s@\@p@n\inb@@kb@x}
\def\book{\cl@s@\@p@n\b@@kb@x\b@@kf@nt\ignore}
\def\publ{\cl@s@\@p@n\p@blb@x} \def\publaddr{\cl@s@\@p@n\p@bl@db@x}
\def\ed{\cl@s@\ifed@s\@p@n\z@\else\ed@true\@p@n\ed@b@x\fi}
\def\eds{\cl@s@\ifed@\@p@n\z@\else\ed@strue\@p@n\ed@b@x\fi}
\def\info{\cl@s@\@p@n\inf@b@x} \def\paperinfo{\cl@s@\@p@n\p@p@nf@b@x}
\def\bookinfo{\cl@s@\@p@n\b@@nf@b@x} 
\def\P@nct{\global\setbox\p@nct@} \def\nopunct{\P@nct\box\voidb@x}
\def\p@@@t#1#2{\ifvoid\p@nct@\else#1\unhbox\p@nct@#2\fi}
\def\sp@@{\penalty-50 \space\hskip\z@ plus.1em}
\def\c@mm@{\p@@@t,\sp@@} \def\sp@c@{\p@@@t\empty\sp@@}
\def\p@tb@x#1#2{\ifvoid#1\else#2\@nb@x#1\fi}
\def\@nb@x#1{\unhbox#1\P@nct\lastbox}
\def\endr@f@{\cl@s@\nopunct
 \R@fb@x\hbox{\unhbox\r@f@b@x \p@tb@x\b@b@x\empty
 \ifvoid\j@@rb@x\ifvoid\inb@@kb@x\ifvoid\p@p@rb@x\ifvoid\b@@kb@x
  \ifvoid\p@p@nf@b@x\ifvoid\b@@nf@b@x
  \p@tb@x\v@lb@x\c@mm@ \ifvoid\y@@rb@x\else\sp@c@(\@nb@x\y@@rb@x)\fi
  \p@tb@x\is@b@x\c@mm@ \p@tb@x\p@g@b@x\c@mm@ \p@tb@x\inf@b@x\c@mm@
  \else\p@tb@x \b@@nf@b@x\c@mm@ \p@tb@x\v@lb@x\c@mm@ \p@tb@x\is@b@x\sp@c@
  \ifvoid\ed@b@x\else\sp@c@(\@nb@x\ed@b@x,\space\ifed@ ed.\else eds.\fi)\fi
  \p@tb@x\p@blb@x\c@mm@ \p@tb@x\p@bl@db@x\c@mm@ \p@tb@x\y@@rb@x\c@mm@
  \p@tb@x\p@g@b@x{\c@mm@\ifp@g@ p\p@@nt\else pp\p@@nt\fi}%
  \p@tb@x\inf@b@x\c@mm@\fi
  \else \p@tb@x\p@p@nf@b@x\c@mm@ \p@tb@x\v@lb@x\c@mm@
  \ifvoid\y@@rb@x\else\sp@c@(\@nb@x\y@@rb@x)\fi
  \p@tb@x\is@b@x\c@mm@ \p@tb@x\p@g@b@x\c@mm@ \p@tb@x\inf@b@x\c@mm@\fi
  \else \p@tb@x\b@@kb@x\c@mm@
  \p@tb@x\b@@nf@b@x\c@mm@ \p@tb@x\p@blb@x\c@mm@
  \p@tb@x\p@bl@db@x\c@mm@ \p@tb@x\y@@rb@x\c@mm@
  \ifvoid\p@g@b@x\else\c@mm@\@nb@x\p@g@b@x p\fi \p@tb@x\inf@b@x\c@mm@ \fi
  \else \c@mm@\@nb@x\p@p@rb@x\ic@\p@tb@x\p@p@nf@b@x\c@mm@
  \p@tb@x\v@lb@x\sp@c@ \ifvoid\y@@rb@x\else\sp@c@(\@nb@x\y@@rb@x)\fi
  \p@tb@x\is@b@x\c@mm@ \p@tb@x\p@g@b@x\c@mm@\p@tb@x\inf@b@x\c@mm@\fi
  \else \p@tb@x\p@p@rb@x\c@mm@\ic@\p@tb@x\p@p@nf@b@x\c@mm@
  \c@mm@\@nb@x\inb@@kb@x \p@tb@x\b@@nf@b@x\c@mm@ \p@tb@x\v@lb@x\sp@c@
  \p@tb@x\is@b@x\sp@c@
  \ifvoid\ed@b@x\else\sp@c@(\@nb@x\ed@b@x,\space\ifed@ ed.\else eds.\fi)\fi
  \p@tb@x\p@blb@x\c@mm@ \p@tb@x\p@bl@db@x\c@mm@ \p@tb@x\y@@rb@x\c@mm@
  \p@tb@x\p@g@b@x{\c@mm@\ifp@g@ p\p@@nt\else pp\p@@nt\fi}%
  \p@tb@x\inf@b@x\c@mm@\fi
  \else\p@tb@x\p@p@rb@x\c@mm@\ic@\p@tb@x\p@p@nf@b@x\c@mm@\p@tb@x\j@@rb@x\c@mm@
  \p@tb@x\v@lb@x\sp@c@ \ifvoid\y@@rb@x\else\sp@c@(\@nb@x\y@@rb@x)\fi
  \p@tb@x\is@b@x\c@mm@ \p@tb@x\p@g@b@x\c@mm@ \p@tb@x\inf@b@x\c@mm@ \fi}}
\def\m@r@f#1#2{\endr@f@\ifvoid\p@nct@\else\R@fb@x\hbox{\unhbox\r@f@b@x
 #1\unhbox\p@nct@\penalty-200\enskip#2}\fi\egroup\b@g@nr@f}
\def\endref{\endr@f@\ifvoid\p@nct@\else\R@fb@x\hbox{\unhbox\r@f@b@x.}\fi
 \parindent\r@f@nd
 \r@fit@m{\ifvoid\k@yb@x\else\if@fn@m\r@fn@{\unhbox\k@yb@x}\else
 \if@fl@b\unhbox\k@yb@x\else\r@fk@y{\unhbox\k@yb@x}\fi\fi\fi}\unhbox\r@f@b@x
 \endgraf\egroup\endgroup}
\def\moreref{\m@r@f;\empty}
\def\transl{\m@r@f;{\unskip\space
 {\sl English translation\ic@}:\penalty-66 \space}}
\def\endRefs{\endgraf\goodbreak\endgroup}

\hyphenation{acad-e-my acad-e-mies af-ter-thought anom-aly anom-alies
an-ti-deriv-a-tive an-tin-o-my an-tin-o-mies apoth-e-o-ses
apoth-e-o-sis ap-pen-dix ar-che-typ-al as-sign-a-ble as-sist-ant-ship
as-ymp-tot-ic asyn-chro-nous at-trib-uted at-trib-ut-able bank-rupt
bank-rupt-cy bi-dif-fer-en-tial blue-print busier busiest
cat-a-stroph-ic cat-a-stroph-i-cally con-gress cross-hatched data-base
de-fin-i-tive de-riv-a-tive dis-trib-ute dri-ver dri-vers eco-nom-ics
econ-o-mist elit-ist equi-vari-ant ex-quis-ite ex-tra-or-di-nary
flow-chart for-mi-da-ble forth-right friv-o-lous ge-o-des-ic
ge-o-det-ic geo-met-ric griev-ance griev-ous griev-ous-ly
hexa-dec-i-mal ho-lo-no-my ho-mo-thetic ideals idio-syn-crasy
in-fin-ite-ly in-fin-i-tes-i-mal ir-rev-o-ca-ble key-stroke
lam-en-ta-ble light-weight mal-a-prop-ism man-u-script mar-gin-al
meta-bol-ic me-tab-o-lism meta-lan-guage me-trop-o-lis
met-ro-pol-i-tan mi-nut-est mol-e-cule mono-chrome mono-pole
mo-nop-oly mono-spline mo-not-o-nous mul-ti-fac-eted mul-ti-plic-able
non-euclid-ean non-iso-mor-phic non-smooth par-a-digm par-a-bol-ic
pa-rab-o-loid pa-ram-e-trize para-mount pen-ta-gon phe-nom-e-non
post-script pre-am-ble pro-ce-dur-al pro-hib-i-tive pro-hib-i-tive-ly
pseu-do-dif-fer-en-tial pseu-do-fi-nite pseu-do-nym qua-drat-ic
quad-ra-ture qua-si-smooth qua-si-sta-tion-ary qua-si-tri-an-gu-lar
quin-tes-sence quin-tes-sen-tial re-arrange-ment rec-tan-gle
ret-ri-bu-tion retro-fit retro-fit-ted right-eous right-eous-ness
ro-bot ro-bot-ics sched-ul-ing se-mes-ter semi-def-i-nite
semi-ho-mo-thet-ic set-up se-vere-ly side-step sov-er-eign spe-cious
spher-oid spher-oid-al star-tling star-tling-ly sta-tis-tics
sto-chas-tic straight-est strange-ness strat-a-gem strong-hold
sum-ma-ble symp-to-matic syn-chro-nous topo-graph-i-cal tra-vers-a-ble
tra-ver-sal tra-ver-sals treach-ery turn-around un-at-tached
un-err-ing-ly white-space wide-spread wing-spread wretch-ed
wretch-ed-ly Brown-ian Eng-lish Euler-ian Feb-ru-ary Gauss-ian
Grothen-dieck Hamil-ton-ian Her-mit-ian Jan-u-ary Japan-ese Kor-te-weg
Le-gendre Lip-schitz Lip-schitz-ian Mar-kov-ian Noe-ther-ian
No-vem-ber Rie-mann-ian Schwarz-schild Sep-tem-ber}

\let\nopagenumber\p@gen@false \let\putpagenumber\p@gen@true

\outer\def\myRefs{\myR@fs} \r@st@re\proclaim
\def\bye{\par\vfill\supereject\cl@selbl\cl@secd\b@e} \r@endd@\b@e
 \let\Key\key \def\endpro{\par\endproclaim}
\let\d@c@\document \def\document{\d@c@\tenpoint}

\newtoks\@@tp@t \@@tp@t\output
\output=\@ft@{\let\{\noexpand\the\@@tp@t}
\let\{\relax

\newif\ifVersion \Versiontrue
\def\p@n@l#1{\ifnum#1=\z@\else\penalty#1\relax\fi}

\def\s@ct#1#2{\ifVersion
 \skip@\lastskip\ifdim\skip@<1.5\bls\vskip-\skip@\p@n@l{-200}\vsk.5>%
 \p@n@l{-200}\vsk.5>\p@n@l{-200}\vsk.5>\p@n@l{-200}\vsk-1.5>\else
 \p@n@l{-200}\fi\ifdim\skip@<.9\bls\vsk.9>\else
 \ifdim\skip@<1.5\bls\vskip\skip@\fi\fi
 \vtop{\twelvepoint\raggedright\s@cf@nt\vp1\vsk->\vskip.16ex
 \s@twd@\parindent{#1}%
 \ifdim\parindent>\z@\adv\parindent.5em\fi\hang\textindent{#1}#2\strut}
 \else
 \p@sk@p{-200}{.8\bls}\vtop{\s@cf@nt\s@twd@\parindent{#1}%
 \ifdim\parindent>\z@\adv\parindent.5em\fi\hang\textindent{#1}#2\strut}\fi
 \nointerlineskip\nobreak\vtop{\strut}\nobreak\vskip-.6\bls\nobreak}

\def\s@bs@ct#1#2{\ifVersion
 \skip@\lastskip\ifdim\skip@<1.5\bls\vskip-\skip@\p@n@l{-200}\vsk.5>%
 \p@n@l{-200}\vsk.5>\p@n@l{-200}\vsk.5>\p@n@l{-200}\vsk-1.5>\else
 \p@n@l{-200}\fi\ifdim\skip@<.9\bls\vsk.9>\else
 \ifdim\skip@<1.5\bls\vskip\skip@\fi\fi
 \vtop{\elevenpoint\raggedright\s@bf@nt\vp1\vsk->\vskip.16ex%
 \s@twd@\parindent{#1}\ifdim\parindent>\z@\adv\parindent.5em\fi
 \hang\textindent{#1}#2\strut}
 \else
 \p@sk@p{-200}{.6\bls}\vtop{\s@bf@nt\s@twd@\parindent{#1}%
 \ifdim\parindent>\z@\adv\parindent.5em\fi\hang\textindent{#1}#2\strut}\fi
 \nointerlineskip\nobreak\vtop{\strut}\nobreak\vskip-.8\bls\nobreak}

\def\gadv{\global\adv} \def\gad#1{\gadv#1\@ne} \def\gadneg#1{\gadv#1-\@ne}

\newcount\t@@n \t@@n=10 \newbox\testbox

\newcount\Sno \newcount\Lno \newcount\Fno

\def\pr@cl#1{\r@st@re\pr@c@\pr@c@{#1}\global\let\pr@c@\relax}

\def\l@L#1{\l@bel{#1}L} \def\l@F#1{\l@bel{#1}F} \def\<#1>{\l@b@l{#1}F}

\def\tagg#1{\tag"\rlap{\rm(#1)}\kern.01\p@"}
\def\Tag#1{\tag{\l@F{#1}}} \def\Tagg#1{\tagg{\l@F{#1}}}

\def\Th#1{\pr@cl{Theorem\xspace\l@L{#1}}\ignore}
\def\Lm#1{\pr@cl{Lemma\xspace\l@L{#1}}\ignore}
\def\Cr#1{\pr@cl{Corollary\xspace\l@L{#1}}\ignore}
\def\Df#1{\pr@cl{Definition\xspace\l@L{#1}}\ignore}
\def\Cj#1{\pr@cl{Conjecture\xspace\l@L{#1}}\ignore}
\def\Prop#1{\pr@cl{Proposition\xspace\l@L{#1}}\ignore}

\def\Pf#1.{\demo{Proof #1}\bgroup\ignore}
\def\edemo{\par\egroup\enddemo} \def\epf{\qed\edemo}
\def\Rem{\demo{\sl Remark}} \def\Ex{\demo{\bf Example}}

\def\Proof#1.{\demo{\let\{\relax Proof #1}\def\t@st@{#1}%
 \ifx\t@st@\empty\else\xdef\@@wr##1##2##3##4{##1{##2##3{\the\cdn@}{##4}}}%
 \wr@@c@{\the\cdn@}{Proof #1}\@@wr\wr@@c\string\subcd{\the\pageno}\fi
 \bgroup\ignore}

\def\Ap@x{Appendix}
\def\Appendix{\Sno=64 \t@@n\@ne \wr@@c{\string\Appencd}
 \def\sf@rm{\char\the\Sno} \def\sf@rm@{\Ap@x\space\sf@rm} \def\sf@rm@@{\Ap@x}
 \def\s@ct@n##1##2{\s@ct\empty{\setbox\z@\hbox{##1}\ifdim\wd\z@=\z@
 \if##2*\sf@rm@@\else\if##2.\sf@rm@@.\else##2\fi\fi\else
 \if##2*\sf@rm@\else\if##2.\sf@rm@.\else\sf@rm@.\enspace##2\fi\fi\fi}}}
\def\Appcd#1#2#3{\gad\Cdentry\global\cdentry\z@\def\Ap@@{\hglue-\l@ftcd\Ap@x}
 \ifx\@ppl@ne\empty\def\l@@b{\@fwd@@{#1}{\space#1}{}}
 \if*#2\entcd{}{\Ap@@\l@@b}{#3}\else\if.#2\entcd{}{\Ap@@\l@@b.}{#3}\else
 \entcd{}{\Ap@@\l@@b.\enspace#2}{#3}\fi\fi\else
 \def\l@@b{\@fwd@@{#1}{\c@l@b{#1}}{}}\if*#2\entcd{\l@@b}{\Ap@x}{#3}\else
 \if.#2\entcd{\l@@b}{\Ap@x.}{#3}\else\entcd{\l@@b}{#2}{#3}\fi\fi\fi}

\let\s@ct@n\s@ct
\def\s@ct@@[#1]#2{\@ft@\xdef\csname @#1@S@\endcsname{\sf@rm}\wr@@x{}%
 \wr@@x{\string\labeldef{S}\space{\?#1@S?}\space{#1}}%
 {
 \s@ct@n{\sf@rm@}{#2}}\wr@@c{\string\Entcd{\?#1@S?}{#2}{\the\pageno}}}
\def\s@ct@#1{\wr@@x{}{
 \s@ct@n{\sf@rm@}{#1}}\wr@@c{\string\Entcd{\sf@rm}{#1}{\the\pageno}}}
\def\s@ct@e[#1]#2{\@ft@\xdef\csname @#1@S@\endcsname{\sf@rm}\wr@@x{}%
 \wr@@x{\string\labeldef{S}\space{\?#1@S?}\space{#1}}%
 {
 \s@ct@n\empty{#2}}\wr@@c{\string\Entcd{}{#2}{\the\pageno}}}
\def\s@cte#1{\wr@@x{}{
 \s@ct@n\empty{#1}}\wr@@c{\string\Entcd{}{#1}{\the\pageno}}}
\def\theSno#1#2{\dff\?#1@S?{#2}%
 \wr@@x{\string\labeldef{S}\space{#2}\space{#1}}\fi}

\newif\ifd@bn@\d@bn@true
\def\Section{\gad\Sno\ifd@bn@\Fno\z@\Lno\z@\fi\@fn@xt[\s@ct@@\s@ct@}
\def\section{\gad\Sno\ifd@bn@\Fno\z@\Lno\z@\fi\@fn@xt[\s@ct@e\s@cte}
\let\Sect\Section 
\def\subsection{\@fn@xt*\subs@ct@\subs@ct}
\def\subs@ct#1{{\s@bs@ct\empty{#1}}\wr@@c{\string\subcd{#1}{\the\pageno}}}
\def\subs@ct@*#1{\vsk->\nobreak
 {\s@bs@ct\empty{#1}}\wr@@c{\string\subcd{#1}{\the\pageno}}}

\def\l@b@l#1#2{\def\n@@{\csname #2no\endcsname}%
 \if*#1\gad\n@@ \@ft@\xdef\csname @#1@#2@\endcsname{\l@f@rm}\else\def\t@st{#1}%
 \ifx\t@st\empty\gad\n@@ \@ft@\xdef\csname @#1@#2@\endcsname{\l@f@rm}%
 \else\@ft@\ifx\csname @#1@#2@mark\endcsname\relax\gad\n@@
 \@ft@\xdef\csname @#1@#2@\endcsname{\l@f@rm}%
 \@ft@\gdef\csname @#1@#2@mark\endcsname{}%
 \wr@@x{\string\labeldef{#2}\space{\?#1@#2?}\space\ifnum\n@@<10 \space\fi{#1}}%
 \fi\fi\fi}
\def\labeldef#1#2#3{\dff\?#3@#1?{#2}}
\def\Labeldef#1#2#3{\dff\?#3@#1?{#2}\@ft@\gdef\csname @#3@#1@mark\endcsname{}}

\def\l@bel#1#2{\l@b@l{#1}{#2}\?#1@#2?}

\newcount\c@cite
\def\?#1?{\csname @#1@\endcsname}
\def\[{\@fn@xt:\c@t@sect\c@t@}
\def\c@t@#1]{{\c@cite\z@\@fwd@@{\?#1@L?}{\adv\c@cite1}{}%
 \@fwd@@{\?#1@F?}{\adv\c@cite1}{}\@fwd@@{\?#1?}{\adv\c@cite1}{}%
 \relax\ifnum\c@cite=\z@{\bf ???}\wrs@x{No label [#1]}\else
 \ifnum\c@cite=1\let\@@PS\relax\let\@@@\relax\else\let\@@PS\underbar
 \def\@@@{{\rm<}}\fi\@@PS{\?#1?\@@@\?#1@L?\@@@\?#1@F?}\fi}}
\def\(#1){{\rm(\c@t@#1])}}
\def\c@t@s@ct#1{\@fwd@@{\?#1@S?}{\?#1@S?\relax}%
 {{\bf ???}\wrs@x{No section label {#1}}}}
\def\c@t@sect:#1]{\c@t@s@ct{#1}} \let\SNo\c@t@s@ct

\newdimen\l@ftcd \newdimen\r@ghtcd \let\nlc\relax
\newcount\Cdentry \newcount\cdentry \let\prentcd\relax \let\postentcd\relax

\def\d@tt@d{\leaders\hbox to 1em{\kern.1em.\hfil}\hfill}
\def\entcd#1#2#3{\gad\cdentry\prentcd\item{\l@bcdf@nt#1}{\entcdf@nt#2}\alb
 \kern.9em\hbox{}\kern-.9em\d@tt@d\kern-.36em{\p@g@cdf@nt#3}\kern-\r@ghtcd
 \hbox{}\postentcd\par}
\def\Entcd#1#2#3{\gad\Cdentry\global\cdentry\z@
 \def\l@@b{\@fwd@@{#1}{\c@l@b{#1}}{}}\vsk.2>\entcd{\l@@b}{#2}{#3}}
\def\subcd#1#2{{\adv\leftskip.333em\entcd{}{\s@bcdf@nt#1}{#2}}}
\def\Refcd#1#2{\def\t@@st{#1}\ifx\t@@st\empty\ifx\r@fl@ne\empty\relax\else
 \R@fcd{\r@fl@ne}{#2}\fi\else\R@fcd{#1}{#2}\fi}
\def\R@fcd#1#2{\sk@@p{.6\bls}\entcd{}{\hglue-\l@ftcd\R@fcdf@nt #1}{#2}}
\def\Refline{\def\r@fl@ne} \def\Refempty{\let\r@fl@ne\empty}
\def\Appencd{\par\adv\leftskip-\l@ftcd\adv\rightskip-\r@ghtcd\@ppl@ne
 \adv\leftskip\l@ftcd\adv\rightskip\r@ghtcd\let\Entcd\Appcd}
\def\appline{\def\@ppl@ne} \def\Appempty{\let\@ppl@ne\empty}
\def\Appline#1{\def\@ppl@ne{\s@bs@ct{}{#1}}}
\def\Leftcd#1{\adv\leftskip-\l@ftcd\s@twd@\l@ftcd{\c@l@b{#1}\enspace}
 \adv\leftskip\l@ftcd}
\def\Rightcd#1{\adv\rightskip-\r@ghtcd\s@twd@\r@ghtcd{#1\enspace}
 \adv\rightskip\r@ghtcd}
\def\C@nt{Contents} \def\Ap@s{Appendices} \def\R@fcs{References}
\def\contents{\@fn@xt*\cont@@\cont@}
\def\cont@{\@fn@xt[\cnt@{\cnt@[\C@nt]}}
\def\cont@@*{\@fn@xt[\cnt@@{\cnt@@[\C@nt]}}
\def\cnt@[#1]{\c@nt@{M}{#1}{44}{\s@bs@ct{}{\@ppl@f@nt\Ap@s}}}
\def\cnt@@[#1]{\c@nt@{M}{#1}{44}{}}
\def\endco{\par\penalty-500\vsk>\vskip\z@\endgroup}
\def\readcd{\@np@t{\jobname.cd}}
\def\Cde{\@fn@xt*\Cde@@\Cde@}
\def\Cde@{\@fn@xt[\Cd@{\Cd@[\C@nt]}}
\def\Cde@@*{\@fn@xt[\Cd@@{\Cd@@[\C@nt]}}
\def\Cd@[#1]{\cnt@[#1]\readcd\endco}
\def\Cd@@[#1]{\cnt@@[#1]\readcd\endco}
\def\contlabeldef{\def\c@l@b}

\long\def\c@nt@#1#2#3#4{\s@twd@\l@ftcd{\c@l@b{#1}\enspace}
 \s@twd@\r@ghtcd{#3\enspace}\adv\r@ghtcd1.333em
 \def\@ppl@ne{#4}\def\r@fl@ne{\R@fcs}\s@ct{}{#2}\B@gr@\parindent\z@\let\nlc\nl
 \let\nl\relax\parskip.2\bls\adv\leftskip\l@ftcd\adv\rightskip\r@ghtcd}

\def\writecd{\immediate\openout\@@cd\jobname.cd \def\wr@@c{\write\@@cd}
 \def\cl@secd{\immediate\write\@@cd{\string\endinput}\immediate\closeout\@@cd}
 \def\closecd{\cl@secd\global\let\cl@secd\relax}}
\let\cl@secd\relax \def\wr@@c#1{} \let\closecd\relax

\def\dff{\@ft@\d@f} \def\d@f{\@ft@\def}
\def\edff{\@ft@\ed@f} \def\ed@f{\@ft@\edef}
\def\gdff{\@ft@\gd@f} \def\gd@f{\@ft@\gdef}
\def\defi#1#2{\def#1{#2}\wr@@x{\string\def\string#1{#2}}}

\def\qed{\hbox{}\nobreak\hfill\nobreak{\m@th$\,\square$}}
\def\back#1 {\strut\kern-.33em #1\enspace\ignore} 
\def\Text#1{\crcr\noalign{\alb\vsk>\normalbaselines\vsk->\vbox{\nt #1\strut}%
 \nobreak\nointerlineskip\vbox{\strut}\nobreak\vsk->\nobreak}}

\def\hcor#1{\advance\hoffset by #1}
\def\vcor#1{\advance\voffset by #1}
\let\bls\baselineskip \let\ignore\ignorespaces
\ifx\ic@\undefined \let\ic@\/\fi
\def\vsk#1>{\vskip#1\bls} \let\adv\advance
\def\vv#1>{\vadjust{\vsk#1>}\ignore}
\def\vvn#1>{\vadjust{\nobreak\vsk#1>\nobreak}\ignore}
\def\vvv#1>{\vskip\z@\vsk#1>\nt\ignore}
\def\vvgood{\vadjust{\penalty-500}}

\def\nngood{\noalign{\penalty-500}}
\def\Goodbreak{\par\penalty-\@m}
\def\vvGood{\vadjust{\penalty-\@m}}

\def\nnGood{\noalign{\penalty-\@m}}
\def\wwgood#1:#2>{\vv#1>\vvgood\vv#2>\vv0>}
\def\vskgood#1:#2>{\vsk#1>\goodbreak\vsk#2>\vsk0>}
\def\cnngood#1:#2>{\cnn#1>\nngood\cnn#2>}
\def\cnnGood#1:#2>{\cnn#1>\nnGood\cnn#2>}
\def\ragood{\vadjust{\vskip\z@ plus 12pt}\vvgood}

\def\hfilll{\hskip\z@ plus 1filll} \def\vfilll{\vskip\z@ plus 1filll}

 \def\setparindent{\edef\Parindent{\the\parindent}}
\def\Type{\vsk.5>\bgroup\parindent\z@\tt\rightskip\z@ plus1em minus1em%
 \spaceskip.3333em \xspaceskip.5em\relax}
\def\endType{\vsk.5>\egroup\nt} 

\let\Hat\widehat \let\Tilde\widetilde \let\dollar\$ \let\ampersand\&
\let\sss\scriptscriptstyle  
\let\vp\vphantom \let\hp\hphantom \let\nt\noindent
\let\cline\centerline \let\lline\leftline \let\rline\rightline
\def\nn#1>{\noalign{\vskip#1\p@@}} \def\NN#1>{\openup#1\p@@}
\def\cnn#1>{\noalign{\vsk#1>}}
 
\let\Lim\lim \def\lim{\Lim\limits} \let\Sum\sum \def\sum{\Sum\limits}
\def\Plus{\bigoplus\limits} 
\let\Prod\prod \def\prod{\Prod\limits} \let\Int\int \def\int{\Int\limits}

\def\tsum{\mathop{\tsize\Sum}\limits} 
 \def\&{.\kern.1em}
\def\nl{\leavevmode\hfill\break} \def\~{\leavevmode\@fn@xt~\m@n@s\@md@@sh}
\def\@md@@sh{\@fn@xt-\d@@sh\@md@sh} \def\@md@sh{\raise.13ex\hbox{--}}
\def\m@n@s~{\raise.15ex\mbox{-}} \def\d@@sh-{\raise.15ex\hbox{-}}
\def\dash{\raise.15ex\hbox{-}} 

\let\percent\% \def\%#1{\ifmmode\mathop{#1}\limits\else\percent#1\fi}
\let\@ml@t\" \def\"#1{\ifmmode ^{(#1)}\else\@ml@t#1\fi}
\let\@c@t@\' \def\'#1{\ifmmode _{(#1)}\else\@c@t@#1\fi}
\let\colon\: \def\:{^{\vp{\topsmash|}}} 

\def\xspace{\kern.34em}

\def\>{\RIfM@\mskip.666667\thinmuskip\relax\else\kern.111111em\fi}
\def\}{\RIfM@\mskip-.666667\thinmuskip\relax\else\kern-.111111em\fi}
\def\){\RIfM@\mskip.333333\thinmuskip\relax\else\kern.0555556em\fi}
\def\]{\RIfM@\mskip-.333333\thinmuskip\relax\else\kern-.0555556em\fi}
\def\d@t{.\hskip.3em} \def\.{\d@t\ignore}

\let\texspace\ \def\ {\ifmmode\alb\fi\texspace}

\let\n@wp@ge\newpage \def\newpage{\endgraf\n@wp@ge}
\let\=\m@th \def\mbox#1{\hbox{\m@th$#1$}}
\def\mtext#1{\text{\m@th$#1$}} \def\^#1{\text{\m@th#1}}
\def\Line#1{\kern-.5\hsize\line{\m@th$\dsize#1$}\kern-.5\hsize}
\def\Lline#1{\kern-.5\hsize\lline{\m@th$\dsize#1$}\kern-.5\hsize}
\def\Cline#1{\kern-.5\hsize\cline{\m@th$\dsize#1$}\kern-.5\hsize}
\def\Rline#1{\kern-.5\hsize\rline{\m@th$\dsize#1$}\kern-.5\hsize}

\def\Ll@p#1{\llap{\m@th$#1$}} \def\Rl@p#1{\rlap{\m@th$#1$}}
 \def\Cl@p#1{\llap{\m@th$#1$\hss}}
\def\Llap#1{\mathchoice{\Ll@p{\dsize#1}}{\Ll@p{\tsize#1}}{\Ll@p{\ssize#1}}%
 {\Ll@p{\sss#1}}}
\def\Clap#1{\mathchoice{\Cl@p{\dsize#1}}{\Cl@p{\tsize#1}}{\Cl@p{\ssize#1}}%
 {\Cl@p{\sss#1}}}
\def\Rlap#1{\mathchoice{\Rl@p{\dsize#1}}{\Rl@p{\tsize#1}}{\Rl@p{\ssize#1}}%
 {\Rl@p{\sss#1}}}
 
\def\LRtph#1#2{\setbox\z@\hbox{#1}\dimen\z@\wd\z@\hbox{\hbox to\dimen\z@{#2}}}
\def\LRph#1#2{\LRtph{\m@th$#1$}{\m@th$#2$}}

\def\Lto#1{\setbox\z@\mbox{\tsize{#1}}%
 \mathrel{\mathop{\hbox to\wd\z@{\rightarrowfill}}\limits#1}}
\def\Lgets#1{\setbox\z@\mbox{\tsize{#1}}%
 \mathrel{\mathop{\hbox to\wd\z@{\leftarrowfill}}\limits#1}}
\def\vpb#1{{\vp{\big(}}^{\]#1}} \def\vpp#1{{\vp{\big]}}_{#1}}
\def\lbc{\mathopen{[\![}} \def\rbc{\mathclose{]\!]}}

\def\lsb{\mathopen{\hbox{\ninepoint[}}}
\def\rsb{\mathclose{\hbox{\ninepoint]}}}
\def\bigll#1{\mathopen{\lower.12ex\hbox{\twelvepoint$\big#1\n@space$}}}
\def\bigrr#1{\mathclose{\lower.12ex\hbox{\twelvepoint$\big#1\n@space$}}}

\let\alb\allowbreak 
\def\ald{\noalign{\alb}} \let\alds\allowdisplaybreaks

\let\o\circ \let\x\times \let\ox\otimes 
\let\sub\subset  \let\tabs\+
\let\le\leqslant \let\ge\geqslant
\let\der\partial \let\8\infty \let\*\star
\let\bra\langle \let\ket\rangle
  \let\map\mapsto 
\let\hto\hookrightarrow

 \def\vert{\ |\ }

\let\lb\lbrace \let\rb\rbrace

\edef\h@xms{\hexnumber@\msafam}
\mathchardef\rarh"0\h@xms4B
\mathchardef\larh"0\h@xms4C

\def\lsym#1{#1\alb\ldots\relax#1\alb}
\def\lc{\lsym,}  \def\lx{\lsym\x} \def\lox{\lsym\ox}
\def\llc{\,,\alb\ {\ldots\ ,}\alb\ }
 
\def\End{\mathop{\roman{End}\>}\nolimits}
\def\Hom{\mathop{\roman{Hom}\>}\nolimits}
\def\Sym{\mathop{\roman{Sym}\)}\nolimits}

\def\qdet{\mathop{\roman{q\)det}}\nolimits}
\def\tr{\mathop{\roman{tr}}\nolimits}

\def\res{\mathop{\roman{res}\>}\limits}
\def\sign{\mathop{\roman{sign}\)}\limits}

\def\id{\roman{id}}  
 \def\1{^{-1}} \let\underscore\_ \def\_#1{_{\Rlap{#1}}}
\def\vst#1{{\lower1.9\p@@\mbox{\bigr|_{\raise.5\p@@\mbox{\ssize#1}}}}}
\def\vrp#1:#2>{{\vrule height#1 depth#2 width\z@}}
\def\vru#1>{\vrp#1:\z@>} \def\vrd#1>{\vrp\z@:#1>}
\def\qqq{\qquad\quad} 
\def\sscr#1{\raise.3ex\mbox{\sss#1}} \def\@@PS{\bold{OOPS!!!}}

\def\intcl{\mathop
 {\Rlap{\raise.3ex\mbox{\kern.12em\curvearrowleft}}\int}\limits}
\def\intcr{\mathop
 {\Rlap{\raise.3ex\mbox{\kern.24em\curvearrowright}}\int}\limits}

\def\pms{\raise.25ex\mbox{\ssize\pm}\>}
\def\mps{\raise.25ex\mbox{\ssize\mp}\>}
 \def\mss{{\sscr-}}

\let\dl\delta \let\Dl\Delta 
\let\epe\epsilon \let\eps\varepsilon \let\epsilon\eps

\let\zt\zeta
\let\tht\theta

 \let\La\Lambda

\let\si\sigma 
 
\let\pho\phi \let\phi\varphi

\def\C{\Bbb C}
\def\R{\Bbb R}
\def\Z{\Bbb Z}

\def\AA{\Bbb A}
\def\BB{\Bbb B}

\def\FF{\Bbb F}

\newif\ifh@ph@
\def\h@phh{\ifh@ph@\)-\else\space\fi}

\let\hyphenhands\h@ph@true
\let\nohyphenhands\h@ph@false

\def\lhs/{the left\h@phh hand side} \def\rhs/{the right\h@phh hand side}
\def\Lhs/{The left\h@phh hand side} \def\Rhs/{The right\h@phh hand side}

\def\difl/{differential} \def\dif/{difference}
\def\cf.{cf.\ \ignore} \def\Cf.{Cf.\ \ignore}
\def\egv/{eigenvector} \def\eva/{eigenvalue} \def\eq/{equation}
\def\gby/{generated by} \def\wrt/{with respect to} \def\st/{such that}
\def\resp/{respectively} \def\sol/{solution} \def\wt/{weight}
\def\offdiag/{off\)-\)diag\-o\-nal} \def\off/{offdiagonal}
\def\pol/{polynomial}
\def\rat/{rational} \def\trig/{trigonometric} \let\tri\trig
\def\fn/{function} \def\var/{variable} \def\raf/{\rat/ \fn/}
\def\inv/{invariant} \def\hol/{holomorphic} \def\hof/{\hol/ \fn/}
\def\mer/{meromorphic} \def\mef/{\mer/ \fn/} \def\mult/{multiplicity}
\def\sym/{symmetric} \def\symg/{\sym/ group}
\def\perm/{permutation} \def\part/{partition}
\def\rep/{representation} \def\irr/{irreducible} \def\irrep/{\irr/ \rep/}
\def\hom/{homomorphism} \def\aut/{automorphism} \def\iso/{isomorphism}
\def\lex/{lexicographical} \def\as/{asymptotic} \def\asex/{\as/ expansion}
\def\ndeg/{nondegenerate} \def\neib/{neighbourhood} \def\deq/{\dif/ \eq/}
\def\hw/{highest \wt/} \def\gv/{generating vector} \def\eqv/{equivalent}
\def\msd/{method of steepest descend} \def\pd/{pairwise distinct}
\def\wlg/{without loss of generality} \def\Wlg/{Without loss of generality}
\def\onedim/{one\)-dim\-en\-sion\-al} \def\fdim/{fin\-ite\)-dim\-en\-sion\-al}
\def\qcl/{quasiclassical} \def\hwv/{\hw/ vector} \let\fd\fdim
\def\hgeom/{hyper\-geo\-met\-ric} \def\hint/{\hgeom/ integral}
\def\hwm/{\hw/ module} \def\Vmod/{Verma module} \def\emod/{evaluation module}
\def\anf/{analytic \fn/} \def\anco/{analytic continuation}
\def\qg/{quantum group} \def\eqg/{elliptic \qg/} \def\dqg/{dynamical \qg/}
\def\qaff/{quantum affine algebra} \def\qloop/{quantum loop algebra}
 \def\eval/{evaluation} \def\ehom/{\eval/ \hom/} \def\epoint/{\eval/ point}
\def\tmx/{transfer matrix} \def\tmcs/{transfer matrices}

\hyphenation{coef-fi-cient orthog-o-nal}

\def\h@ph{\discretionary{}{}{-}}
\def\$#1$-{\kern\mathsurround\text{\m@th$#1$}\h@ph}
\def\prfx#1#2{\kern\mathsurround\text{\m@th$#1$#2}}

\def\Rm/{\prfx{R\)}-mat\-rix} \def\Rms/{\prfx{R\)}-mat\-ri\-ces}
\def\Ba/{Bethe ansatz} \def\Bv/{Bethe vector} \def\Bae/{\Ba/ \eq/}
\def\KZv/{Knizh\-nik\)-Zamo\-lod\-chi\-kov} \def\KZvB/{\KZv/\)-Bernard}
\def\KZ/{{\sl KZ\/}} \def\qKZ/{{\sl qKZ\/}}
\def\KZB/{{\sl KZB\/}} \def\qKZB/{{\sl qKZB\/}}
\def\qKZo/{\qKZ/ operator} \def\qKZc/{\qKZ/ connection}
\def\KZe/{\KZ/ \eq/} \def\qKZe/{\qKZ/ \eq/} \def\qKZBe/{\qKZB/ \eq/}
\def\XXX/{{\sl XXX\/}} \def\XXZ/{{\sl XXZ\/}} \def\XYZ/{{\sl XYZ\/}}
\def\YB/{Yang\)-Baxter \eq/} \def\GZ/{Gel\-fand\)-Zet\-lin}
\def\PBW/{Poincar\'e\>-Birkhoff\>-Witt}
\def\PBWi/{Poincar\'e\)-Birkhoff\>-\]Witt}

\def\SPb/{St\&Peters\-burg}
\def\home/{\SPb/ Branch of Steklov Mathematical Institute}
\def\homeaddr/{Fontanka 27, \SPb/ \,191023, Russia}
\def\homemail/{vt\@ pdmi.ras.ru}
\def\absence/{On leave of absence from \home/}
\def\support/{Supported in part by}
\def\imail/{vt\@ math.iupui.edu}
\def\DMS/{Department of Mathematical Sciences}
\def\IUPUI/{Indiana University Purdue University at Indianapolis}
\def\IUPUIaddr/{Indianapolis, IN 46202, USA}
\def\UNC/{Department of Mathematics, University of North Carolina}
\def\ChH/{Chapel Hill} \def\UNCaddr/{\ChH/, NC 27599, USA}
\def\avemail/{anv\@ email.unc.edu}

\def\Aomoto/{K\&Aomoto}
\def\Cher/{I\&Che\-red\-nik}
\def\Dri/{V\]\&G\&Drin\-feld}
\def\Fadd/{L\&D\&Fad\-deev}
\def\Feld/{G\&Felder}
\def\Fre/{I\&B\&Fren\-kel}
\def\Etingof/{P\]\&Etin\-gof}
\def\Gustaf/{R\&A\&Gustafson}
\def\Izergin/{A\&G\&Izer\-gin}
\def\Jimbo/{M\&Jimbo}
\def\Kazh/{D\&Kazh\-dan}
\def\Khor/{S\&Khorosh\-kin}
\def\Kor/{V\]\&E\&Kore\-pin}
\def\Kulish/{P\]\&P\]\&Ku\-lish}
\def\Lusz/{G\&Lusztig}
\def\Miwa/{T\]\&Miwa}
\def\MN/{M\&Naza\-rov}
\def\Mukhin/{E\&Mukhin}
\def\Reshet/{N\&Reshe\-ti\-khin} \def\Reshy/{N\&\]Yu\&Reshe\-ti\-khin}
\def\Rimanyi/{R\&Rim\'anyi}
\def\SchV/{V\]\&\]V\]\&Schecht\-man} \def\Sch/{V\]\&Schecht\-man}
\def\Skl/{E\&K\&Sklya\-nin}
\def\Smirnov/{F\]\&A\&Smir\-nov}
\def\Takh/{L\&A\&Takh\-tajan}
\def\VT/{V\]\&Ta\-ra\-sov} \def\VoT/{V\]\&O\&Ta\-ra\-sov}
\def\Varch/{A\&\]Var\-chenko} \def\Varn/{A\&N\&\]Var\-chenko}
\def\Zhel/{D\&P\]\&Zhe\-lo\-benko}

\def\AiA/{Al\-geb\-ra i Ana\-liz}
\def\DAN/{Do\-kla\-dy AN SSSR}
\def\FAA/{Funk\.Ana\-liz i ego pril.}
\def\Izv/{Iz\-ves\-tiya AN SSSR, ser\.Ma\-tem.}
\def\TMF/{Teo\-ret\.Ma\-tem\.Fi\-zi\-ka}
\def\UMN/{Uspehi Matem.\ Nauk}

\def\Adv/{Adv\.Math.}
\def\AMS/{Amer\.Math\.Society}
\def\AMSa/{AMS \publaddr Providence RI}
\def\AMST/{\AMS/ Transl.,\ Ser\&\)2}
\def\AMSTr/{\AMS/ Transl.,} \def\Ser2{Ser\&\)2}
\def\Astq/{Ast\'erisque}
\def\ContM/{Contemp\.Math.}
\def\CMP/{Comm\.Math\.Phys.}
\def\DMJ/{Duke\.Math\.J.}
\def\IJM/{Int\.J\.Math.}
\def\IMRN/{Int\.Math\.Res.\ Notices}
\def\Inv/{Invent\.Math.} 
\def\JMP/{J\.Math\.Phys.}
\def\JPA/{J\.Phys.\ A}
\def\JSM/{J\.Soviet Math.}
\def\JSP/{J\.Stat\.Phys.}
\def\LMJ/{Leningrad Math.\ J.}
\def\LpMJ/{\SPb/ Math.\ J.}
\def\LMP/{Lett\.Math\.Phys.}
\def\NMJ/{Nagoya Math\.J.}
\def\Nucl/{Nucl\.Phys.\ B}
\def\OJM/{Osaka J\.Math.}
\def\RIMS/{Publ\.RIMS, Kyoto Univ.}
\def\SIAM/{SIAM J\.Math\.Anal.}
\def\SMNS/{Selecta Math., New Series}
\def\TMP/{Theor\.Math\.Phys.}
\def\ZNS/{Zap\.nauch\.semin.\ LOMI}

\def\ASMP/{Advanced Series in Math.\ Phys.{}}

\def\Birk/{Birkh\"auser}
\def\CUP/{Cambridge University Press} \def\CUPa/{\CUP/ \publaddr Cambridge}
\def\Spri/{Springer\)-Verlag} \def\Spria/{\Spri/ \publaddr Berlin}
\def\WS/{World Scientific} \def\WSa/{\WS/ \publaddr Singapore}

\newbox\lefthbox \newbox\righthbox

\let\sectsep. \let\labelsep. \let\contsep. \let\labelspace\relax
\let\sectpre\relax \let\contpre\relax
\def\sf@rm{\the\Sno} \def\sf@rm@{\sectpre\sf@rm\sectsep}
\def\c@l@b#1{\contpre#1\contsep}
\def\l@f@rm{\ifd@bn@\sf@rm\labelsep\fi\labelspace\the\n@@}

\def\sectformdef{\def\sf@rm}

\let\DoubleNum\d@bn@true \let\SingleNum\d@bn@false

\def\NoNewNum{\let\writeldf\relax\def\l@b@l##1##2{\if*##1%
 \@ft@\xdef\csname @##1@##2@\endcsname{\mbox{*{*}*}}\fi}}
\def\NoNewTime{\def\todaydef##1{\def\today{##1}}
 \def\nowtimedef##1{\def\nowtime{##1}}}
\def\NoInput{\let\Input\input\let\writeldf\relax}
\def\Fixed{\NoNewTime\NoInput}

\newbox\dtlb@x
\def\DateTimeLabel{\global\setbox\dtlb@x\vbox to\z@{\ifMag\eightpoint\else
 \ninepoint\fi\sl\vss\rline\today\rline\nowtime}
 \global\headline{\hfil\box\dtlb@x}}

\def\sectfont#1{\def\s@cf@nt{#1}} \sectfont\bf
\def\subsectfont#1{\def\s@bf@nt{#1}} \subsectfont\it
\def\Entcdfont#1{\def\entcdf@nt{#1}} \Entcdfont\relax
\def\labelcdfont#1{\def\l@bcdf@nt{#1}} \labelcdfont\relax
\def\pagecdfont#1{\def\p@g@cdf@nt{#1}} \pagecdfont\relax
\def\subcdfont#1{\def\s@bcdf@nt{#1}} \subcdfont\it
\def\applefont#1{\def\@ppl@f@nt{#1}} \applefont\bf
\def\Refcdfont#1{\def\R@fcdf@nt{#1}} \Refcdfont\bf

\def\reffont#1{\def\r@ff@nt{#1}} \reffont\rm
\def\keyfont#1{\def\k@yf@nt{#1}} \keyfont\rm
\def\paperfont#1{\def\p@p@rf@nt{#1}} \paperfont\it
\def\bookfont#1{\def\b@@kf@nt{#1}} \bookfont\it
\def\volfont#1{\def\v@lf@nt{#1}} \volfont\bf
\def\issuefont#1{\def\iss@f@nt{#1}} \issuefont{no\p@@nt}
\def\volume{vol\p@@nt}

\def\adjustmid#1{\kern-#1\p@\alb\hskip#1\p@\relax}
\def\adjustend#1{\adjustnext{\kern-#1\p@\alb\hskip#1\p@}}

\newif\ifcd

\let\goodbm\relax  \let\cnngood\relax 
\let\goodbu\relax  \let\uugood\relax \def\ungood{}
 \def\vvm#1>{\ignore} \def\vvnm#1>{\ignore} \def\cnnm#1>{}
\def\cnnu#1>{} \def\vvu#1>{\ignore} \def\vvnu#1>{\ignore} \def\vsku#1>{}
\let\uline\relax  \def\uupage{}

\def\wwmgood#1:#2>{\ifMag\vv#1>\cnngood\vv#2>\vv0>\fi}
\def\vskmgood#1:#2>{\ifMag\vsk#1>\goodbm\vsk#2>\vsk0>\fi}
\def\vskugood#1:#2>{\ifUS\vsk#1>\goodbu\vsk#2>\vsk0>\fi}
\def\vskm#1:#2>{\ifMag\vsk#1>\else\vsk#2>\fi}
\def\vvmm#1:#2>{\ifMag\vv#1>\else\vv#2>\fi}
\def\vvnn#1:#2>{\ifMag\vvn#1>\else\vvn#2>\fi}
\def\nnm#1:#2>{\ifMag\nn#1>\else\nn#2>\fi}
\def\kerm#1:#2>{\ifMag\kern#1em\else\kern#2em\fi}
\def\keru#1>{\ifUS\kern#1em\fi}

\tenpoint

\Fixed

\ifMag\UStrue\hfuzz4pt\else\USfalse\fi

\loadbold
\loadeusm

\ifUS
\PaperUS
\else
\PaperA4
\fi

\font@\Beufm=eufm10 scaled 1440
\font@\Eufm=eufm8 scaled 1440
\newfam\Frakfam \textfont\Frakfam\Eufm

\sectfont{\elevenpoint\bf}
\def\Pfs#1.{\demo{Proof\/ {\rm(}\]sketch\/{\rm)} #1}\bgroup\ignore}

\def\lwg{\lsym\wedge}

\def\vsh{\raise .15ex\mbox{\>\ssize\vdash\]}}

\def\HY{\mathop{\roman{HY}}\nolimits}

\def\grad{\mathop{\roman{grad}\)}\nolimits}
\def\lbc{\mathopen{[\][}} \def\rbc{\mathclose{]\]]}}

\let\opi\varpi

\def\Kbar{\Rlap{\;\overline{\!\}\phantom K\}}\>}K}

\def\Rbar{\Rlap{\,\)\overline{\!\]\phantom R\]}\)}R}
\def\Qbar{\Rlap{\,\)\overline{\!\]\phantom Q\}}\>}Q}

\def\tbar{\bar{\]t\)}}

\def\FF{\Bbb F}

\def\mb{\bold m}

\def\vb{\bold v}
\def\wb{\bold w}
\def\xb{\bold x}

\def\ab{\boldsymbol a}
\def\bb{\boldsymbol b}

\def\mb{\boldsymbol m}

\def\Ac{\Cal A}

\def\Gc{\Cal G}

\def\Zc{\Cal Z}

\def\Rc{\check R}

\def\Ae{\eusm A}
\def\Be{\eusm B}
\def\De{\eusm D}
\def\Fe{\eusm F}
\def\Ie{\eusm I}
\def\Je{\eusm J}
\def\Ke{\eusm K}
\def\Me{\eusm M}
\def\Pe{\eusm P}
\def\Qe{\eusm Q}
\def\Se{\eusm S}
\def\Te{\eusm T}
\def\Ue{\eusm U}
\def\We{\eusm W}
\def\Xe{\eusm X}
\def\Ye{\eusm Y}

\def\Dg{\frak D}

\def\hg{\frak h}
\def\Mg{\frak M}
\def\ng{\frak n}

\def\gl{\frak{gl}}

\def\Ah{\Hat A}
\def\Bh{\Hat B}
\def\Ch{\Hat C}
\def\Dh{\Hat D}
\def\BBh{\Rlap{\Hat{\phantom{\BB}\]}\)}\BB}

\def\Gch{\Rlap{\>\Hat{\phantom{\}\Gc}}}\Gc}

\def\mh{\Hat m}

\def\Qh{\Hat Q}
\def\Rh{\Hat R}

\def\Teh{\Hat{\Te}}

\def\Qt{\Tilde Q}
\def\Rt{\Rlap{\>\Tilde{\}\phantom R\>}\}}R}

\def\psit{\Tilde\psi}
\def\Sti{\Tilde S}
\def\tti{\tilde t}

\let\wg\wedge

\def\Ekab{(E\"1_{ab}\!\lsym+E\"k_{ab})\vst{\Vws k}}
\def\Rw#1{R^{\)\wg#1}}
\def\Rwe#1{R_{\)\wg#1,\)\wg1}}
\def\Rwk{R_{\)\wg k,\)\wg1}}
\def\Rwwk{R_{\)\wg1,\)\wg k}}
\def\Rbbw#1{R^{\)\vp\wgb\raise.16ex\mbox{\ssize\smash{\wgbb}}#1}}

\def\Qw#1{Q^{\)\wg#1}}
\def\Qbw#1{\Qbar^{\)\wg#1}}
\def\Tw#1{T^{\wg#1}}
\def\Vw#1{V^{\wg#1}}
\def\Vws#1{V^{\]\wg#1}}
\def\Ww#1{W^{\wg#1}}
\def\Wws#1{W^{\]\wg#1}}
\def\Wtw#1{\bigl(\vb_1\]\wg\Ww{#1}\bigr)}
\def\Wtws#1{(\vb_1\wg\Wws{#1})\>}

\let\card\#
\def\+#1{^{\)\bra\]#1\]\ket\}}}
\def\*#1{^{\bra\]#1\]\ket\]}}
\def\##1{^{\)\lb\]#1\]\rb}}
\def\3#1{^{\]\lb\]#1\]\rb}}
\def\2#1{^{\](#1)}}
\def\4#1{^{[\)#1\)]}}
\def\0#1#2{{\vp{#2}}^{#1\}}#2}
\def\7#1#2{{\vp{#2}}^{#1}#2}
\def\6#1{_{\lb#1\rb}}
\def\9#1{_{\]\lb#1\rb}}

\def\Ixe{\Ie_{\)\xi\],Q}} \def\Jxe{\Je_{\)\xi\],\)K}}
\def\Ixik{I_{\)\xi\],\)k,Q}} 
\def\Jxik{J_{\)\xi\],\)k,\)K\)}} 
\def\Kxe{\Ke_{\)\xi\],\)K}} \def\Qxe{\Qe_{\)\xi\],Q}}
\def\Mxi{\Me_{\>\xi\],\)K}}
\def\Uxe{\Ue_{\)\xi\],\)k,Q}} \def\Uxi{\Uxe\:}
\def\Wxe{\We_{\]\xi\],\)k,\)K}} \def\Wxi{\Wxe\:}

\def\Xxi#1{\Xe_{\)\xi\],Q}^{#1}}

\def\Yti{\Yn\txinv} 
\def\Ytti{\bigl(\)\Yn/\Ynp\bigr)\txinv}

\def\txinv{\bigl[\)\txn\bigr]
\bigl[\}\bigl[(t^1_1)\vpb{-1}\}\lc(t^{N-1}_{\xi^{N-1}})\vpb{-1}\bigr]\}\bigr]}
\def\trVN{{\tr_{V^{\]\ox N\]}}\:\}}}
\def\Zxi#1{\Zc_{\)\xi\],#1,\)K}}

\def\uder{\lbc\)u\1\},\der_u\rbc}
\def\udzt{\lbc\)u\1\},\der_u\),\zt\>\rbc}

\def\one{1}
\def\ngp{\ng_{\sss+}} \def\ngm{\ng_{\sss-}} \def\ngpm{\ng_{\sss\pm}}

 \def\gln{\gl_N} \def\glnn{\gl_{N-1}} \def\glnr{\gl_N^{\>\R}}
\def\glnx{\gln[\)x\)]} \def\glnnx{\glnn[\)x\)]}
\def\ngpx{\ngp[\)x\)]}  \def\ngpmx{\ngpm[\)x\)]}
\def\ngpxo{\ngp^{\]\bra\]N-1\]\ket\}}[\)x\)]}

\def\Un{U(\gln\})} \def\Uqn{U_q(\gln\})}
\def\Unx{U(\gln\lsb\)x\)\rsb)} \def\Unnx{U(\glnn\lsb\)x\)\rsb)}
 \def\Unmx{U(\ngm\lsb\)x\)\rsb)}
\def\Yn{Y(\gln\})}
\def\Uhn{U_q(\Tilde{\gln\}})} 
\def\Uhnm{U^\mss_q\}(\Tilde{\gln\}})} 
\def\Ynn{Y(\glnn\})} 
\def\Ynp{Y_{\sss\}+\]}(\gln\})} \def\Ynx{Y_{\sss\}\x\]}(\gln\})}
\def\Ynnp{Y_{\sss\}+\]}(\glnn\})}

\def\prodl{\mathop{\overleftarrow\prod}\limits}
\def\prodr{\mathop{\overarrow\prod}\limits}

\def\Symr{\mathop{\0{R\>}{\roman S}\roman{ym}\)}\nolimits}

\def\vstt#1#2{{%
 \lower1.9\p@@\mbox{\bigr|^{(#2)}_{\raise.5\p@@\mbox{\ssize#1}}}}}

\def\Cn{\C^N} \def\Cnn{\C^{N-1}}

\def\an{a=1\lc N} \def\bn{b=1\lc N} \def\abn{1\le a<b\le N}
\def\ak{\mathchoice{a_1\lc a_k}{a_1\lc a_k}{a_1\]\lc\)a_k}{\@@PS}}

\def\Lan{\La_1\lc\La_n} 
\def\LaN{\La^{\]1}\}\lc\La^{\]N}}

\def\Mnd{M_1\]\ldots\)M_n} \def\SMox{S_{M_1}\!\]\lox S_{M_n}}
\def\Mn{M_1\lc M_n} \def\Mox{M_1\}\lox M_n} \def\Moxs{M_1\]\lox\]M_n}
\def\Mnr{M^{\)\R}_1\lc M^{\)\R}_n}
\def\Moxr{M^{\)\R}_1\!\lsym{\mathrel{{\ox_{\R}\:\}}}}\]M^{\)\R}_n}
 \def\Moxxs{M_n\]\lox\]M_1} \def\zni{z_n\lc z_1}
\def\Moxz{M_1(z_1)\lox M_n(z_n)} \def\Mozx{M_n(z_n)\lox M_1(z_1)}
\def\ik{i=\lk} \def\lk{1\ldots\)k} \def\kN{k=1\lc N} \def\koN{k=0\lc N}
 \def\inn{i=1\lc n} \def\jn{j=1\lc n}
\def\kn{k=1\lc n}

\def\lr{1\ldots\)r} \def\ur{u_1\lc u_r} 
 \def\xr{x_1\lc x_r}
  
\def\zn{z_1\lc z_n}
\def\xin{\xi^1\]\lc\xi^{N-1}} 
 
\def\vbn{\vb_1\lc\vb_N} \def\wbn{\wb_1\lc\wb_{N-1}}

\def\txn{t^1_1\lc t^{\)N-1}_{\xi^{N-1}}}
\def\ttxn{\tti^1_1\lc\tti^{\)N-1}_{\xi^{N-1}}}
\def\txin{t^1_1\lc t^1_{\xi^1}\llc t^{\)N-1}_1\}\lc t^{\)N-1}_{\xi^{N-1}}}

\def\txine{t^1_1\lc t^1_{\xi^1}}

\def\vox{v_1\]\lox v_n}

\def\cvar/{coloured \var/}
\def\wtf/{\wt/ \fn/} \def\vwtf/{vector\)-valued \wtf/}
\def\Coeffs/{The coef\-fi\-cients of the series}
\def\coeffs/{the coef\-fi\-cients of the series}
\def\wcoeff/{with coef\-fi\-cients in}
\def\ainv/{anti\)-invo\-lu\-tion}

\def\hmod/{\prfx{\hg\)}-mod\-ule}
\def\hnmod/{\prfx{\hg_N\)}-mod\-ule} \def\hnnmod/{\prfx{\hg_{N-1}\)}-mod\-ule}
\def\gnmod/{\$\gln$-mod\-ule} \def\gnxmod/{\$\glnx\)$-mod\-ule}
\def\Uqnmod/{\$\Uqn\)$-mod\-ule}
\def\Ynmod/{\$\Yn\)$-mod\-ule} \def\Ynnmod/{\$\Ynn\)$-mod\-ule}
\def\Uhnmod/{\$\Uhn\)$-mod\-ule} \def\Uhnmmod/{\$\Uhnm\)$-mod\-ule}
\def\Eqnmod/{\$\Eqn\)$-mod\-ule} \def\Eqnnmod/{\$\Eqnn\)$-mod\-ule}

\ifMag \ifUS
 \let\goodbu\goodbreak  \let\uugood\vvgood \let\uline\nl
 \let\vvu\vv \let\vvnu\vvn \let\cnnu\cnn \let\ungood\nngood \let\vsku\vsk
 \def\uupage{\ifvmode\newpage\else\vadjust{\vfill\eject}\fi} \else
 \let\goodbm\goodbreak  \let\cnngood\vvgood \let\cnnm\cnn
  \let\vvm\vv \let\vvnm\vvn  \fi
 \let\goodbreak\relax  \let\vvgood\relax
  \def\nngood{}  \fi

\csname beta.def\endcsname

\labeldef{F} {1\labelsep \labelspace 1}  {Mg}
\labeldef{F} {1\labelsep \labelspace 2}  {Me}

\labeldef{F} {3\labelsep \labelspace 1}  {inv}
\labeldef{F} {3\labelsep \labelspace 2}  {YB}
\labeldef{L} {3\labelsep \labelspace 1}  {RRA}
\labeldef{F} {3\labelsep \labelspace 3}  {AR}
\labeldef{F} {3\labelsep \labelspace 4}  {RQ}
\labeldef{L} {3\labelsep \labelspace 2}  {Rwk}
\labeldef{L} {3\labelsep \labelspace 3}  {inw}

\labeldef{S} {4} {Yn}
\labeldef{F} {4\labelsep \labelspace 1}  {Tab}
\labeldef{F} {4\labelsep \labelspace 2}  {Tabcd}
\labeldef{F} {4\labelsep \labelspace 3}  {RTT}
\labeldef{F} {4\labelsep \labelspace 4}  {Dl}
\labeldef{F} {4\labelsep \labelspace 5}  {DlT}
\labeldef{F} {4\labelsep \labelspace 6}  {EeT}
\labeldef{L} {4\labelsep \labelspace 1}  {Taa0}
\labeldef{L} {4\labelsep \labelspace 2}  {DlTaa}
\labeldef{F} {4\labelsep \labelspace 7}  {opi}
\labeldef{F} {4\labelsep \labelspace 8}  {Dlopi}
\labeldef{F} {4\labelsep \labelspace 9}  {AT}
\labeldef{F} {4\labelsep \labelspace 10} {Tw}
\labeldef{F} {4\labelsep \labelspace 11} {RTTw}
\labeldef{L} {4\labelsep \labelspace 3}  {DlTw}
\labeldef{F} {4\labelsep \labelspace 12} {qdet}
\labeldef{L} {4\labelsep \labelspace 4}  {center}
\labeldef{F} {4\labelsep \labelspace 13} {Te}
\labeldef{L} {4\labelsep \labelspace 5}  {Tekl}
\labeldef{L} {4\labelsep \labelspace 6}  {max}
\labeldef{L} {4\labelsep \labelspace 7}  {Q=1}
\labeldef{L} {4\labelsep \labelspace 8}  {Tetr}
\labeldef{F} {4\labelsep \labelspace 14} {Dgm}
\labeldef{L} {4\labelsep \labelspace 9}  {DgT}
\labeldef{F} {4\labelsep \labelspace 15} {DgmT}
\labeldef{F} {4\labelsep \labelspace 16} {DgNT}
\labeldef{F} {4\labelsep \labelspace 17} {AQA}
\labeldef{L} {4\labelsep \labelspace 10} {ADTA}
\labeldef{L} {4\labelsep \labelspace 11} {Tesym}

\labeldef{S} {5} {Ba}
\labeldef{F} {5\labelsep \labelspace 1}  {BBh}
\labeldef{L} {5\labelsep \labelspace 1}  {rega}
\labeldef{L} {5\labelsep \labelspace 2}  {thtBB}
\labeldef{F} {5\labelsep \labelspace 2}  {BB}
\labeldef{F} {5\labelsep \labelspace 3}  {TTtt}
\labeldef{F} {5\labelsep \labelspace 4}  {sixi}
\labeldef{L} {5\labelsep \labelspace 3}  {BBS}
\labeldef{F} {5\labelsep \labelspace 5}  {Ie}
\labeldef{F} {5\labelsep \labelspace 6}  {Ixik}
\labeldef{F} {5\labelsep \labelspace 7}  {Xxi}
\labeldef{L} {5\labelsep \labelspace 4}  {main}
\labeldef{F} {5\labelsep \labelspace 8}  {TeB}
\labeldef{F} {5\labelsep \labelspace 9}  {1-X}

\labeldef{S} {6} {emods}
\labeldef{F} {6\labelsep \labelspace 1}  {Tez}
\labeldef{F} {6\labelsep \labelspace 2}  {BBtz}
\labeldef{F} {6\labelsep \labelspace 3}  {Bae}
\labeldef{F} {6\labelsep \labelspace 4}  {Xxiz}
\labeldef{L} {6\labelsep \labelspace 1}  {main2}
\labeldef{F} {6\labelsep \labelspace 5}  {TeBox}
\labeldef{F} {6\labelsep \labelspace 6}  {1-Xz}
\labeldef{L} {6\labelsep \labelspace 2}  {Bsing}

\labeldef{S} {7} {Sform}
\labeldef{L} {7\labelsep \labelspace 1}  {opitau}
\labeldef{F} {7\labelsep \labelspace 1}  {opitauz}
\labeldef{L} {7\labelsep \labelspace 2}  {tauTe}
\labeldef{L} {7\labelsep \labelspace 3}  {RLM}
\labeldef{F} {7\labelsep \labelspace 2}  {RX}
\labeldef{F} {7\labelsep \labelspace 3}  {Rvw}
\labeldef{L} {7\labelsep \labelspace 4}  {RLMinv}
\labeldef{L} {7\labelsep \labelspace 5}  {RYB}
\labeldef{F} {7\labelsep \labelspace 4}  {R8}
\labeldef{F} {7\labelsep \labelspace 5}  {tauR}
\labeldef{F} {7\labelsep \labelspace 6}  {RMox}
\labeldef{F} {7\labelsep \labelspace 7}  {tauRM}
\labeldef{F} {7\labelsep \labelspace 8}  {RXXR}
\labeldef{F} {7\labelsep \labelspace 9}  {SMz}
\labeldef{L} {7\labelsep \labelspace 6}  {Shapz}

\labeldef{S} {8} {current}
\labeldef{F} {8\labelsep \labelspace 1}  {Lab}
\labeldef{F} {8\labelsep \labelspace 2}  {LL}
\labeldef{F} {8\labelsep \labelspace 3}  {detL}
\labeldef{F} {8\labelsep \labelspace 4}  {DG}
\labeldef{F} {8\labelsep \labelspace 5}  {Gc12}
\labeldef{L} {8\labelsep \labelspace 1}  {ADLA}
\labeldef{L} {8\labelsep \labelspace 2}  {Gckl}
\labeldef{L} {8\labelsep \labelspace 3}  {K=0}
\labeldef{F} {8\labelsep \labelspace 6}  {tau}
\labeldef{L} {8\labelsep \labelspace 4}  {Gcsym}
\labeldef{F} {8\labelsep \labelspace 7}  {FF}
\labeldef{L} {8\labelsep \labelspace 5}  {FFpol}
\labeldef{F} {8\labelsep \labelspace 8}  {Je}
\labeldef{F} {8\labelsep \labelspace 9}  {Jxik}
\labeldef{L} {8\labelsep \labelspace 6}  {main3}
\labeldef{F} {8\labelsep \labelspace 10} {GcF}

\labeldef{S} {9} {cemods}
\labeldef{F} {9\labelsep \labelspace 1}  {Labz}
\labeldef{F} {9\labelsep \labelspace 2}  {Gcz}
\labeldef{L} {9\labelsep \labelspace 1}  {Shap}
\labeldef{F} {9\labelsep \labelspace 3}  {Bae2}
\labeldef{F} {9\labelsep \labelspace 4}  {MxiL}
\labeldef{F} {9\labelsep \labelspace 5}  {Zxiz}
\labeldef{F} {9\labelsep \labelspace 6}  {FFtz}
\labeldef{L} {9\labelsep \labelspace 2}  {main4}
\labeldef{F} {9\labelsep \labelspace 7}  {GcFox}
\labeldef{L} {9\labelsep \labelspace 3}  {Fsing}

\labeldef{S} {10} {filt}
\labeldef{L} {10\labelsep \labelspace 1}  {gradY}
\labeldef{F} {10\labelsep \labelspace 1}  {gradopi}
\labeldef{F} {10\labelsep \labelspace 2}  {TL}
\labeldef{F} {10\labelsep \labelspace 3}  {ST}
\labeldef{F} {10\labelsep \labelspace 4}  {STx}
\labeldef{L} {10\labelsep \labelspace 2}  {SkQ}
\labeldef{F} {10\labelsep \labelspace 5}  {SeGc}
\labeldef{F} {10\labelsep \labelspace 6}  {DS}
\labeldef{F} {10\labelsep \labelspace 7}  {SD}

\labeldef{S} {11} {proof}
\labeldef{F} {11\labelsep \labelspace 1}  {pix}
\labeldef{F} {11\labelsep \labelspace 2}  {psih}
\labeldef{L} {11\labelsep \labelspace 1}  {psiY}
\labeldef{L} {11\labelsep \labelspace 2}  {phipsi}
\labeldef{L} {11\labelsep \labelspace 3}  {BB1}
\labeldef{F} {11\labelsep \labelspace 3}  {BBe}
\labeldef{F} {11\labelsep \labelspace 4}  {RN-1}
\labeldef{F} {11\labelsep \labelspace 5}  {ABCD}
\labeldef{F} {11\labelsep \labelspace 6}  {Sxw}
\labeldef{F} {11\labelsep \labelspace 7}  {BhB}
\labeldef{F} {11\labelsep \labelspace 8}  {AhB}
\labeldef{F} {11\labelsep \labelspace 9}  {DhB}
\labeldef{F} {11\labelsep \labelspace 10} {BBR}
\labeldef{F} {11\labelsep \labelspace 11} {AB}
\labeldef{F} {11\labelsep \labelspace 12} {DB}
\labeldef{F} {11\labelsep \labelspace 13} {Ract}
\labeldef{L} {11\labelsep \labelspace 4}  {RSym}
\labeldef{F} {11\labelsep \labelspace 14} {SymR}
\labeldef{L} {11\labelsep \labelspace 5}  {ADB}
\labeldef{F} {11\labelsep \labelspace 15} {ABB}
\labeldef{F} {11\labelsep \labelspace 16} {DBB}
\labeldef{F} {11\labelsep \labelspace 17} {Qbar}
\labeldef{F} {11\labelsep \labelspace 18} {TAD}
\labeldef{F} {11\labelsep \labelspace 19} {BBB}
\labeldef{F} {11\labelsep \labelspace 20} {TeBB}
\labeldef{F} {11\labelsep \labelspace 21} {BeQ}
\labeldef{L} {11\labelsep \labelspace 6}  {DRR}
\labeldef{L} {11\labelsep \labelspace 7}  {ARR}
\labeldef{F} {11\labelsep \labelspace 22} {Apsih}
\labeldef{F} {11\labelsep \labelspace 23} {AhRpsih}
\labeldef{F} {11\labelsep \labelspace 24} {Dres}
\labeldef{F} {11\labelsep \labelspace 25} {TeX}
\labeldef{F} {11\labelsep \labelspace 26} {TeB2}
\labeldef{F} {11\labelsep \labelspace 27} {UxikQ}

\labeldef{S} {12} {proof2}
\labeldef{L} {12\labelsep \labelspace 1}  {psiU}
\labeldef{L} {12\labelsep \labelspace 2}  {phipsi2}
\labeldef{L} {12\labelsep \labelspace 3}  {FF1}
\labeldef{F} {12\labelsep \labelspace 1}  {FFe}
\labeldef{L} {12\labelsep \labelspace 4}  {gradB}
\labeldef{F} {12\labelsep \labelspace 2}  {SeX}
\labeldef{F} {12\labelsep \labelspace 3}  {Pel}
\labeldef{L} {12\labelsep \labelspace 5}  {Bew}
\labeldef{F} {12\labelsep \labelspace 4}  {BP}
\labeldef{F} {12\labelsep \labelspace 5}  {Fe}
\labeldef{L} {12\labelsep \labelspace 6}  {T1S2}
\labeldef{F} {12\labelsep \labelspace 6}  {GcZ}

\labeldef{S} {\char 65} {naive}
\labeldef{F} {\char 65\labelsep \labelspace 1}  {DgNTz}
\labeldef{F} {\char 65\labelsep \labelspace 2}  {detLz}
\labeldef{L} {\char 65\labelsep \labelspace 1}  {limit}
\labeldef{F} {\char 65\labelsep \labelspace 3}  {TLlim}
\labeldef{F} {\char 65\labelsep \labelspace 4}  {SGlim}
\labeldef{F} {\char 65\labelsep \labelspace 5}  {Dlim}
\labeldef{F} {\char 65\labelsep \labelspace 6}  {Slim}
\labeldef{F} {\char 65\labelsep \labelspace 7}  {BFlim}
\labeldef{F} {\char 65\labelsep \labelspace 8}  {QKlim}
\labeldef{F} {\char 65\labelsep \labelspace 9}  {Mlim}

\labeldef{S} {\char 66} {dyn}
\labeldef{L} {\char 66\labelsep \labelspace 1}  {dynHam}

\labeldef{S} {\char 67} {real}
\labeldef{L} {\char 67\labelsep \labelspace 1}  {realZ}
\labeldef{L} {\char 67\labelsep \labelspace 2}  {singR}
\labeldef{L} {\char 67\labelsep \labelspace 3}  {SMz+}
\labeldef{L} {\char 67\labelsep \labelspace 4}  {realX}

\document

\center
\hrule height 0pt
\vsk.3>

{\twelvepoint\bf \bls1.2\bls
Bethe Eigenvectors of Higher Transfer Matrices
\par}
\vsk1.4>
\=
E\&Mukhin$^{\,*}$, \VT/$^{\,\star}$, \;and \;\Varch/$^{\,\diamond}$
\vsk1.4>
{\it
$^{*,\star}\}$\DMS/\\\IUPUI/\\\IUPUIaddr/
\vsk.24>
$^\star\}$\home/\\\homeaddr/
\vsk.5>
$^\diamond\]$\UNC/\\ \UNCaddr/
\vsk1.5>
\sl April 2006}
\endcenter

\ftext{\=\bls12pt\vsk-.1>\nt
$\]{\vru1.8ex>}^*${\tenpoint\sl E-mail\/{\rm:} mukhin\@math.iupui.edu}\vv.1>\nl
$\]^\star\)$\support/ by RFFI grant 05\)\~\)01\~\)00922\nl
\hp{$\]^\star\)$}%
{\tenpoint\sl E-mail\/{\rm:} \imail/\,{\rm,} \homemail/}\vv.1>\nl
${\]^\diamond\)}$\support/ NSF grant DMS\)\~\)0244579\nl
\vv-1.2>\hp{$^\diamond$}{\tenpoint\sl E-mail\/{\rm:} \avemail/}}

\vsk1.5>
\Abstract
We consider the \XXX/-\)type and Gaudin quantum integrable models associated
with the Lie algebra $\gln$. The models are defined on a tensor product
${\Mox}$ of \irr/ \gnmod/s. For each model, there exist $N$ one\>-parameter
families of commuting operators on ${\Mox}$, called the \tmcs/. We show that
the \Bv/s for these models, given by the algebraic nested \Ba/, are \egv/s
of higher \tmcs/ and compute the corresponding \eva/s.
\endAbs

\vsk.5>\vsk0>

\Sect{Introduction}
In this paper we consider quantum integrable models associated with
the Lie algebra $\gln$, the \XXX/-\)type and Gaudin models. The models are
defined on a tensor product ${\Mox}$ of \irr/ \gnmod/s. For each model,
there exist $N$ one\>-parameter families of mutually commuting operators on
${\Mox}$, called the \tmcs/. In this paper the \tmcs/ for the \XXX/-\)type
model are denoted by $\Te^{\Moxs}_{1,Q}(u\);z)\lc\Te^{\Moxs}_{N\},Q}(u\);z)$,
see Section~\[:emods], and the \tmcs/ for the Gaudin model are denoted by
$\Gc^{\Moxs}_{1,K}(u\);z)\lc\Gc^{\Moxs}_{N\},K}(u\);z)$, see
Section~\[:cemods].
\vsk.2>
The important problem is to find \egv/s and \eva/s of the \tmcs/.
If $\Mn$ are \hwm/s, this can be done by the algebraic Bethe ansatz method.
The ${N\]=2}$ case is fairly well understood, for example, see~\cite{G}\),
\cite{KBI}\). However, the ${N\]>2}$ case is much more technically involved
and less studied.
\vsk.1>
The purpose of this paper is to obtain \egv/s and \eva/s of the \tmcs/
for ${N>2}$ using the algebraic nested \Ba/ method. For the \tmx/
$\Te^{\Moxs}_{1,Q}(u\);z)$ this was done in~\cite{KR1}\). The \tmx/
$\Te^{\Moxs}_{N\},Q}(u\);z)$ comes from the central elements of the Yangian
$\Yn$ and is proportional to the identity operator on the whole tensor product
${\Mox}$ \cite{KS}\). In this paper we extend the technique developed
in \cite{KR1} to show that the \Bv/s are \egv/s of the higher \tmcs/
$\Te^{\Moxs}_{k,Q}(u)$, ${k>1}$, and to compute the \eva/s, see
Theorems~\[main] and~\[main2]. The formula for \eva/s of the higher \tmcs/
on \Bv/s were conjectured in \cite{KR2}\). However, no proof can be found
in literature.
\vsk.1>
For the Gaudin model, the \tmx/ $\Gc^{\Moxs}_{1,K}(u\);z)$ comes from
the central elements of the current algebra $\glnx$ and is proportional
to the identity operator on the whole tensor product ${\Mox}$. The \tmx/
$\Gc^{\Moxs}_{2,K}(u\);z)$ is the generating \fn/ of the famous Gaudin
Hamiltonians which also appear as the operators in the \KZv/ \eq/s.
Eigenvectors and \eva/s of $\Gc^{\Moxs}_{2,K}(u\);z)$ can be obtained from
those of $\Te^{\Moxs}_{1,Q}(u\);z)$ by taking a suitable limit. For more
straightforward approach see~\cite{Ju}\). For a construction of \egv/s and
\eva/s of the Gaudin Hamiltonians as a quasiclassical limit of \sol/s to
the \KZv/ \eq/s see~\cite{RV}\).
\vsk.1>
The \tmcs/ for the \XXX/-\)type model were constructed in \cite{KS} more
than 20~years ago. Suprisingly, an explicit form of the higher \tmcs/
for the Gaudin model became known only recently. In this paper we use the
approach suggested in~\cite{T} that gives the \tmcs/ of the Gaudin model as a
suitable limit of proper linear combinations of the \tmcs/ of the \XXX/-\)type
model. We trace this limit in the proof of the Theorem~\[main]. This allows us
to show that the Bethe vectors of the Gaudin model are \egv/s of the higher
\tmcs/ and to compute the \eva/s; see Theorems~\[main3] and~\[main4].
\vsk.1>
There is another approach to higher \tmcs/ of the Gaudin model and the Bethe
ansatz, based on the \rep/ theory of affine Lie algebras at the critical level
\cite{FFR}\). It is interesting and important to identify the \tmcs/
$\Gc^{\Moxs}_{1,K}(u\);z)\lc\Gc^{\Moxs}_{N\},K}(u\);z)$ among the higher
Gaudin Hamiltonians constructed in~\cite{FFR}\).
\vsk.5> 
Following \cite{T}\), it is convenient to combine the \tmcs/
of the \XXX/-\)type model into a single \dif/ operator
\vvn.2>
$$
\Dg^{\Moxs}_{N,Q}(u\);\der_u\);z)\,=\,
1\)-\)\Te^{\Moxs}_{1,Q}(u\);z)\,e^{-\der_u}
\]\lsym+\)(-1)^N\,\Te^{\Moxs}_{N\},Q}(u\);z)\,e^{-N\der_u}\),
$$
where $\der_u\}=\]d/du$, cf.~\(DgNTz)\), \(DgNT)\). The operator
$\Dg^{\Moxs}_{N,Q}(u\);\der_u\);z)$ acts on \fn/s of $u$ with values
in $\Mox\)$. Then any common \egv/ ${v\in\Mox}$ of the \tmcs/,
\vvn.2>
$$
\Te^{\Moxs}_{k,Q}(u\);z)\,v\,=\,\Ae^v_k(u)\,v\,,\qqq\Rlap{\kN\,,}
\vv.2>
$$
where $\Ae^v_1(u)\lc\Ae^v_{\]N}(u)$ are suitable scalar \fn/s,
defines the scalar \dif/ operator
$$
\Mg^v(u\);\der_u)\,=\,
1\)-\)\Ae^v_1(u)\,e^{-\der_u}\]\lsym+\)(-1)^N\)\Ae^v_{\]N}(u)\,e^{-N\der_u}\).
$$
If $v=\BB^{\vox}_{\)\xi}(t\);z)$ is the \Bv/ of the \XXX/-\)type model,
cf.~\(BBtz)\), \(BB)\), then the operator $\Mg^v(u\);\der_u)$ coincides with
the operator $\Mg_{\>\xi\],Q}(u\),\]\der_u\);t\);z\);\La)$, cf.~\(1-Xz)\),
which is written as a product,
$$
\Mg_{\>\xi\],Q}(u\),\]\der_u\);t\);z\);\La)\,=\,
\bigl(1-\Xxi{\)1}(u\);t\);z\);\La)\,e^{-\der_u}\bigr)\)\ldots\)
\bigl(1-\Xxi N(u\);t\);z\);\La)\,e^{-\der_u}\bigr)\,,
\vv.2>
\Tag{Mg}
$$
where the \fn/s $\Xxi{\>a}(u\);t\);z\);\La)$, see~\(Xxiz)\),
are explicitly defined in terms of the \sol/ $t$ of the \Bae/s \(Bae)\).
This gives a formula for the corresponding \eva/s
$\Ae^v_1(u\);z)\lc\Ae^v_{\]N}(u)$.
\vsk.5> 
It is also convenient to combine the \tmcs/ of the Gaudin model into
a single \difl/ operator
\vvn.2>
$$
\De^{\Moxs}_{\]K}(u\);\der_u\);z)\,=\,
\der_u^N\}-\)\Gc^{\Moxs}_{1,K}(u\);z)\,\der_u^{N\]-1}\]\lsym+\)
(-1)^N\,\Gc^{\Moxs}_{N\},K}(u\);z)\,,
$$
cf.~\(detLz)\), \(detL)\). Then any common \egv/ ${v\in\Mox}$ of the \tmcs/,
\vvn.1>
$$
\Gc^{\Moxs}_{k,K}(u\);z)\,v\,=\,\Zc^{\)v}_k(u)\,v\,,\qqq\Rlap{\kN\,,}
\vv.1>
$$
where $\Zc^{\)v}_{\)1}(u)\lc\Zc^{\)v}_N(u)$ are suitable scalar \fn/s,
defines the scalar \difl/ operator
$$
\Me^{\)v}(u\);\der_u)\,=\,\der_u^N\}-\)\Zc^{\)v}_{\)1}(u)\,\der_u^{N\]-1}
\]\lsym+\)(-1)^N\,\Zc^{\)v}_N(u)\,.
$$
If $v=\FF^{\vox}_{\)\xi}(t\);z)$ is the \Bv/ of the Gaudin model,
cf.~\(FFtz)\), \(FF)\), then the operator $\Me^{\)v}(u\);\der_u)$
coincides with the operator $\Mxi(u\),\der_u\);t\);z\);\La)$, cf.~\(MxiL)\),
which is written as a product,
$$
\Mxi(u\),\der_u\);t\);z\);\La)\,=\,
\bigl(\der_u\]-\Ye^{\)1}(u\);t\);z\);\La)\bigr)\)\ldots\)
\bigl(\der_u\]-\Ye^N(u\);t\);z\);\La)\bigr)\,,
\Tag{Me}
$$
where the \fn/s $\Ye^a(u\);t\);z\);\La)$ are explicitly defined in terms of
the \sol/ $t$ of the \Bae/s \(Bae2)\). This gives a formula for the \eva/s
$\Zc^{\)v}_{\)1}(u)\lc\Zc^{\)v}_{\]N}(u)$.
\vsk.5> 
The Gaudin model can be obtained from the \XXX/-\)type model in the limit
\vvn.1>
as $u\),z$ tend to infinity. In this limit, the \dif/ operators
$\Dg^{\Moxs}_{N,Q}(u\);\der_u\);z)$ and
$\Mg_{\>\xi\],Q}(u\),\]\der_u\);t\);z\);\La)$ turn into the \difl/ operators
$\De^{\Moxs}_{\]K}(u\);\der_u\);z)$ and $\Mxi(u\),\der_u\);t\);z\);\La)$,
\vvn-.1>
$$
\gather
\Dg^{\Moxs}_{N,Q}(\eps\1u\),\eps\)\der_u\);\eps\1\]z)\,=\,
\eps^{\)N}\>\De^{\Moxs}_{\]K}\](u\);\der_u\);z)\)+\)O(\eps^{\)N\]+1})\,,
\\
\nn6>
\Mg^{\Moxs}_{N,Q}(\eps\1u\),\eps\)\der_u\);\eps\1t\);\eps\1\]z\);\La)\,=\,
\eps^{\)N}\)\Me^{\Moxs}_{\]K}\](u\);\der_u\);t\);z\);\La)\)+\)
O(\eps^{\)N\]+1})\,,
\endgather
$$
as ${Q=1+\eps\)K}$ and ${\eps\to 0}$. Moreover, the \Bae/s \(Bae) of
the \XXX/-\)type model become in this limit the \Bae/s \(Bae2) of the Gaudin
model, and the \Bv/ $\BB^{\vox}_{\)\xi}(t\);z)$ turns into the \Bv/
$\FF^{\vox}_{\)\xi}(t\);z)$, cf.~Theorem~\[limit].
\vsk.5> 
The operators \(Mg) and \(Me) of the \XXX/-\)type and Gaudin models appeared
in \cite{MV2} and \cite{MV1}\), \resp/. Those operators were called in
\cite{MV1}\), \cite{MV2} the fundamental operators of the corresponding \sol/s
of the \Bae/s. It was shown therein that for \fd/ \gnmod/s $\Mn$, the kernel of
a fundamental operator is generated by \pol/s only and the singular points of
the fundamental operator were described. These facts and Schubert calculus
imply that the number of \sol/s of the corresponding \Bae/s (properly counted)
does not exceed the dimension of ${\Mox}$ \cite{MV1}\), \cite{MV2}\). That was
a step to the proof of the completeness of the \Ba/ conjecture, which says that
the \Bv/s form a basis of ${\Mox}$ under certain conditions.
\vsk.5> %
An important observation for applications is the fact that the \tmcs/ are \sym/
\wrt/ a suitable \sym/ bilinear form on ${\Mox}$, see Theorems~\[Shap]
and~\[Shapz]\). Under certain assumptions, this fact allows us to conclude that
the eigenvalues of the \tmcs/ take real values for real $u$. In particular,
this implies that the \pol/s generating the kernel of the fundamental operator
of a \sol/ of the \Bae/s can be chosen to have real coefficients only,
see applications to real algebraic geometry in \cite{MTV}\).

\Sect{Basic notation}
We will be using the standard superscript notation for embeddings of tensor
factors into tensor products. If $\Ac_1\lc\Ac_k$ are unital associative
algebras, and $a\in\]\Ac_i$, then
$$
a\"i\)=\,\one^{\ox(i\)-1)}\ox\)a\)\ox\)\one^{\ox(k-i)}\)\in\>\Ac_1\lox\Ac_k\,.
\vv-.1>
$$
If $a\in\]\Ac_i$ and $b\in\]\Ac_j$, then $(a\ox b)\"{ij}\]=\)a\"i\>b\"j\}$,
etc.
\Ex
Let $k=2$. Let $\Ac_1\),\)\Ac_2$ be two copies of the same algebra $\Ac$.
Then for any $a\),b\in\]\Ac$ we have $a\"1\}=a\ox\one$, $b\"2\}=\one\ox b$,
$(a\ox b)\"{12}\]=\)a\ox b$ and $(a\ox b)\"{21}\]=\)b\ox a$.
\enddemo
Fix a positive integer $N$.
All over the paper we use the convention ${V\}=\Cn}\}$, and we identify
elements of $\End(V)$ with ${N\!\x\!N}\}$ matrices using the standard basis
of $\Cn\}$.
\vsk.1>
Let $e_{ab}$, $a\),\bn$, be the standard
generators of the Lie algebra $\gln\}$:
$$
[\)e_{ab}\,,\)e_{cd}\)]\,=\,\dl_{bc}\,e_{ad}-\)\dl_{ad}\,e_{cb}\,.
$$
A vector $v$ in a \gnmod/ is a \wt/ vector of \wt/ $(\LaN)$ if
${e_{aa}\)v=\La^av}$ for any $\an\}$. A vector $v$ is called
a \em{singular vector} if ${e_{ab}\)v=0}$ for any $\abn\}$.
\vsk.2>
Let ${\AA\6k\}\in\)\End(V^{\ox k})}$ be the skew-symmetrization projector:
\vvn-.2>
$$
\AA\6k v_1\]\lox v_k\>=\,{1\over k\)!}\,
\sum_{\si\in S_k\!}\>\sign(\si)\,v_{\si_1}\!\lox v_{\si_k}\,.
\vv-.4>
$$
We denote by $\Vw k\}$ the image of \)$\AA\6k\!$, and for any ${Q\in\End(V)}$
\vvn-.1>
set ${\Qw k\}=\)Q^{\)\ox k}\vst{\Vws k}}\}$. In particular, $\Qw N\}=\det\)Q$.
\vsk.1>
The space $V\}$ is considered as a \gnmod/ with the natural action:
$e_{ab}\)\map E_{ab}$, where ${E_{ab}\in\End(V)}$ is the matrix with
the only nonzero entry equal to $1$ at the intersection of the \]\$a\)$-th row
and \]\$b\)$-th column. The \gnmod/ $V\}$ is called the \em{vector \rep/}.
The space $\Vw k\}$ is a \gnmod/ with the action $e_{ab}\)\map\Ekab\}$.
\vsk.1>
We will use the following notation for products of noncommuting factors:
\vvn-.1>
$$
\alignat2
& \)\prodr_{1\le i\le k}\!X_i\>=\,X_1\ldots\)X_k\,, &&
\)\prodl_{1\le i\le k}\!X_i\>=\,X_k\ldots\)X_1\,,
\\
\nn4>
& \prodr_{1\le i<j\le k}=\)\prodr_{1\le i\le k}\,\prodr_{i<j\le k} &\qqq &
\prodl_{1\le i<j\le k}=\)\prodl_{1\le j\le k}\,\prodl_{1\le i<j}
\endalignat
$$

\Sect{\Rms/}
Let ${P\)=\!\sum_{a,b=1}^N\!E_{ab}\ox E_{ba}}$. It is the flip map:
\vvn.06>
${P(x\ox y)\)=\)y\ox x}$ for any ${x\),y\in V}\}$,
and $A\"{21}\]=\)P\]A\>P$ for any $A\in\End(V^{\ox2})$.
\vsk.2>
The \em{\rat/ \Rm/} is $R(u)\)=\)u\)+P\in\)\End(V^{\ox 2})$.
It satisfies the inversion relation
$$
R(u)\>R\"{21}(\]-\)u)\,=\,1-u^2
\Tag{inv}
$$
and the \YB/
\vvn.3>
$$
R\"{12}(u-v)\>R\"{13}(u)\>R\"{23}(v)\,=\,
R\"{23}(v)\>R\"{13}(u)\>R\"{12}(u-v)\,.
\vv.3>
\Tag{YB}
$$
\Rem
The \Rm/ $R(u)$ is \sym/: $R(u)\)=\)R\"{21}(u)$. However, we will not use
this property and will write all relations, for instance \(inv)\),
in the form that naturally extends to more general \Rms/.
\enddemo
Below we use the well known \em{fusion procedure} to define \Rms/
$\Rw{k,\wg l}(u)$ acting in the tensor products $\Vw k\}\ox\Vw{\)l}\}$,
and describe their properties. The proofs can be found in
\cite{MNO, Sections~1\),\,2\)}\), \cite{N, Sections 2\),\,3\)}\). Notice that
our notation differs from those in \cite{MNO}\), \cite{N}\), instead of $R(u)$
the \Rm/ $1-P/u$ is used therein.
\Lm{RRA}
$\dsize{\kern-.5em\prodr_{1\le i<j\le k}\!\!R\"{ij}(i-j)\,=
\!\!\prodl_{1\le i<j\le k}\!\!R\"{ij}(i-j)\,=\,
\AA\6k}\>(\]-1)^k\)\prod_{j=1}^k\,(-\)j)^{k-j+1}\]$.
\vvn.5>
\endpro
The \YB/ \(YB) and Lemma~\[RRA] imply that
\vvn.1>
$$
\align
\AA\"{\lk}\6k \AA\"{k+1\ldots\)k+\)l}\6l &
\prodr_{1\le i\le k}\,\prodl_{1\le j\le l}\!R\"{i,j+k}(u+i-j-k+l\))\,={}
\Tagg{AR}
\\
\nn6>
{}=\)\biggl(\){} &
\prodl_{1\le i\le k}\,\prodr_{1\le j\le l}\!R\"{i,j+k}(u+i-j-k+l\))\}\biggr)
\ \AA\"{\lk}\6k \AA\"{k+1\ldots\)k+\)l}\6l\).\kern-1.8em
\\
\cnn-.5>
\endalign
$$
Set
\vvn-.5>
$$
\Rw{k,\)\wg\)l}(u)\,=
\prodl_{1\le i\le k}\,\prodr_{1\le j\le l}\!R\"{i,j+k}(u+i-j-k+l\))
\vst{\Vw k\]\ox\Vw l}\in\End(\Vw k)\ox\End(\Vw{\)l})\,.
$$
The introduced \Rms/ satisfy the inversion relation
\vvn-.3>
$$
\Rw{k,\)\wg\)l}(u)\>\bigl(\Rw{l,\)\wg\)k}(\]-\)u)\bigr)\2{21}\)=\,
\prod_{i=1}^k\,\prod_{j=1}^l\,\bigl(1-(u-i+j)^2\)\bigr)\,,
\vv-.6>
$$
and the \YB/
$$
\align
\bigl(\Rw{k,\)\wg\)l}(u-v)\bigr)\2{12} &
\bigl(\Rw{k,\)\wg m}(u)\bigr)\2{13}
\bigl(\Rw{\)l,\)\wg m}(v)\bigr)\2{23}={}
\\
\nn3>
{}=\,{} &
\bigl(\Rw{\)l,\)\wg m}(v)\bigr)\2{23}
\bigl(\Rw{k,\)\wg m}(u)\bigr)\2{13}
\bigl(\Rw{k,\)\wg\)l}(u-v)\bigr)\2{12}\,.
\endalign
$$
In addition, for any $Q\in\End(V)$,
$$
\bigl[\Rw{k,\)\wg\)l}(u)\>,\)\Qw k\!\ox\)\Qw l\>\bigr]\)=\,0\,.
\Tag{RQ}
$$
\vsk-.3>
Set
\vvn-1.2>
$$
\gather
\Rwk(u)\,=\,u\>+\sum_{a,b=1}^N\>\Ekab\}\ox E_{ba}\in\)\End(\Vw k)\ox\End(V)\,,
\\
\nn-6>
\Text{and}
\nn-2>
\Rwwk(u)\,=\,u+k-1\,+\sum_{a,b=1}^N\>E_{ab}\ox
(E\"1_{ba}\!\lsym+E\"k_{ba})\vst{\Vws k}\in\)\End(V)\ox\End(\Vw k)\,.
\endgather
$$
In particular, $\Rwe{N\}}(u)\)=\)u+1$ and $R_{\)\wg1,\)\wg N}(u)\)=\)u+N$.
\Lm{Rwk}
$\dsize\Rw{k,\)\wg1}(u)\)=\)\Rwk(u)\>\prod_{i=1}^{k-1}\>(u-i)\,,\quad
\Rw{1,\)\wg k}(u)\)=\)\Rwwk(u)\>\prod_{i=0}^{k-2}\>(u+i)$.
\vvn.3>
\endpro
\Cr{inw}
$\Rwk(u)\>\bigl(\Rwwk(\]-\)u)\bigr)\2{21}\]=\>(u+1)\>(k-u)$.
\endpro

\Sect[Yn]{Yangian $\Yn$}
In this section we collect some known facts from the \rep/ theory of
the Yangian $\Yn$ that will be used in the paper. We refer the reader to
reviews \cite{MNO} and \cite{Mo} for proofs and details. Notice that the series
$T_{ab}(u)$ in \(Tab) corresponds to the series $T_{ba}(u)$ in~\cite{MNO}\),
\cite{Mo}\).
\vsk.2>
The Yangian $\Yn$ is the unital associative algebra with generators
$\{T_{ab}\#s\}$, $a\),\bn$ and $s=1,2,\ldots$.
Organize them into generating series:
$$
T_{ab}(u)\,=\,\dl_{ab}\,+\sum_{s=1}^\8\,T_{ab}\#s\)u^{-s}\),
\qqq\Rlap{a\),\bn\,.}
\vv-.1>
\Tag{Tab}
$$
The defining relations in $\Yn)$ have the form
\vvn.2>
$$
(u-v)\>\bigl[\)T_{ab}(u)\>,T_{cd}(v)\)\bigr]\,=\,
T_{cb}(v)\>T_{ad}(u)-T_{cb}(u)\>T_{ad}(v)\,,
\vv.2>
\Tag{Tabcd}
$$
for all ${a,b,c,d=1\lc N}$.
\vsk-.1>
Combine all series \(Tab) together into a series
\vvn-.16>
${T(u)\>=\!\sum_{a,b=1}^N\!E_{ab}\ox T_{ab}(u)}$ \wcoeff/
${\End(\Cn)\ox\alb\Yn}$. Relations \(Tabcd) can be written as
\vvn.16>
the following equality for series \wcoeff/
\,$\End(\Cn)\ox\End(\Cn)\ox\Yn\)$:
\vvn.2>
$$
R\"{12}(u-v)\>T\"{13}(u)\>T\"{23}(v)\,=\,
T\"{23}(v)\>T\"{13}(u)\>R\"{12}(u-v)\,.
\Tag{RTT}
$$
\vsk.2>
The Yangian $\Yn$ is a Hopf algebra. In terms of generating series \(Tab)\),
the coproduct $\Dl:\Yn\to\Yn\ox\Yn$ reads as follows:
\vvn-.4>
$$
\Dl\bigl(T_{ab}(u)\bigr)\,=\,\sum_{c=1}^N\,T_{cb}(u)\ox T_{ac}(u)\,,
\qqq\Rlap{a\),\bn\,.}\kern-2em
\vv-.4>
\Tag{Dl}
$$
This formula amounts to an equality
$$
(\id\ox\Dl)\>\bigl(T(u)\bigr)\)=\)T\"{13}(u)\>T\"{12}(u)
\vv.1>
\Tag{DlT}
$$
for series \wcoeff/ $\End(V)\ox\Yn\ox\Yn$.
\vsk.1>
There is a one\>-parameter family of \aut/s ${\rho_x\}:\Yn\to\Yn}$ that is
defined in terms of the series $T(u)$ by the rule
${\rho_x\bigl(T(u)\bigr)\)=\>T(u-x)}$;
in \rhs/, $(u-x)\1\}$ has to be expanded as a power series in $u\1\}$.
\vsk.2>
The Yangian $\Yn$ contains the universal enveloping algebra $\Un$ as a Hopf
subalgebra. The embedding is given by \)$e_{ab}\map T_{ba}\#1\]$ for any
$a\),\bn$. We identify $\Un$ with its image in $\Yn$ under this embedding.
It is clear from relations \(Tabcd) that for any $a\),\bn$,
\vvn-.3>
$$
\bigl[\)E_{ab}\ox\one+\one\ox e_{ab}\>,T(u)\)\bigr]\>=\>0\,.
\Tag{EeT}
$$
\vsk-.1>
The \em{\ehom/} ${\epe:\Yn\to\Un}$ is given by the rule
${\epe\):\)T_{ab}\#1\]\map\>e_{ba}}$ for any $a\),\bn$, and
$\epe\):\)T_{ab}\#s\]\map\>0$ for any $s>1$ and all $a\),b$.
Both the \aut/s $\rho_x$ and the \hom/ $\epe$ restricted to the subalgebra
$\Un$ are the identity maps.
\vsk.1>
For a \gnmod/ $M$ denote by $M(x)$ the \Ynmod/ induced from $M$ by the \hom/
$\epe\circ\rho_x$. The module $M(x)$ is called an \em{\emod/} over $\Yn$
with the \em{\epoint/} $x$.
\vsk.2>
Denote by $\Ynp$ the left ideal in $\Yn$ generated by \coeffs/ $T_{ab}(u)$ for
$1\le b<a\le N\}$. If $A\>,B\in\Yn$ and $A-\]B\in\Ynp$, then we will write
$A\)\simeq B$. Relations \(Tabcd) imply the following lemma.
\Lm{Taa0}
\atem
For any $a\),\bn\},$ \coeffs/ ${[\)T_{aa}(u)\>,T_{bb}(v)\)]}$ belong to $\Ynp$.
\bitem
For any ${Z\in\Ynp}$ and any ${\an}\}$, \coeffs/ $Z\>T_{aa}(u)$ belong to
$\Ynp$.
\uugood
\endpro
\nt
Therefore, for any $\an\}$, \coeffs/ $T_{aa}(u)$ act on the space $\Yn/\Ynp$
by multiplication from the right, and those actions commute.
\Lm{DlTaa}
\Coeffs/ \)$\Dl\bigl(T_{aa}(u)\bigr)\)-\)T_{aa}(u)\ox T_{aa}(u)$, $\an$,
belong to $\Yn\ox\Ynp\,+\>\Ynp\ox\Yn$.
\endpro
\nt
The proof is straightforward, cf.~\(Dl)\).
\vsk.5>
A vector $v$ in a \Ynmod/ is called \em{singular} \wrt/ the action of $\Yn$
if ${\Ynp\>v\)=\)0}$. A singular vector $v$ that is an \egv/ for the action
of $T_{11}(u)\lc T_{N\}N}(u)$ is called a \em{\wt/ singular vector}.
\Ex
Let $M$ be a \gnmod/ and ${v\in M}$ a singular vector of \wt/ $(\LaN)$.
Then for any $x\in\C$, the vector $v$ is a \wt/ singular vector \wrt/
the action of $\Yn$ in the \emod/ $M(x)$.
\enddemo
There is a Hopf algebra \ainv/ ${\opi:\Yn\to\Yn}$ that is defined
in terms of the series $T_{ab}(u)$ by the rule
\vvn-.3>
$$
\opi\bigl(T_{ab}(u)\bigr)\>=\,T_{ba}(u)\,,\qquad a\),\bn\,.
\vv-.2>
\Tag{opi}
$$
In particular, for any ${X\]\in\Yn}$ we have
\vvn-.6>
$$
\Dl\bigl(\opi(X)\bigr)\,=\,
\bigll((\]\opi\ox\opi)\bigl(\Dl(X)\bigr)\]\bigrr)\2{21}\).
\Tag{Dlopi}
$$
\vsk-.3>
Consider the series ${T\"{k,\)k+1}(u)\>\ldots\>T\"{1,\)k+1}(u-k+1)}$ \wcoeff/
$\End(V^{\ox k})\ox\Yn$. Relations \(RTT) and Lemma~\[RRA] imply that
\vvn-.1>
$$
\align
\AA\"{\lk}\6k\,T\"{1,\)k+1}(u-k+1)\>\ldots\>{}& T\"{k,\)k+1}(u)\,={}
\Tag{AT}
\\
\nn4>
{}={}\,T\"{k,\)k+1}(u)\>\ldots\>{}&T\"{1,\)k+1}(u-k+1)\;\AA\"{\lk}\6k\).
\\
\cnn-.3>
\endalign
$$
Therefore, \coeffs/ ${T\"{k,\)k+1}(u)\>\ldots\>T\"{1,\)k+1}(u-k+1)}$ preserve
$\Vw k\}$. Define the series $\Tw k(u)$ \wcoeff/ $\End(\Vw k)\ox\Yn$:
$$
\Tw k(u)\)=\)T\"{k,\)k+1}(u)\>\ldots\>T\"{1,\)k+1}(u-k+1)\)
\vstt{\Vw k}{\lk}\}\,.
\vv-.2>
\Tag{Tw}
$$
In particular, $\Tw 1(u)\)=\)T(u)$. Formulae \(RTT) and \(AR) imply
that the introduced series satisfy the following relations:
$$
\align
\bigl(\Rw{k,\)\wg\)l}(u-v)\bigr)\2{12} &
\bigl(\Tw k(u)\bigr)\2{13}\bigl(\Tw l(v)\bigr)\2{23}\,={}
\Tag{RTTw}
\\
\nn4>
{}=\,{} & \bigl(\Tw l(v)\bigr)\2{23}\bigl(\Tw k(u)\bigr)\2{13}
\bigl(\Rw{k,\)\wg\)l}(u-v)\bigr)\2{12}\).
\\
\cnn-.2>
\endalign
$$
\Lm{DlTw}
$(\id\ox\Dl)\>\bigl(\Tw k(u)\bigr)\,=\,
\bigl(\Tw k(u)\bigr)\2{13}\>\bigl(\Tw k(u)\bigr)\2{12}\}$.
\endpro
\Pf.
The statement follows from the coproduct formula~\(DlT) and
the definition~\(Tw) of $\Tw k(u)$.
\uugood
\epf
Since the space $\Vw N\}$ is \onedim/, \coeffs/ $\Tw N(u)$ belong to the
Yangian $\Yn$. The series $\Tw N(u)$ is called the \em{quantum determinant} and
is denoted by $\qdet\)T(u)$. Explicit formulae for the quantum determinant that
follow from \(AT) and \(Tw) are
\vvn-.4>
$$
\align
\qdet T(u)\, &{}=\,\sign(\si)\)\sum_{\tau\in S_N\!}\)\sign(\tau)\;
T_{\si_1\],\)\tau_1}(u)\)\ldots\>T_{\si_N\},\)\tau_N}(u-N+1)
\Tag{qdet}
\\
\nn6>
&{}=\,\sign(\si)\)\sum_{\bb\in S_N\!}\)\sign(\tau)\;
T_{\tau_N\},\)\si_N}(u-N+1)\)\ldots\>T_{\tau_1\],\)\si_1}(u)\,,
\\
\cnn.1>
\endalign
$$
where the \perm/ ${\si\in\]S_N}$ is arbitrary. Relations~\(RTTw) and
Lemma~\[Rwk] imply that the coefficients of $\qdet T(u)$ are central in $\Yn$.
\Th{center}
\back\;\cite{MNO}
The coefficients of $\qdet T(u)$ are free generators of the center of
the Yangian $\Yn$.
\endpro
Let $\tr_W\:\]:\)\End(W)\to\)\C$ be the trace map.
For any $Q\in\End(V)$ define the series
\vvn.3>
$$
\Te_{k,Q}(u)\,=\,({\tr_{\Vws k}\:\}}\ox\id\))\bigl(\Qw k\>\Tw k(u)\bigr)\,,
\vv.3>
\Tag{Te}
$$
${\kN}\}$, \wcoeff/ $\Yn$. We set ${\Te_{0,Q}(u)\)=\)1}$ by definition.
It is easy to see that ${\Te_{N\},Q}(u)\)=\)\det\)Q\,\qdet T(u)}$. The series
$\Te_{1,Q}(u)\lc\Te_{N\},Q}(u)$ are called the \em{\tmcs/}.
They have been introduced in \cite{KS}\).
\Prop{Tekl}
\back\;\cite{KS}
We have \,$\bigl[\)\Te_{k,Q}(u)\>,\Te_{l,Q}(v)\)\bigr]\)=\)0$.
\endpro
\Pf. The statement follows from formulae \(RTTw) and \(RQ)\):
\vvn-.1>
$$
\alignat2
\Te_{k,Q}(u)\>\Te_{l,Q}(v)\, &{}=\,
({\tr_{\Vws k}\:\}}\ox{\tr_{\Vws l}\:\}}\ox\id\))\Bigl(\}{} &&
\bigl(\Qw k\>\Tw k(u)\bigr)\2{13}\bigl(\Qw l\>\Tw l(v)\bigr)\2{23}\Bigr)\,={}
\\
\nn4>
& {}=\,({\tr_{\Vws k}\:\}}\ox{\tr_{\Vws l}\:\}}\ox\id\))\Bigl(\}{} &&
\bigl(\Qw k\bigr)\21\bigl(\Qw l\bigr)\22\bigl(\Tw k(u)\bigr)\2{13}
\bigl(\Tw l(v)\bigr)\2{23}\Bigr)\,={}
\\
\nn4>
& {}=\,({\tr_{\Vws k}\:\}}\ox{\tr_{\Vws l}\:\}}\ox\id\))\Bigl(\}{} &&
\bigl(\Qw k\bigr)\21\bigl(\Qw l\bigr)\22
\bigll(\]\bigl(\Rw{k,\)\wg\)l}(u-v)\bigr)\2{12}\bigrr)^{\}-1}\x{}
\\
\nn1>
&&\Llap{{}\x\,{}} &\bigl(\Tw k(v)\bigr)\2{23}\bigl(\Tw l(u)\bigr)\2{13}
\bigl(\Rw{k,\)\wg\)l}(u-v)\bigr)\2{12}\Bigr)={}
\\
\ald
\nn4>
& {}=\,({\tr_{\Vws k}\:\}}\ox{\tr_{\Vws l}\:\}}\ox\id\))\Bigl(\}{} &&
\bigll(\]\bigl(\Rw{k,\)\wg\)l}(u-v)\bigr)\2{12}\bigrr)^{\}-1}
\bigl(\Qw k\bigr)\21\bigl(\Qw l\bigr)\22\)\x{}
\\
\nn1>
&&\Llap{{}\x\,{}} &\bigl(\Tw k(v)\bigr)\2{23}\bigl(\Tw l(u)\bigr)\2{13}
\bigl(\Rw{k,\)\wg\)l}(u-v)\bigr)\2{12}\Bigr)={}
\\
\nn4>
& {}=\,({\tr_{\Vws k}\:\}}\ox{\tr_{\Vws l}\:\}}\ox\id\))\Bigl(\}{} &&
\bigl(\Qw k\bigr)\21\bigl(\Qw l\bigr)\22\bigl(\Tw l(v)\bigr)\2{23}
\bigl(\Tw k(u)\bigr)\2{13}\Bigr)\,={}
\\
\nn4>
& {}=\,({\tr_{\Vws k}\:\}}\ox{\tr_{\Vws l}\:\}}\ox\id\))\Bigl(\}{} &&
\bigl(\Qw l\>\Tw l(v)\bigr)\2{23}\bigl(\Qw k\>\Tw k(u)\bigr)\2{13}\Bigr)\,={}
\\
\nn8>
& {}=\,\Te_{l,Q}(v)\>\Te_{k,Q}(u)\,.
\\
\cnn-.7>
\Text{\qed}
\endalignat
$$
\edemo
\Th{max}
\back\;\cite{NO}
If the matrix ${Q\in\End(V)}$ has a simple spectrum, then \coeffs/
\>${\Te_{1,Q}(u)\lc\Te_{N-1,Q}(u)}$ and \>$\qdet T(u)$ generate
a maximal commutative subalgebra in $\Yn$.
\endpro
Recall that the Yangian $\Yn$ contains $\Un$ as a subalgebra, the embedding
being given by \)$e_{ab}\map T_{ba}\#1\}$.
\Prop{Q=1}
\back\;\cite{KS}
Let $Q$ be the identity matrix. Then for any $\kN$, the coefficients of
\>$\Te_{k,Q}(u)$ commute with the subalgebra $\Un\sub\Yn$.
\endpro
\Pf.
Similarly to \(EeT)\), for any $a\),\bn\}$, we have
$$
\bigl[\)e_{ab}\vst{\Vws k}\}\ox\one+\one\ox e_{ab}\>,\Tw k(u)\)\bigr]
\>=\,0\,,
\vv-.1>
$$
which implies that for the identity matrix $Q$,
\vvn-.3>
$$
\bigl[\)e_{ab}\>,\Te_{k,Q}(u)\)\bigr]\>=\,-\>({\tr_{\Vws k}\:\}}\ox\id\))
\Bigl(\bigl[\)e_{ab}\vst{\Vws k}\}\ox\one\>,\Tw k(u)\)\bigr]\Bigr)\,=\,0
\vv-.3>
$$
due to the cyclic property of the trace.
\epf
\Lm{Tetr}
For any $m=k\lc N\}$, and any distinct $i_1\lc i_k\]\in\lb\)1\lc m\)\rb$,
we have
$$
\align
& \Te_{k,Q}(u)\,={}
\\
\nn6>
& \!\!{}=\,{m\)!\,(N\]-m)\)!\over k\)!\,(N\]-k)\)!}\;
({\tr_{V^{\]\ox m}}\:\}}\ox\id\))\bigl(Q\"{i_1}\!\]\ldots\>Q\"{i_k}\>
\>T\"{i_1,\)m+1}(u)\>\ldots\>T\"{i_k,\)m+1}(u-k+1)\;
\AA\"{1\ldots m}\6m\)\bigr)\,.\!
\\
\cnn-.5>
\endalign
$$
\endpro
\Pf.
If $m=k$ and ${i_j\]=j}$ for all ${j=1\lc k}$, the claim is \eqv/ to
formula~\(Te)\). In general, we are also using the relations
\>${P\"{ij}\)Q\"i\)T\"{i,\)m+1}(u)\)=\)T\"{j,\)m+1}(u)\,Q\"jP\"{ij}}$ and
${P\"{ij}\AA\6m\]=-\)\AA\6m\}=\AA\6mP\"{ij}}\}$, ${1\le i<j\le m}$,
for the flip maps, and the formula
\vvn-.3>
$$
\gather
(\id^{\ox k}\!\]\ox
\tr_{V^{\]\ox(\]m-k)}}\:)\>\AA\6m\>=\,
{m\)!\,(N\]-m)\)!\over k\)!\,(N\]-k)\)!}\ \AA\6k\>.
\\
\cnn->
\Text{\qed}
\cnn-.2>
\endgather
$$
\edemo
Consider the algebra $\Yn\uder$ of formal \difl/ operators;
the indeterminates $u\1\}$ and $\der_u\!$ do not commute obeying the following
relation $\der_u{\cdot}\>u^{-s}\}=u^{-s}\>\der_u\]-s\)u^{-s-1}\}$ instead.
Following \cite{T}\), for any $m=1\lc N\}$, set
\vvn-.2>
$$
\align
& \hp{{}=\,(\tr}\Dg_{m,Q}(u\),\der_u)\,={}
\Tagg{Dgm}
\\
\nn6>
& {}=\,({\tr_{V^{\]\ox m}}\:\}}\ox\id\))\Bigl(\]
\bigl(1-Q\"1\)T\"{1,\)m+1}(u)\,e^{-\der_u}\bigr)\)\ldots\)
\bigl(1-Q\"mT\"{m,\)m+1}(u)\,e^{-\der_u}\bigr)\;
\AA\"{1\ldots m}\6m\)\Bigr)\kern-.4em
\\
\cnn-.2>
\endalign
$$
and \,$\Dg_{\)0,Q}(u\),\der_u)\)=\)1$.
\Prop{DgT}
\back\;\cite{T}
For any $m=1\lc N\}$, we have
\vvn-.2>
$$
\Dg_{m,Q}(u\),\]\der_u)\,=\,{1\over(N\]-m)\)!}\;
\sum_{k=0}^m\,(-1)^k\,{(N\]-k\))\)!\over(m-k\))\)!}
\ \Te_{k,Q}(u)\,e^{-k\der_u}\).
\vv-.8>
\Tag{DgmT}
$$
\endpro
\Pf.
Expand the product under the trace in formula \(Dgm)\),
\vvn-.3>
$$
\align
& \Dg_{m,Q}(u\),\]\der_u)\,=\,
\sum_{k=0}^m\;\sum_{1\le i_1\lsym<i_k\le m\!\!}(-1)^k\)\x{}
\\
\nn8>
&\ \,\){}\x\,({\tr_{V^{\]\ox m}}\:\}}\ox\id\))\bigl(Q\"{i_1}\)T\"{i_1,\)m+1}(u)
\>\ldots\>Q\"{i_k}\)T\"{i_k,\)m+1}(u-k+1)\;\AA\"{1\ldots m}\6m\)\bigr)
\,e^{-k\der_u}\).
\kern-2em
\endalign
$$
Then Lemma~\[Tetr] yields the claim.
\epf
\Rem
Formula \(DgmT) for ${m=N}$ shows that the expression $\Dg_{N,Q}(u\),\]\der_u)$
serves as a generating \fn/ of the \tmcs/ ${\Te_{1,Q}(u)\lc\Te_{N,Q}(u)}$:
$$
\Dg_{N,Q}(u\),\]\der_u)\,=\,
\sum_{k=0}^N\,(-1)^k\,\Te_{k,Q}(u)\,e^{-k\der_u}\).
\Tag{DgNT}
$$
\enddemo
\Rem
Taking into account formula \(AT) and the equality ${(\AA\6m\])^2\]=\)\AA\6m}$,
one can show that
$$
\align
\bigl(1 &{}-Q\"1\)T\"{1,\)m\]+1}(u)\,e^{-\der_u}\bigr)\)\ldots\)
\bigl(1-Q\"mT\"{m\],\)m\]+1}(u)\,e^{-\der_u}\bigr)\;
\AA\"{1\ldots m}\6m\)={}
\Tagg{AQA}
\\
\nn5>
& \]{}=\,\AA\"{1\ldots m}\6m
\bigl(1-Q\"1\)T\"{1,\)m\]+1}(u)\,e^{-\der_u}\bigr)\)\ldots\)
\bigl(1-Q\"mT\"{m\],\)m\]+1}(u)\,e^{-\der_u}\bigr)\;\AA\"{1\ldots m}\6m\).
\kern-2.4em
\endalign
$$
Since $\AA\6N\}$ is a \onedim/ projector, formula \(AQA) implies the following
proposition.
\Prop{ADTA}
We have
$$
\align
\bigl(& {}\)\id^{\ox N}\!\]\ox\)\Dg_{\]N\],Q}\](u\),\der_u)\bigr)\;
\AA\"{1\ldots N}\6N\> \,={}
\\
\nn5>
& \!\]{}=\,\bigl(1-Q\"1\)T\"{1,\)N\]+1}(u)\,e^{-\der_u}\bigr)\)\ldots\)
\bigl(1-Q\"NT\"{N\],\)N\]+1}(u)\,e^{-\der_u}\bigr)\;
\AA\"{1\ldots N}\6N\>.
\\
\cnn-.5>
\endalign
$$
\endpro
\enddemo
Let $Q^\dag\}$ be the matrix transposed to $Q$;
for example, $E^\dag_{ab}\}=E_{ba}$.
\Prop{Tesym}
We have \,$\opi\bigl(\Te_{k,Q}(u)\bigr)\>=\>\Te_{k,Q^\dag}(u)$,
\,where $\opi$ is given by \(opi)\).
\endpro
\Pf.
The claim follows from Lemma~\[Tetr] for $\Te_{k,Q}(u)$ with ${m=k}$ and
\vvn.1>
${i_j\}=k-j+1}\]$, formula \(AT), two more observations,
\vvn.1>
$\opi\bigl(T(u)\bigr)\)=\)\bigl(T(u)\bigr)^{\}\dag}\}$ and
$\bigl(\AA\"{\lk}\6k\bigr)^{\}\dag}\}=\)\AA\"{\lk}\6k\}$,
and the standard properties of the trace:
$$
\alds
\align
& \,\opi\bigl(\Te_{k,Q}(u)\bigr)\,={}
\\
\nn5>
& {}=\,{1\over k\)!}\;({\tr_{V^{\]\ox k}}\:\}}\ox\opi)\bigl(\)
Q\"k\!\]\ldots\>Q\"1\>T\"{k,\)k+1}(u)\>\ldots\>T\"{1,\)k+1}(u-k+1)\;
\AA\"{\lk}\6k\)\bigr)\,={}
\\
\nn5>
& {}=\,{1\over k\)!}\;({\tr_{V^{\]\ox k}}\:\}}\ox\opi)\bigl(\)
Q\"k\!\]\ldots\>Q\"1\>\AA\"{\lk}\6k\>T\"{1,\)k+1}(u-k+1)\>\ldots\>
T\"{k,\)k+1}(u)\bigr)\,={}
\\
\nn5>
& {}=\,{1\over k\)!}\;({\tr_{V^{\]\ox k}}\:\}}\ox\id\))\bigl(\)
Q\"k\!\]\ldots\>Q\"1\>\AA\"{\lk}\6k\>\bigl(T\"{k,\)k+1}(u)\bigr)^{\}\dag}\]
\ldots\>\bigl(T\"{1,\)k+1}(u-k+1)\bigr)^{\}\dag}\)\bigr)\,={}
\\
\nn5>
& {}=\,{1\over k\)!}\;({\tr_{V^{\]\ox k}}\:\}}\ox\id\))\Bigl(\]\bigl(\)
Q\"k\!\]\ldots\>Q\"1\>\AA\"{\lk}\6k\>\bigl(T\"{k,\)k+1}(u)\bigr)^{\}\dag}\]
\ldots\>\bigl(T\"{1,\)k+1}(u-k+1)\bigr)^{\}\dag}\)\bigr)^{\}\dag}\Bigr)\,={}
\\
\nn5>
& {}=\,{1\over k\)!}\;({\tr_{V^{\]\ox k}}\:\}}\ox\id\))\bigl(\)
T\"{k,\)k+1}(u)\>\ldots\>{}T\"{1,\)k+1}(u-k+1)\;\AA\"{\lk}\6k\>
\bigl(Q\"k\bigr)^{\}\dag}\]\ldots\>\bigl(Q\"1\bigr)^{\}\dag}\)\bigr)\,={}
\\
\nn5>
& {}=\,{1\over k\)!}\;({\tr_{V^{\]\ox k}}\:\}}\ox\id\))\bigl(\)
\bigl(Q\"k\bigr)^{\}\dag}\]\ldots\>\bigl(Q\"1\bigr)^{\}\dag}\,
T\"{k,\)k+1}(u)\>\ldots\>{}T\"{1,\)k+1}(u-k+1)\;\AA\"{\lk}\6k\)\bigr)\,={}
\\
\nn6>
& {}=\,\Te_{k,Q^\dag}(u)\,.
\\
\cnn->
\Text{\qed}
\endalign
$$
\edemo

\Sect[Ba]{The \Ba/}
In this section we formulate the main technical result of this paper,
Theorem~\[main]. Equality \(TeB) is, in a sense, a universal formulation
of the \eva/ problem for Bethe \egv/s of the transfermatrix $\Te_{k,Q}(u)$.
Given a \Ynmod/ $M$ and a \wt/ singular vector $v\in\]M$, we will apply formula
\(TeB) to $v$ and obtain the Bethe \egv/ for the action of $\Te_{k,Q}(u)$ in
$M$ and the formula for the corresponding \eva/. In Section~\[:emods] we will
do it in detail for $M$ being a tensor product of \hw/ \emod/s and ${v\in M}$
being the tensor product of \hwv/s.
\vsk.2>
Let ${\xi=(\xin)}$ be a collection of nonnegative integers.
Set ${\xi^{\)<a}\}=\xi^1\]\lsym+\xi^{a-1}}\}$, $\an$, and
${|\)\xi\)|=\xi^1\]\lsym+\xi^{N-1}\}=\)\xi^{\)<N}}\}$.
\vvn.16>
Consider a series in $|\)\xi\)|$ \var/s $\txin$ \wcoeff/ $\Yn$:
\vvn-.4>
$$
\align
\BBh_{\)\xi}(\txn)\,&{}=\,(\tr_{V^{\]\ox|\xi|}}\:\!\}{}\ox\id\))\>\biggl(\)
T\"{1,\)|\xi|+1}(t^1_1)\ldots\)T\"{|\xi|,\)|\xi|+1}(t^{N-1}_{\xi^{N-1}})\,\x{}
\Tagg{BBh}
\\
\nn5>
& \>{}\x\prodr_{(a,i)<(b,j)\}}\]\!
R\"{\xi^{<b}\]+\)j,\)\xi^{<a}\]+\)i\)}(t^b_j\]-t^a_i)\;
E_{21}^{\)\ox\xi^1}\!\}\lox E_{N\},\)N-1}^{\)\ox\xi^{N-1}}\ox\one\]\biggr)\,,
\kern-2em
\\
\cnn-.1>
\endalign
$$
where the pairs are ordered lexicographically, $(a,i)<(b,j)$ if $a<b$,
or $a=b$ and $i<j$, and the product is taken over all two element subsets
of the set $\lb\>(c,k)\,\vert\,c=1\lc N-1\}\,,\ \,k=1\lc \xi^c\)\rb$.
Further on, we will abbreviate, $t=(\txn)$.
\vsk-.3>\vsk0>
\Lm{rega}
The series \,$\BBh_{\)\xi}(t)$ is divisible by
${\prod_{a=1}^N\prod_{\,1\le i<j\le\xi^a\!\!}\>(t^a_j\]-t^a_i\]+1)}$ in
\vvn-.3>
\uline
$\Yti$.
\endpro
\Pf.
\bls1.14\bls
Using the \YB/ \(YB) we can rearrange the factors
in the product of \Rms/ in formula \(BBh) \st/ the product
$\!{\prodr_{1\le a<N}\prodr_{\;1\le i<j\le\xi^a\!}
\]R\"{\xi^{<a}\]+\)j,\)\xi^{<a}\]+\)i\)}(t^a_j\]-t^a_i)}$ appears on the right.
Then formula $R(u)\>{E_{bc}\}\ox E_{bc}}\)=\)(u+1)\>{E_{bc}\}\ox E_{bc}}$
yields the claim.
\epf
Let $\tht:\Yn\to\Yn/\Ynp$ be the canonical projection.
\vsk-.3>\vsk0>
\Lm{thtBB}
The series \,$\tht\bigl(\)\BBh_{\)\xi}(t)\bigr)$ is divisible by
${\prod_{a=1}^{N-2}\)\prod_{b=a+2}^{N-1}\)\prod_{i=1}^{\xi^a\!}\>
\prod_{j=1}^{\xi^b\!}\,(t^b_j\]-t^a_i)}$ in
\uline
$\Ytti$.
\endpro
\nt
The proof is given in Section~\[:proof] after Proposition~\[BB1].
\vsk.5>
Set
\vvn-1.2>
$$
\BB_{\)\xi}(t)\,=\,\BBh_{\)\xi}(t)\,\,
\prod_{a=1}^{N-1}\,\prod_{1\le i<j\le\xi^a}{1\over t^a_j\]-t^a_i\]+1}\,
\prod_{1\le a<b\)<N}\,\prod_{i=1}^{\xi^a\!}\;\prod_{j=1}^{\xi^b\!}\,
{1\over t^b_j\]-t^a_i}\;,
\vv-.2>
\Tag{BB}
$$
cf.~\(BBh)\). To indicate the dependence on $N\}$, if necessary,
we will write $\BB\*N_{\)\xi}(t)$.
\Rem
The expression $\BB_{\)\xi}(t)$ serves as a creation operator to
the universal weight \fn/ and \Bv/s in tensor products of \emod/s.
See more details in Section~\[:emods].
\enddemo
\Ex
Let $N=2$ and $\xi=(\)\xi^1)$.
Then $\BB_{\)\xi}(t)\)=\)T_{12}(t^1_1)\ldots T_{12}(t^1_{\xi^1}\])$.
\enddemo
\Ex
Let $N=3$ and $\xi=(1,1)$. Then
$$
\BB_{\)\xi}(t)\,=\,T_{12}(t^1_1)\>T_{23}(t^2_1)\,+\,
{1\over t^2_1\]-t^1_1}\;T_{13}(t^1_1)\>T_{22}(t^2_1)\,.
\vv0>
$$
\enddemo
\Ex
Let $N=4$ and $\xi=(1,1,1)$. Then
\vvn-.1>
$$
\alignat2
\BB_{\)\xi}(t)\, &{}=\,T_{12}(t^1_1)\>T_{23}(t^2_1)\>T_{34}(t^3_1)\,+{}
\\
\nn4>
&{}\)+\,{1\over t^2_1\]-t^1_1}\;T_{13}(t^1_1)\>T_{22}(t^2_1)\>T_{34}(t^3_1)
\,+\,{1\over t^3_1\]-t^2_1}\;T_{12}(t^1_1)\>T_{24}(t^2_1)\>T_{33}(t^3_1)\,+{}
\\
\nn3>
&{}+\,{1\over(t^2_1\]-t^1_1)\>(t^3_1\]-t^1_1)}\,
\bigl(\)T_{14}(t^1_1)\>T_{22}(t^2_1)\>T_{33}(t^3_1)\>+\>
T_{13}(t^1_1)\>T_{24}(t^2_1)\>T_{32}(t^3_1)\bigr)\,+{}
\\
\nn3>
&& \Llap{{}+\,
{(t^2_1\]-t^1_1)\>(t^3_1\]-t^1_1)+1\over
(t^2_1\]-t^1_1)\>(t^3_1\]-t^1_1)\>(t^3_1\]-t^2_1)}\;
T_{14}(t^1_1)\>T_{23}(t^2_1)\>T_{32}(t^3_1)\,} & .
\endalignat
$$
\enddemo
In general, $\BB_{\)\xi}(t)$ is a sum of the products
\vvn-.7>
$$
T_{a_1\],\)b_1}(t^1_1)\>\ldots\>
T_{a_{|\xi|}\],\)b_{|\xi|}}(t^{N-1}_{\xi^{N-1}})\;p\)(\txn)\!
\prod_{1\le a<b\)<N}\,\prod_{i=1}^{\xi^a\!}\;
\prod_{j=1}^{\xi^b\!}\,{1\over t^b_j\]-t^a_i}\;
\vv-.4>
\Tag{TTtt}
$$
with various $a_1\lc a_{|\xi|}\]$, $b_1\lc b_{|\xi|}\]$,
and \pol/s $p\)(\txn)$.
\vsk.1>
The direct product of the \symg/s ${S_{\xi^1}\}\lx S_{\xi^{N-1}}}$ acts
on expressions in $|\)\xi\)|$ \var/s, permuting the \var/s with the same
superscript:
$$
\si^1\!\lx\si^{N-1}:\)f(\txn)\,\map\>
f(t^1_{\si^{\]1}_1}\lc t^1_{\si^{\]1}_{\}\xi^1}};\;\ldots\;;
t^{\)N-1}_{\si^{\]N-1}_1}\}\lc t^{\)N-1}_{\si^{\]N-1}_{\}\xi^{N-1}}})\,,
\vv-.4>
\Tag{sixi}
$$
where ${\si^a\!\in S_{\xi^a}}$, $\an-1$.
\Lm{BBS}
\back\;\cite{TV1, Theorem~3.3.4\)}
The expression $\BB_{\)\xi}(t)$ is invariant under the action \(sixi)
of the group ${S_{\xi^1}\}\lx S_{\xi^{N-1}}}$.
\endpro
\vsk-.2>
Fix a diagonal matrix $Q\)=\]\sum_{a=1}^N\]Q_aE_{aa}$. Introduce the series
\vvn-.4>
$$
\align
\Ixe^{\)a,i}(t)\, &{}=\,Q_a\>T_{aa}(t^a_i)\,
\prod_{j=1}^{\xi^{a-1}\!\!}\;(t^a_i\]-t^{a-1}_j\}+1)
\,\prod_{j=1}^{\xi^a}\,(t^a_i\]-t^a_j\]-1)\,
\prod_{j=1}^{\xi^{a+1}\!\!}\;(t^a_i\]-t^{a+1}_j)\,+{}
\Tagg{Ie}
\\
& {}\)+\,Q_{a+1}\>T_{a+1,a+1}(t^a_i)\,
\prod_{j=1}^{\xi^{a-1}\!\!}\;(t^a_i\]-t^{a-1}_j)
\,\prod_{j=1}^{\xi^a}\,(t^a_i\]-t^a_j\]+1)\;
\prod_{j=1}^{\xi^{a+1}\!\!}\;(t^a_i\]-t^{a+1}_j\}-1)\,,\kern-2.2em
\\
\cnnu-.6>
\ungood
\cnnu.6>
\cnn-.2>
\endalign
$$
$\an-1$, \>$i=1\lc\xi^a\}$. Here we set $\xi^0\}=\)\xi^N\!\]=\)0$.
\vsk.2>
Given the following data: integers
${a_1\lc a_{|\xi|+k-1}\),\)b_1\lc b_{|\xi|+k-1}\),\)c\)\in\]\lb\)1\lc N\)\rb}$
and ${i\in\]\lb\)1\lc\xi^c\)\rb}$,
a sequence $s_1\lc s_{|\xi|+\)k-1}$, which is a \perm/ of the sequence
$u\lc{u-k+1}\),\alb\)t^1_1\lc\Hat{t^{\)c}_i}\lc t^{\)N-1}_{\xi^{N-1}}$,
where the \var/ $t^{\)c}_i$ is omitted, and a \pol/ $p\)(u\);\txn)$,
\vvn-.2>
consider the product
$$
\align
& \;T_{a_1\],\)b_1}(s_1)\>\ldots\>
T_{a_{|\xi|+k-1}\],\)b_{|\xi|+k-1}}(s_{|\xi|+\)k-1})
\,\,\Ixe^{\)c,i}(t)\,\x{}
\Tagg{Ixik}
\\
\nn7>
& {}\x\,p\)(u\);\txn)\;\prod_{a=1}^{N-1}\,\prod_{j=1}^{\xi^a\!}\,
\biggl({1\over u-t^a_j}
\prod_{1\le j<l\le\xi^a}\]{1\over t^a_j\]-t^a_l}\biggr)
\prod_{a=1}^{N-2}\,\prod_{j=1}^{\xi^a\!}\;
\prod_{l=1}^{\xi^{a+1}\!\!}\;{1\over(t^{a+1}_l\}-t^a_j)^2}\kern-2em
\\
\cnn-.1>
\endalign
$$
where the factors $(u-t^a_j)\1\}$ in \(Ixik) have to be expanded as power
series in $u\1\!$. We denote by $\Ixik\]$ the span over $\C$ of all products
\(Ixik) with various $a_1\lc a_{|\xi|+\)k-1}\]$, $b_1\lc b_{|\xi|+\)k-1}\]$,
$c\>,i$, $s_1\lc s_{|\xi|+\)k-1}$, and $p\)(u\);\txn)$.
\vsk.1>
Set
\vvn-.7>
$$
\Xxi{\>a}(u\);t)\,=\,Q_a\,T_{aa}(u)\,
\prod_{i=1}^{\xi^{a-1}\!\!}\;{u-t^{a-1}_i\}+1\over u-t^{a-1}_i}
\,\prod_{i=1}^{\xi^a}\,{u-t^a_i\]-1\over u-t^a_i}\;,
\Tag{Xxi}
$$
$\an$. Expressions \(Xxi) are regarded as power series in $u\1\!$ \wcoeff/
$\Yn\bigl[\)\txn\bigr]$.
\Th{main}
Let $Q$ be a diagonal matrix. For any $\kN\}$, we have
\vvn-.2>
$$
\Te_{k,Q}(u)\,\BB_{\)\xi}(t)\,\simeq\,
\BB_{\)\xi}(t)\;\sum_{\ab}\,\prod_{r=1}^k\,\Xxi{\>a_r}(u-r+1\);t)\,+\,
\Uxi(u\);t)\,,\kern-1.8em
\vv-.1>
\Tag{TeB}
$$
where the sum is taken over all \$k\)$-tuples $\ab=(\ak)$ \st/
${1\le a_1\lsym<a_k\le N}\]$, and \>$\Uxi(u\);t)$ belongs to \>$\Ixik$,
\vvn.06>
cf.~\(Ixik)\). By Lemma~\[Taa0], the order of the product of elements
\)$\Xxi{\>a_r}$ \}is irrelevant.
\endpro
\nt
The theorem is proved in Section~\[:proof].
\vsk-.5>
\vsk0>
\Rem
Similarly to formula \(DgNT)\), the expressions
${\sum_{\ab}\,\prod_{r=1}^k\,\Xxi{\>a_r}(u-r+1\);t)}$ can be obtained
as coefficients of a suitable formal \difl/ operator:
$$
\align
\Mg_{\>\xi\],Q}(u\),\]\der_u)\, &{}=\,
\bigl(1-\Xxi{\)1}(u\);t)\,e^{-\der_u}\bigr)\)\ldots\)
\bigl(1-\Xxi N(u\);t)\,e^{-\der_u}\bigr)\,={}
\Tag{1-X}
\\
\nn5>
& {}=\,\sum_{k=0}^N\,(-1)^k\!\!\sum_{1\le a_1\lsym<a_k\le N\!\!\!}\;\;
\prod_{r=1}^k\,\Xxi{\>a_r}(u-r+1\);t)\,e^{-k\der_u}\>.
\endalign
$$
\enddemo
\Rem
Since \coeffs/ $\qdet T(u)$ are central in $\Yn$,
the first formula of \(qdet) for the identity \perm/ $\si$ implies that
\vvn-.3>
$$
\qdet T(u)\;\BB_{\)\xi}(t)\,\simeq\,
\BB_{\)\xi}(t)\,\prod_{a=1}^N\,T_{aa}(u-a+1)\,.
\vv-.2>
$$
Together with the relation ${\Te_{N\},Q}(u)\)=\)\det\)Q\,\qdet T(u)}$,
this proves Theorem~\[main] for $k=N\}$.
\enddemo

\Sect[emods]{Tensor products of \emod/s}
In this section we consider the action of the \tmcs/
${\Te_{1,Q}(u)\lc\Te_{N\},Q}(u)}$ in a tensor product of \emod/s over
the Yangian $\Yn$. We will apply Theorem~\[main] to get \egv/s and \eva/s of
the \tmcs/. In Section~\[:Sform] we will also show that the operators of
the action are \sym/ \wrt/ a certain \sym/ bilinear form, see Theorem~\[Shapz].
\vsk.1>
{\bls1.1\bls
Let $\Mn$ be \gnmod/s. Consider the tensor product ${\Moxz}$ of \eval/
\Ynmod/s. The operators $T_{ab}(u)\vst{\Moxz}\]$, \,$a\),\bn\}$, are
\$\End(\Mox)\)$-valued \raf/s in $u,\zn$ with the denominator
${\prod_{i=1}^n\,(u-z_i)}$, \,and
\>$T_{ab}(u)\vst{\Moxz}\}=\>\dl_{ab}\]+\)O(u\1)$ as ${u\to\8}$.
\vsk.1>}
Further on, we will abbreviate, $z=(\zn)$. For any $\kn$, the operator
\vvn.1>
$$
\Te^{\Moxs}_{k,Q}(u\);z)\,=\,\Te_{k,Q}(u)\vst{\Moxz}
\vv-.4>
\Tag{Tez}
$$
is a \raf/ in $u,\zn$ with the denominator
$\prod_{i=1}^n\,\prod_{j=0}^{k-1}\,(u-j-z_i)$, and
\>$\Te^{\Moxs}_{k,Q}(u\);z)\)=\)\tr_{\Vws k}\:\]\Qw k\>\]+\)O(u\1)$
\vvn.2>
\>as ${u\to\8}$. The operators
$\Te^{\Moxs}_{1,Q}(u\);z)\lc\Te^{\Moxs}_{N\},Q}(u\);z)$ are called
\vvn.16>
the \em{\tmcs/ of the \XXX/\)-type model}
associated with the Lie algebra $\gln$.
\vsk.3>
Let $\Mn$ be \hw/ \gnmod/s with \hw/s $\Lan$, $\La_i\]=(\La^1_i\lc\La^N_i)$,
and \hwv/s $v_1\lc v_n$. Then for any $\zn$, the vector ${\vox}$ is a \wt/
singular vector \wrt/ the action of $\Yn$ in the module ${\Moxz}$ and
\vvn-.6>
$$
T_{aa}(u)\>\vox\,=\,\vox\,\prod_{i=1}^n\,{u-z_i\]+\La^a_i\over u-z_i}\;,
\qquad\Rlap{\an\).}
$$
\vsk-.3>
Let ${\xi=(\xin)}$ be a collection of nonnegative integers. The \fn/
$$
\BB^{\vox}_{\)\xi}(t\);z)\,=\,\BB_{\)\xi}(t)\>\vox\;
\prod_{a=1}^{N-1}\>\prod_{i=1}^{\xi^a\!}\,\)\prod_{j=1}^n\,(t^a_i\]-z_j)\,
\prod_{a=1}^{N-2}\,\prod_{i=1}^{\xi^a\!}\;
\prod_{j=1}^{\xi^{a+1}\!\!}\;(t^{a+1}_j\}-t^a_i)\kern-.06em
\vv-.1>
\Tag{BBtz}
$$
is a \pol/ in $\txn$, $\zn$, cf.~\(TTtt) and Lemma~\[thtBB].
We call it the \em{universal \wt/ \fn/} for the \em{\XXX/\]-\)type} model.
\par
Let $Q\)=\]\sum_{a=1}^N\]Q_aE_{aa}$ be a diagonal matrix.
The system of algebraic equations
\vvn-.5>
$$
\align
& Q_a\,\prod_{j=1}^n\>(t^a_i\]-z_j+\La^a_j)\,
\prod_{j=1}^{\xi^{a-1}\!\!}\;(t^a_i\]-t^{a-1}_j\}+1)
\,\prod_{\tsize{j=1\atop j\ne i}}^{\xi^a}\,(t^a_i\]-t^a_j\]-1)\,
\prod_{j=1}^{\xi^{a+1}\!\!}\;(t^a_i\]-t^{a+1}_j)\,={}
\Tagg{Bae}
\\
\nn->
& {}\!=\,Q_{a+1}\,\prod_{j=1}^n\>(t^a_i\]-z_j+\La^{a+1}_j)\,
\prod_{j=1}^{\xi^{a-1}\!\!}\;(t^a_i\]-t^{a-1}_j)
\,\prod_{\tsize{j=1\atop j\ne i}}^{\xi^a}\,(t^a_i\]-t^a_j\]+1)\;
\prod_{j=1}^{\xi^{a+1}\!\!}\;(t^a_i\]-t^{a+1}_j\}-1)\,,\kern-2.5em
\\
\cnn-.4>
\cnnu-.5>
\ungood
\cnnu.5>
\endalign
$$
$\an-1$, \>$i=1\lc\xi^a\}$, \,$\xi^0\}=\)\xi^N\!\]=\)0$,
\vvn.16>
is called the \em{\Bae/s}. Say that a \sol/ $\tti=(\ttxn)$ of system \(Bae)
is \em{\offdiag/} if ${\tti^a_i\]\ne\)\tti^a_j}$ for any $\an-1$,
\>${1\le i<j\le\xi^a}\}$, and ${\tti^a_i\]\ne\)\tti^{a+1}_j}$ for any
${\an-2}$, \>$i=1\lc\xi^a\!$, $j=1\lc\xi^{a+1}\!$. For any $\an\}$, set
\vvn-.7>
$$
\Xxi{\>a}(u\);t\);z\);\La)\,=\,Q_a\,\prod_{i=1}^n\,{u-z_i+\La^a_i\over u-z_i}
\ \prod_{i=1}^{\xi^{a-1}\!\!}\;{u-t^{a-1}_i\}+1\over u-t^{a-1}_i}
\,\prod_{i=1}^{\xi^a}\,{u-t^a_i\]-1\over u-t^a_i}\;.
\Tag{Xxiz}
$$
\Th{main2}
Let $\Mn$ be \hw/ \gnmod/s with \hw/s $\Lan$ and \hwv/s $v_1\lc v_n$.
Let $Q$ be a diagonal matrix and $\tti=(\ttxn)$ an \offdiag/ \sol/
of system \(Bae)\). Then for any $\kN\}$, we have
\vvn-.6>
$$
\Te^{\Moxs}_{k,Q}(u\);z)\,\BB^{\vox}_{\)\xi}(\tti\);z)\,=\,
\BB^{\vox}_{\)\xi}(\tti\);z)\;\sum_{\ab}\,\prod_{r=1}^k\,
\Xxi{\>a_r}(u-r+1\);\tti\);z\);\La)\,,\kern-1em
\Tag{TeBox}
$$
where the sum is taken over all \$k\)$-tuples $\ab=(\ak)$ \st/
${1\le a_1\lsym<a_k\le N}\]$.
\endpro
\Pf.
Theorem~\[main] and formula \(Ixik) imply that
\vvn.4>
$$
\align
& \Te^{\Moxs}_{k,Q}(u\);z)\,\BB^{\vox}_{\)\xi}(t\);z)\,={}
\\
\nn4>
&\!{}=\BB^{\vox}_{\)\xi}(t\);z)\;\sum_{\ab}\,\prod_{r=1}^k\,
\Xxi{\>a_r}(u-r+1\);t\);z\);\La)\,+\,\Uxe^{\)\vox}(u\);t\);z\);\La)\,,
\\
\cnn-.2>
\endalign
$$
where $\Uxe^{\)\vox}(u\);t\);z\);\La)$ is a ratio of a \pol/ \wcoeff/
${\Mox}$ vanishing at \sol/s of system \(Bae) and the product
\vvn-.1>
$$
\prod_{a=1}^{N-1}\prod_{1\le i<j\le\xi^a}\!(t^a_i\]-t^a_j)\,
\prod_{a=1}^{N-2}\,\prod_{i=1}^{\xi^a\!}\;
\prod_{j=1}^{\xi^{a+1}\!\!}\;(t^{a+1}_j\}-t^a_i)
\vv-.1>
$$
that does not vanish at \offdiag/ points $(\ttxn)$.
This yields formula \(TeBox)\).
\epf
For an \offdiag/ \sol/ $\tti=(\ttxn)$ of system \(Bae)\),
the vector $\BB^{\vox}_{\)\xi}(\tti\);z)$ is called the \em{\Bv/}.
\Rem
Similarly to formula \(1-X)\), the \eva/s
${\sum_{\ab}\,\prod_{r=1}^k\,\Xxi{\>a_r}(u-r+1\);t\);z\);\La)}$
\vvn-.2>
of the \tmcs/ $\Te^{\Moxs}_{k,Q}(u\);z)$ can be obtained as coefficients
of a suitable \dif/ operator:
$$
\align
\Mg_{\>\xi\],Q}(u\),\]\der_u\);t\);z\);\La)\,=\,
\bigl(1-\Xxi{\)1}(u\);t\);z\);\La)\,e^{-\der_u}\bigr)\)\ldots\)
\bigl(1-\Xxi N(u\);t\);z\);\La)\,e^{-\der_u}\bigr)\,={}\}\kern-2.4em &
\Tagg{1-Xz}
\\
\nn5>
{}=\,\sum_{k=0}^N\,(-1)^k\!\!\sum_{1\le a_1\lsym<a_k\le N\!\!\!}\;\;
\prod_{r=1}^k\,\Xxi{\>a_r}(u-r+1\);t\);z\);\La)\,e^{-k\der_u}\>.\kern-2.4em &
\endalign
$$
\enddemo
\Rem
Since \coeffs/ $\qdet T(u)$ are central in $\Yn$,
the first formula of \(qdet) for the identity \perm/ $\si$ implies that
$$
\qdet T(u)\vst{\Moxz}=\,
\prod_{i=1}^n\,\prod_{a=1}^N\,{u-a+1-z_i+\La^a_i\over u-a+1-z_i}\;.
\vv-.1>
$$
Together with relation \,${\Te_{N\},Q}(u)\)=\)\det\)Q\,\qdet T(u)\)
\prod_{i=1}^n\>\prod_{a=1}^N\>(u-a+1-z_i)}$, \,this proves
\vvnu-.4>\uline
Theorem~\[main2] for $k=N\}$.
\enddemo
\Prop{Bsing}
\back\;\cite{KR1}
Let $Q$ be the identity matrix, that is, ${Q_1\]\lsym=Q_N\]=1}$,
and $\tti=(\ttxn)$ an \offdiag/ \sol/ of the \Bae/s \(Bae)\).
Then the \Bv/ $\BB^{\vox}_{\)\xi}(\tti\);z)$ is a singular vector \wrt/
the \,\$\gln$-\)action in $\Mox$ of \wt/ $\bigl(\)
\sum_{i=1}^n\La^1_i\]-\)\xi^1\],\,\sum_{i=1}^n\La^2_i\]+\)\xi^1\}-\)\xi^2\]
\lc\)\sum_{i=1}^n\La^N_i\!+\)\xi^{N\]-1}\)\bigr)$.
\endpro
\Pf.
Clearly, for any $t$ the vector $\BB^{\vox}_{\)\xi}(t\);z)$ has the indicated
\wt/. The fact that the \Bv/ $\BB^{\vox}_{\)\xi}(\tti\);z)$ is singular
is proved in Section~\[:proof].
\epf
\Rem
Let $\Mn$ be any \gnmod/s, not necessarily \hw/ ones, and let ${v\in\Mox}$
be any \wt/ singular vector \wrt/ the action of $\Yn$ in the module ${\Moxz}$.
Define \pol/s $q^v_1(u)\lc q^v_N(u)$ by the rule
\vvn.1>
$$
T_{aa}(u)\>v\;\prod_{i=1}^n\,(u-z_i)\,=\,q^v_a(u)\>v\,,\qqq\Rlap{\an\).}
\vv-.8>
$$
The \fn/
\vvn-.6>
$$
\BB^v_{\)\xi}(t\);z)\,=\,\BB_{\)\xi}(t\);z)\>v\;
\prod_{a=1}^{N-1}\>\prod_{j=1}^{\xi^a\!}\,\)\prod_{i=1}^n\,(t^a_i\]-z_j)\,
\prod_{a=1}^{N-2}\,\prod_{i=1}^{\xi^a\!}\;
\prod_{j=1}^{\xi^{a+1}\!\!}\;(t^{a+1}_j\}-t^a_i)
$$
is a \pol/ in $\txn$, cf.~\(TTtt) and Lemma~\[thtBB]. For $\an\}$, set
$$
\Xxi{\>a,v}(u\);t)\,=\,Q_a\,q^v_a(u)\,\prod_{i=1}^n\,{1\over u-z_i}
\ \prod_{i=1}^{\xi^{a-1}\!\!}\;{u-t^{a-1}_i\}+1\over u-t^{a-1}_i}
\,\prod_{i=1}^{\xi^a}\,{u-t^a_i\]-1\over u-t^a_i}\;.
$$
Similarly to Theorem~\[main2], we can show that for any \offdiag/ \sol/
$\tti$ of the system
$$
\align
Q_a & \,q^v_a(t^a_i)\,
\prod_{j=1}^{\xi^{a-1}\!\!}\;(t^a_i\]-t^{a-1}_j\}+1)
\,\prod_{\tsize{j=1\atop j\ne i}}^{\xi^a}\,(t^a_i\]-t^a_j\]-1)\,
\prod_{j=1}^{\xi^{a+1}\!\!}\;(t^a_i\]-t^{a+1}_j)\,={}
\\
\nn->
& {}\!=\,Q_{a+1}\,q^v_{a+1}(t^a_i)\,
\prod_{j=1}^{\xi^{a-1}\!\!}\;(t^a_i\]-t^{a-1}_j)
\,\prod_{\tsize{j=1\atop j\ne i}}^{\xi^a}\,(t^a_i\]-t^a_j\]+1)\;
\prod_{j=1}^{\xi^{a+1}\!\!}\;(t^a_i\]-t^{a+1}_j\}-1)\,,\kern-1.6em
\\
\cnnu-.5>
\ungood
\cnnu.5>
\endalign
$$
$\an-1$, \>$i=1\lc\xi^a\}$, \)where $\xi^0\}=\)\xi^N\!\]=\)0$, we have
\vvn-.2>
$$
\Te^{\Moxs}_{k,Q}(u\);z)\,\BB^v_{\)\xi}(\tti\);z)\,=\,
\BB^v_{\)\xi}(\tti\);z)\;\sum_{\ab}\,\prod_{r=1}^k\,
\Xxi{\>a_r,v}(u-r+1\);\tti\))\,.
\vv-.4>
$$
Moreover, for any $t$,
\vvn-.5>
$$
\qdet T(u)\;\BB^v_{\)\xi}(t\);z)\,=\,\BB^v_{\)\xi}(t\);z)\;
\prod_{j=0}^{N-1}\,\biggl(q^v_a(u-j)\,
\prod_{i=1}^n\,{1\over u-j-z_i}\biggr)\,.
$$
\enddemo

\Sect[Sform]{Shapovalov form}
Define an \ainv/ ${\tau\]:\Un\to\Un}$ by the rule
${\tau(e_{ab})=e_{ba}}$, $a\),\bn$. Let $M$ be a \hw/ \gnmod/ with \hwv/ $v$.
The \em{Shapovalov form} ${S_M:M\ox M\to\C}$ is the unique symmetric bilinear
form \st/ ${S_M(v\),\]v)=1}$, and
${S_M(Xw_1\),\]w_2)\)=\)S_M\bigl(w_1\),\tau(X)\)w_2\bigr)}$
for any ${X\]\in\Un}$ and any ${w_1,w_2\in M}$.
For an \irr/ module $M$, the form $S_M\}$ is \ndeg/.
\vsk.1>
\Lm{opitau}
For any $X\]\in\Yn$, ${w_1,w_2\in M}$, and $x\in\C$,
we have
$$
S_M(X\vst{M(x)}\]w_1\),\]w_2)\,=\,
S_M\bigl(w_1\),\opi(X)\vst{M(x)}\]w_2\bigr)\,.
\vv-.2>
$$
\endpro
\nt
The proof is straightforward.
\vsk.5>
Let $\Mn$ be \hw/ \gnmod/s. Then formula~\(Dlopi) and Lemma~\[opitau] imply
that for any ${X\]\in\Yn}$ and any ${w_1,w_2\in\Mox}$, we have
$$
\align
& \SMox\bigl(X\vst{\Moxz}\]w_1\),\]w_2\bigr)\,=\,
\Tag{opitauz}
\\
\nn4>
&\){}=\,\SMox\bigll( w_1\),\bigl(\opi(X)\vst{\Mozx}\bigr)\2{n\ldots1}w_2\bigrr)\,,
\\
\cnn-.1>
\endalign
$$
where $S_{M_i}\!$ is the Shapovalov form for the module $M_i$, $\inn$.
\Prop{tauTe}
For any $\kN\}$, and any ${w_1,w_2\in\Mox}$, we have
\vvn.1>
$$
\align
& \SMox\bigl(\Te^{\Moxs}_{k,Q}(u\);\zn)\>w_1\),\]w_2\bigr)\,=\,
\\
\nn6>
&\){}=\,\SMox\bigll( w_1\),\bigl(\)
\Te^{\Moxxs}_{k,Q^\dag}(u\);\zni)\bigr)\2{n\ldots1}w_2\bigrr)\,,\kern-2em
\\
\cnn-.1>
\endalign
$$
where the superscript $\)\dag\]$ stands for the transposition of matrices.
\endpro
\Pf.
The claim follows from formulae~\(Tez)\), \(opitauz) and Proposition~\[Tesym].
\epf
We remind below some facts on \Rms/ for tensor products of \emod/s over
the Yangian $\Yn$. Propositions~\[RLM]\,--\,\[RYB] follow from Drinfeld's
results on the universal \Rms/ for Yangians \cite{D}\).
\Prop{RLM}
Let $L$ and $M$ be \hw/ \gnmod/s with the respective \uline
\hwv/s $v$ and $w$.
For generic $x,y\in\C$ there exists a unique operator
${R_{LM}\:(x-y)\in\End(L\ox\}M)}$ \st/
\vvn-.4>
$$
\gather
R_{LM}\:(x-y)\;X\vst{L(x)\)\ox M(y)}=\,
\bigl(X\vst{M(y)\)\ox L(x)}\bigr)\2{21}\)R_{LM}\:(x-y)
\Tag{RX}
\\
\nn5>
\Text{for any $X\]\in\Yn$, and}
\nn4>
R_{LM}\:(u)\,v\ox w\)=\)v\ox w\,.
\Tag{Rvw}
\endgather
$$
\endpro
\nt
The operator $R_{LM}\:(u)$ is called the \em{\rat/ \Rm/} for the tensor
product $L\ox\}M$.
\Ex
Let \)$L=\Vw l\}$ and \)$M=\Vw m\}$. \)Then
\vvn-.3>
$$
R_{LM}\:(u)\)=\)\Rw{\)l,\)\wg m}(u)\;{u+\max\)(m-l,0)\over u+m}\,\,
\prod_{i=0}^{l-1}\,\prod_{j=0}^{m-1}\>{1\over u+j-i}\;.
\vv-.1>
$$
\enddemo
\Prop{RLMinv}
$R_{LM}\:(u)\>=\>\bigll(\]\bigl(R_{ML}\:(-\)u)\bigr)\2{21}\bigrr)^{\]-1}\}$.
\endpro
\Prop{RYB}
Let $M_1, M_2, M_3$ be \hw/ \$\gln$-modules. Then the operators
$R_{M_iM_j}\:\}(u)$, \)$1\le i<j\le 3$, \)satisfy the \YB/
\vvn.3>
$$
R\"{12}_{M_1M_2}\}(u_1\]-u_2)\>R\"{13}_{M_1M_3}\}(u_1)\>
R\"{23}_{M_2M_3}\}(u_2)\>=\>R\"{23}_{M_2M_3}\}(u_2)\>
R\"{13}_{M_1M_3}\}(u_1)\>R\"{12}_{M_1M_2}\}(u_1\]-u_2)\,.
$$
\endpro
\vsk.4>
By the definition of \emod/s over $\Yn$, formula \(RX) means that the operator
$R_{LM}\:(u)$ commutes with the \$\gln$-action in ${L\ox\}M}$ and for any
$a\),\bn\}$,
$$
R_{LM}\:(u)\>\Bigl(\)u\>e_{ab}\ox 1\,+\tsum_{c=1}^N\>e_{ac}\ox e_{cb}\)\Bigr)
\,=\,\Bigl(\)u\>e_{ab}\ox 1\,+\tsum_{c=1}^N\>e_{cb}\ox e_{ac}\)\Bigr)\>
R_{LM}\:(u)\,.
\vv-.1>
$$
The operator $R_{LM}\:(u)$ is uniquely determined by these properties together
with the normalization condition \(Rvw)\), which implies that the restriction
of $R_{LM}\:(u)$ to any \wt/ subspace of the \$\gln$-module ${L\ox\}M}$ is
a \raf/ of $u$, and
\vvn-.1>
$$
R_{LM}\:(u)\,=\,1+O(u\1)\,,\qqq\Rlap{u\to\8\,.}
\Tag{R8}
$$
In addition, for any ${w_1,w_2\in L\ox\}M}$, we have
\vvn.1>
$$
S_L\}\ox S_M\bigl(R_{LM}\:(u)\>w_1\),\]w_2\bigr)\,=\,
S_L\}\ox S_M\bigl(w_1\),R_{LM}\:(u)\>w_2\bigr)\,.
\vv.1>
\Tag{tauR}
$$
\par
Let $\Mn$ be \hw/ \gnmod/s. Set
$$
R_{\Mnd}\:\}(z)\,=\prodr_{1\le i<j\le k}\!R_{M_iM_j}\"{ij}\}(z_i\]-z_j)\>
\in\)\End(\Mox)\,.\kern-2em
\vv-.4>
\uugood
\Tag{RMox}
$$
By Proposition~\[RLMinv] we have
\vvn-.3>
$$
\bigl(R_{\Mnd}\:\}(z)\bigr)^{\]-1}\)=
\prodl_{1\le i<j\le k}\!R_{M_jM_i}\"{ji}\}(z_j\]-z_i)\,.
$$
For any ${w_1,w_2\in\Mox}$, formula~\(tauR) and Proposition~\[RYB] yield
\vvn.1>
$$
\SMox\bigl(R_{\Mnd}\:\}(z)\>w_1\),\]w_2\bigr)\,=\,
\SMox\bigl(w_1\),R_{\Mnd}\:\}(z)\>w_2\bigr)\,,
\Tag{tauRM}
$$
and for any $X\]\in\Yn$, formula~\(RX) implies
$$
R_{\Mnd}\:\}(z)\;X\vst{\Moxz}=\,
\bigl(X\vst{\Mozx}\bigr)\2{n\ldots1}\)R_{\Mnd}\:\}(z)\,.
\Tag{RXXR}
$$
\vsk.2>
Define a bilinear form $S_{\Moxs}^{\>z}\!$ on $\Mox\}$ by the rule
\vvn.2>
$$
S_{\Moxs}^{\>z}\}(w_1\),w_2)\,=\,
(\SMox)\bigl(w_1\),\)R_{\Mnd}(z)\>w_2\bigr)\,.
\vv.2>
\Tag{SMz}
$$
For \irr/ modules $\Mn$ and generic $z$, the form $S_{\Moxs}^{\>z}\!$ is
\ndeg/. By equality~\(tauRM)\), the form $S_{\Moxs}^{\>z}\!$ is \sym/.
\Th{Shapz}
For a \sym/ matrix $Q\}$, ${Q\]=Q^\dag}$, the operators
\vvn.06>
\>$\Te^{\Moxs}_{k,Q}(u\);z)$, ${\kN}\}$, are \sym/ \wrt/
the \sym/ bilinear form $S_{\Moxs}^{\>z}\}$.
\endpro
\Pf.
Proposition~\[tauTe] and formula~\(RXXR) yield
\vvn.2>
$$
S_{\Moxs}^{\>z}\}\bigl(\)\Te^{\Moxs}_{k,Q}(u\);z)\>w_1\),w_2\bigr)\,=\,
S_{\Moxs}^{\>z}\}\bigl(w_1\),\)\Te^{\Moxs}_{k,Q^\dag}(u\);z)\>w_2\bigr)
\vv.1>
$$
for any ${\kN}\}$. This proves the statement.
\epf

\Sect[current]{Current algebra $\glnx$}
Let $\glnx$ be the Lie algebra of \pol/s \wcoeff/ $\gln\}$ with the pointwise
commutator. We call it the \em{current algebra}. In this section we construct,
following~\cite{T}\), analogues of the \tmcs/ $\Te_{1,Q}(u)\lc\Te_{N\},Q}(u)$
for the algebra $\Unx$, and formulate the counterparts of Theorems~\[main] and
\[main2].
\vsk.1>
We will be using the standard generators of $\glnx$,
$e_{ab}\3s=\)e_{ab}\>x^{\)s}\}$, $a\),\bn\}$, and $s=0,1,\ldots$,
which commute as follows: ${[\)e_{ab}\3r\),\)e_{cd}\3s\)]\)=\)
\dl_{bc}\>e_{ad}\3{r+s}\}-\)\dl_{ad}\>e_{cb}\3{r+s}}\}$. We will identify
the Lie algebra $\gln\}$ with the subalgebra in $\glnx\}$ of constant \pol/s,
thus making no \dif/ between the elements $e_{ab}\30\}$ and the generators
$e_{ab}\}$ of $\gln$.
\vsk.2>
Let $\ngp\}=\!\!\Plus_{\abn}\!\!\C\)e_{ab}\]$ and
\)$\ngm\}=\!\!\Plus_{\abn}\!\!\C\)e_{ba}\}$ be the standard nilpotent
\vv.06>
subalgebras in $\gln$, and let $\ngpmx$ be the corresponding subalgebras
\vvn.1>
in $\glnx$. If $A\>,B\in\Unx$ and $A-\]B\in\Unx\>\ngpx$, then we will write
$A\)\simeq B$.
\vsk.2>
A vector $v$ in a \gnxmod/ is called \em{singular} \wrt/ the action of $\glnx$
if ${\ngpx\>v\)=\)0}$. A singular vector $v$ that is an \egv/ for the action
of $e_{aa}\3s$, $\an\}$, $s=0,1,\ldots$, is called a \em{\wt/ singular vector}.
\vsk.2>
Let $M$ be a \gnmod/. For any $y\in\C$ denote by $M(y)$ the corresponding
\eval/ \gnxmod/, that is, $e_{ab}\3s\]\vst{M(y)}\}=y^s\)e_{ab}\vst M$,
$a\),\bn\}$, \)$s=0,1,\ldots$.
\vsk.1>
We organize elements $e_{ab}\3s\}$, $a\),\bn\}$, \>$s=0,1,\ldots$,
into generating series
\vvn-.3>
$$
L_{ab}(u)\,=\,\sum_{s=0}^\8\,e_{ba}\#s\)u^{-s-1}\),
\qqq\Rlap{a\),\bn\,.}
\vv-.8>
\Tag{Lab}
$$
and combine those series together into a series
\vvn-.16>
${L(u)\>=\!\sum_{a,b=1}^N\!E_{ab}\ox L_{ab}(u)}$ \wcoeff/
${\End(\Cn)\ox\alb\Unx}$. Notice the flip of subscripts of the generators
in the sum~\(Lab)\); this is done to have better notational correspondence
further with the Yangian case. The commutation relations in $\Unx$ can be
written in the following form:
\vvn.2>
$$
(u-v)\>\bigl[\)L_{ab}(u)\>,L_{cd}(v)\)\bigr]\,=\,
\dl_{bc}\)\bigl(L_{ad}(u)-L_{ad}(v)\bigr)\)-\)
\dl_{ad}\)\bigl(L_{cb}(u)-L_{cb}(v)\bigr)\,,
\vv.2>
$$
that amounts to the equality
\vvn.2>
for series \wcoeff/ \)${\End(\Cn)\ox\End(\Cn)\ox\Unx}$:\ifUS\kern-14pt\fi
\vvn.3>
$$
(u-v)\>\bigl[\)L\"{13}(u)\>,L\"{23}(v)\)\bigr]\,=\,
\bigl[\)L\"{13}(u)+L\"{23}(v)\>,P\"{12}\)\bigr]\,.
\vv.3>
\Tag{LL}
$$
Here ${P\]\in\End(\Cn)\ox\End(\Cn)}$ is the flip map.
\vsk.2>
Fix a matrix $K\]\in\End(V)$ and consider the formal \difl/ operator
$$
\align
& \De_{\]K}\](u\),\der_u)\,={}
\Tagg{detL}
\\
\nn4>
{}=\,{} & (\trVN\ox\id\))\Bigl(\]
\bigl(\der_u\]-K\"1\}-L\"{1,\)N+1}(u)\bigr)\)\ldots\)
\bigl(\der_u\]-K\"N\}-L\"{N\],\)N+1}(u)\bigr)\;
\AA\"{1\ldots N}\6N\)\Bigr)\,.\kern-1em
\endalign
$$
of order $N\}$. Let $\Gc_{0,K}(u)\)\lc\Gc_{N\},K}(u)\in\)\Unx\lbc\)u\1\rbc$
be its coefficients:
$$
\De_{\]K}\](u\),\der_u)\,=\,
\sum_{k=0}^N\,(-1)^k\,\Gc_{\)k,K}(u)\,\der_u^{\>N\]-k}\).
\Tag{DG}
$$
The series $\Gc_{1,K}(u)\)\lc\Gc_{N\},K}(u)$ are called the \em{Gaudin \tmcs/}.
\vsk-.5>\vsk0>
\Ex
Clearly, ${\Gc_{0,K}(u)\)=\)1}$ because
\>${\tr_{V^{\]\ox N\]}}\:\AA\6N=\)1}$.
\vvn-.5>
Let ${K=\!\sum_{a,b=1}^N\!K_{ab}\)E_{ab}}\}$.
\uline We also have that
\vvn-.6>
$$
\gather
\Gc_{1,K}(u)\,=\,({\tr_V\:\}}\ox\id\))\bigl(K\"1\}+L\"{12}(u)\bigr)\,=\,
\sum_{a=1}^N\,\bigl(K_{aa}\]+L_{aa}(u)\bigr)\,,
\Tag{Gc12}
\\
\nn6>
{\align
\Gc_{2,K}(u)\, &{}=\,({\tr_{V^{\]\ox2\]}}\:\}}\ox\id\))
\Bigl(\]\bigll(\]\bigl(K\"1\}+L\"{13}(u)\bigr)\>\bigl(K\"2\}+L\"{23}(u)\bigr)-
\dot L\"{23}(u)\bigrr)\,\AA\"{12}\62\)\Bigr)\,={}
\\
\nn4>
& {}=\,{1\over2}\,\Bigl(\]\bigl(\)\Gc_{1,K}(u)\bigr)^2\]-\)
(N-1)\>\dot\Gc_{1,K}(u)\>-\sum_{a,b=1}^N
\bigl(K_{ab}\]+L_{ab}(u)\bigr)\bigl(K_{ba}\}+L_{ba}(u)\bigr)\]\Bigr)\,,
\endalign}
\\
\cnn-.4>
\endgather
$$
where \>$\dot G\)=\)d\)G/d\)u\>$.
\vvu-.8>
\uugood
\vvu.8>
\enddemo
The next two propositions were obtained in \cite{T} from the Yangian case.
We will give their proofs at the end of Section~\[:filt].
\Prop{ADLA}
We have
$$
\AA\"{1\ldots N}\6N\,\De\"{N\]+1}_{\]K}\](u\),\der_u)\,=\,
\bigl(\der_u\]-K\"1\}-L\"{1,\)N+1}(u)\bigr)\)\ldots\)
\bigl(\der_u\]-K\"N\}-L\"{N\],\)N+1}(u)\bigr)\;\AA\"{1\ldots N}\6N\).\!
\vv-.5>
$$
\endpro
\Prop{Gckl}
We have \,$\bigl[\)\Gc_{\)k,K}(u)\>,\Gc_{\)l,K}(v)\)\bigr]\)=\)0$.
\endpro
\Prop{K=0}
Let $K$ be the zero matrix. Then for any $\kN$, the coefficients of
\>$\Gc_{\)k,K}(u)$ commute with the subalgebra $\Un\sub\Unx$.
\endpro
\Pf.
Similarly to \(EeT)\), for any $a\),\bn\}$, we have
$\bigl[\)E_{ab}\ox\one+\one\ox e_{ab}\>,L(u)\)\bigr]\)=\)0$.
Hence, the claim follows from formulae \(detL)\), \(DG)\),
and the cyclic property of the trace, because
$\bigl[\)E\"1_{ab}\}\lsym+E\"N_{ab}\],\)\AA\"{1\ldots N}\6N\)\bigr]\)=\)0$.
\epf

We extend the \ainv/ ${\tau\]:\Un\to\Un}$ to the algebra $\Unx$ by the rule
$$
\tau(e_{ab}\3s)\>=\>e_{ba}\3s\,,\qqq\Rlap{a\),\bn\),\quad s=0,1,\ldots\,.}
\Tag{tau}
$$
The extension is consistent with the standard embedding $\Un\hto\Unx$.
\Prop{Gcsym}
We have ${\tau\bigl(\)\Gc_{k,K}(u)\bigr)\>=\>\Gc_{k,K^\dag}(u)}$,
where the superscript $\)\dag\]$ stands for the transposition of matrices.
\endpro
\nt
The statement is proved in Section~\[:filt].
\vsk.5>
Let ${\xi=(\xin)}$ be a collection of nonnegative integers.
Set ${|\)\xi\)|=\xi^1\]\lsym+\xi^{N-1}}\}$.
\vvn.16>
Consider an expression in $|\)\xi\)|$ \var/s $\txin$ \wcoeff/ $\Unmx$:
\vvn-.3>
$$
\alignat2
& \FF_{\]\xi}(\txn)\,={}
\Tagg{FF}
\\
\nn4>
& \!\}{}=\,\Sym_{\)t^1_1\lc t^1_{\xi^1}}\}\ldots\,
\Sym_{\)t^{\)N-1}_1\}\lc t^{\)N-1}_{\xi^{N-1}}}\biggl\lb
\,\sum_{\mb}\,\prodr_{\abn\!}\biggl(\){1\over m_{ab}\)!}\; && \]
\prod_{i=1}^{m_{ab}}\,L_{ab}(t^a_{\mh^a_{\}ab}\]+\)i})\,\x{}\kern-2em
\\
\nn3>
&& \Llap{{}\x{}\,} & \prod_{c=a}^{b-2}\,{1\over
t^{c+1}_{\mh^{c+1}_{\}ab}+\)i}\}-t^c_{\mh^c_{\}ab}+\)i}}\biggr)\biggr\rb\,,
\kern-3em
\endalignat
$$
where the sum is taken over all arrays of nonnegative integers
${\mb=(m_{ab})_{\)1\le a<b\le N}}$ \st/
$\sum_{c=1}^a\>\sum_{b=a+1}^N\!m_{\)c\)b}=\>\xi^a\}$, and
\>$\mh^c_{\]ab}=\!\!\sum_{\tsize{(r,s)<(a,b)\atop s>c}\!}\!\!m_{\)rs}$.
\vsk-.72>
\vsk0>
\Lm{FFpol}
The product
\vvn-.1>
\,$\dsize{\FF_{\]\xi}(\txn)\,\prod_{a=1}^{N-2}\,\prod_{i=1}^{\xi^a\!}\;
\prod_{j=1}^{\xi^{a+1}\!\!}\;(t^{a+1}_j\}-t^a_i)}$ \)is a power series in
\>${(t^1_1)\vpb{-1}\}\lc(t^{N-1}_{\xi^{N-1}})\vpb{-1}}\!$ \wcoeff/ $\Unmx$.
\endpro
\nt
The proof is straightforward.
\vsku-.8>
\goodbu
\vsku.8>
\par
Let $K$ be a diagonal matrix, $K=\]\sum_{a=1}^N\]K_aE_{aa}$.
Introduce the expressions
\vvn-.3>
$$
\align
\Jxe^{\)a,i}(t)\, &{}=\,K_a\]-K_{a+1}+\)
L_{aa}(t^a_i)\)-\)L_{a+1,a+1}(t^a_i)\,+{}
\Tagg{Je}
\\
\nn5>
& {}\)+\,\sum_{j=1}^{\xi^{a-1}\!\!}\;{1\over t^a_i\]-t^{a-1}_j}\;-\,
2\)\sum_{\tsize{j=1\atop j\ne i}}^{\xi^a}\,{1\over t^a_i\]-t^a_j}\;+\,
\sum_{j=1}^{\xi^{a+1}\!\!}{1\over t^a_i\]-t^{a+1}_j}\;,
\\
\cnn-.2>
\endalign
$$
$\an-1$, \>$i=1\lc\xi^a\}$. Here we set $\xi^0\}=\)\xi^N\!\]=\)0$.
\vsk.2>
Denote \>$\dsize{L\4r_{ab}(u)\)=\){d^r\over d\)u^r}L_{ab}(u)}$\>.
\vvn.2>
Given the following data:
\Item
an integer $m$ not excceding ${|\)\xi\)|+k-1}$,
integers $a_1\lc a_m\),\alb\)b_1\lc b_m\),\){c\)\in\]\lb\)1\lc N\)\rb}$ and
${i\in\]\lb\)1\lc\xi^c\)\rb}$,
\iitem
a sequence $s_1\lc s_m$, where each $s_j$ is either $u$ or one of the \var/s
$t^1_1\lc\Hat{t^{\)c}_i}\lc t^{\)N-1}_{\xi^{N-1}}$, where the \var/ $t^{\)c}_i$
is excluded,
\iitem
nonnegative integers $r_1\lc r_m$ \st/ ${r_j\]=0}$ if ${s_j\ne u}$,
\vvn-.5>
\,and $\sum_{s_j=u}\}r_j\le k$,
\iitem
a \pol/ $p\)(u\);\txn)$,
\vsk.1>\nt
consider the product
$$
\align
L\4{r_1}_{a_1b_1}(s_1)\> &{}\ldots\>L\4{r_m}_{a_mb_m}(s_m)
\,\,\Jxe^{\)c,i}(t)\,\x{}
\Tagg{Jxik}
\\
\nn5>
{}\x{}& \,p\)(u\);\txn)\;
\prod_{a=1}^{N-1}\,\prod_{j=1}^{\xi^a\!}\,{1\over(u-t^a_j)^k}\;
\prod_{a=1}^{N-2}\,\prod_{j=1}^{\xi^a\!}\;
\prod_{l=1}^{\xi^{a+1}\!\!}\;{1\over t^{a+1}_l\}-t^a_j}\kern-2em
\\
\cnn-.2>
\endalign
$$
where the factors $(u-t^a_j)^{-k}\}$ in \(Jxik) have to be expanded as power
series in $u\1\!$. We denote by $\Jxik\]$ the span over $\C$ of all products
\(Jxik) with various $m$, $a_1\lc a_m\]$, $b_1\lc b_m\]$, $r_1\lc r_m\]$,
$c\>,i$, $s_1\lc s_m$, and $p\)(u\);\txn)$.
\vsk.2>
Consider the formal \difl/ operator
\vvn-.3>
$$
\Mxi(u\),\der_u\);t)\,=
\prodr_{1\le a\le N\!}\biggl(\der_u-K_a-\)L_{aa}(u)\)-\>
\sum_{i=1}^{\xi^{a-1}\!\!}\,{1\over u-t^{a-1}_i\}}\;+\>
\sum_{i=1}^{\xi^a}\,{1\over u-t^a_i}\biggr)\,.
$$
Let \>$\Zxi{\)1}(u\);t)\lc\Zxi{\)N\}}(u\);t)$ be its coefficients:
\vvn-.2>
$$
\Mxi(u\),\der_u\);t)\,=\,
\sum_{k=0}^N\,(-1)^k\,\Zxi{\>k}(u\);t)\,\der_u^{\>N\]-k}\).
\vv-.2>
$$
Expressions $\Zxi{\)1}(u\);t)\lc\Zxi{\)N\}}(u\);t)$ are regarded as power
\vvn.1>
series in $u\1\!$ \wcoeff/ $\Unx\bigl[\)\txn\bigr]$.
\Th{main3}
Let $K$ be a diagonal matrix. For any $\kN\}$, we have
\vvn.3>
$$
\Gc_{\)k,K}(u)\,\FF_{\]\xi}(t)\,\simeq\,
\FF_{\]\xi}(t)\,\Zxi{\>k}(u\);t)\>+\>\Wxi(u\);t)\,,
\vvn.2>
\Tag{GcF}
$$
where \>$\Wxi(u\);t)$ belongs to $\Jxik\}$.
\endpro
\nt
The theorem is proved in Section~\[:proof2].
\vsku-.6>
\goodbu
\vsku.6>
\vsku0>
\Rem
It is easy to see that \coeffs/ $\Gc_{1,K}(u)$ are central in $\Unx$ and
$\Zxi{\)1}(u\);t)=\Gc_{1,K}(u)$ for any $\xi$. This implies Theorem~\[main3]
for $k=1$.
\enddemo

\Sect[cemods]{Tensor products of \emod/s of the current algebra}
Let $\Mn$ be \gnmod/s. Consider the tensor product $\Moxz$ of \eval/
\gnxmod/s. Then
\vvn-.6>
$$
L_{ab}(u)\vst{\Moxz}=\,\sum_{i=1}^n\){e_{ba}\"i\over u-z_i}\;,\qqq\Rlap{\an\),}
\vv-.3>
\Tag{Labz}
$$
where we consider
$e_{ba}\"i\]=\)\one^{\ox(i\)-1)}\]\ox e_{ab}\]\ox\)\one^{\ox(n-i)}$ as
an operator acting in $\Mox$, cf.~\(Lab)\). By formulae \(detL)\), \(DG)\),
for any $\kn$, the operator
\vvn.3>
$$
\Gc^{\)\Moxs}_{k,K}(u\);z)\,=\,\Gc_{k,K}(u)\vst{\Moxz}
\vvn-.3>
\Tag{Gcz}
$$
is a \raf/ in $u,\zn$ with the denominator $\prod_{i=1}^n\,(u-z_i)^k\}$ and
\vvn-.1>
$$
\Gc^{\)\Moxs}_{k,K}(u\);z)\,=\,\tr K^{\wg k}+\,O(u\1)\,,
\qqq\Rlap{u\to\8\,.}
\vv.1>
$$
The operators $\Gc^{\Moxs}_{1,K}(u\);z)\lc\Gc^{\Moxs}_{N\},K}(u\);z)$
\vvn.16>
are called the \em{\tmcs/ of the Gaudin model} associated with the Lie algebra
$\gln$.
\vsk.2>
Let $\tau$ be the \ainv/ \(tau) on $\Unx$. Recall that for a \hw/ \gnmod/ $M$
we denote by $S_M$ the Shapovalov form on $M$.
\vsk.1>
Let $\Mn$ be \hw/ \gnmod/s. For any $\zn$, ${X\]\in\Unx}$,
and ${w_1,w_2\in\Mox}$, we have
$$
\align
& \SMox\bigl(X\vst{\Moxz}\]w_1\),\]w_2\bigr)\,={}
\\
\nn5>
& \){}=\,
\SMox\bigl(w_1\),\tau(X)\vst{\Moxz}\]w_2\bigr)\,.
\endalign
$$
Then Proposition~\[Gcsym] and formula \(Gcz) yield that for any $\kN\}$,
\vvn.3>
$$
\SMox\bigl(\>\Gc^{\)\Moxs}_{k,K}(u\);z)\>w_1\),\]w_2\bigr)\,=\,
\SMox\bigl(w_1\),\>\Gc^{\)\Moxs}_{k,K^\dag}(u\);z)\>w_2\bigr)\,.
$$
In particular, we have the following result.
\Th{Shap}
For a \sym/ matrix $K\}$, ${K\]=K^\dag}$, the operators
\>$\Gc^{\)\Moxs}_{k,K}(u\);z)$, ${\kN}\}$, are \sym/ \wrt/
the tensor product ${\SMox}\!\}$ of the Shapovalov forms.
\endpro
\vsk.1>
Let $\Mn$ be \hw/ \gnmod/s with \hw/s $\Lan$, $\La_i\]=(\La^1_i\lc\La^N_i)$,
and \hwv/s $v_1\lc v_n$. Then for any $\zn$, the vector ${\vox}$ is a \wt/
\vvn.06>
singular vector \wrt/ the action of $\glnx$ in the module ${\Moxz}$ and
\vvn-.5>
$$
L_{aa}(u)\>\vox\,=\,\vox\,\sum_{i=1}^n\){\La^a_i\over u-z_i}\;,
\qqq\Rlap{\an\).}
$$
\vsk-.6>
Fix a diagonal matrix $K=\]\sum_{a=1}^N\]K_aE_{aa}$. Let ${\xi=(\xin)}$
be a collection of nonnegative integers. Set $\xi^0\}=\)\xi^N\!\]=\)0$.
The system of equations
\vvn.2>
$$
\sum_{j=1}^n\,{\La^a_j\]-\La^{a+1}_j\over t^a_i\]-z_j}\;+\,
\sum_{j=1}^{\xi^{a-1}\!\!}\;{1\over t^a_i\]-t^{a-1}_j}\;-\,
2\)\sum_{\tsize{j=1\atop j\ne i}}^{\xi^a}\,{1\over t^a_i\]-t^a_j}\;+\,
\sum_{j=1}^{\xi^{a+1}\!\!}{1\over t^a_i\]-t^{a+1}_j}\;=\,K_{a+1}\]-K_a\,,
\kern-1.2em
\Tag{Bae2}
$$
$\an-1$, \>$i=1\lc\xi^a\}$, is called the \em{\Bae/s}. We always assume that
for a \sol/ $\tti=(\ttxn)$ of system \(Bae2) any denominator in \eq/s does not
equal zero unless the corresponding numerator equals zero.
\vsk.1>
Consider the \difl/ operator
\vvn-.3>
$$
\Mxi(u\),\der_u\);t\);z\);\La)\,=
\prodr_{1\le a\le N\!}\biggl(\der_u-K_a-\>\sum_{i=1}^n\){\La^a_i\over u-z_i}
\;-\>\sum_{i=1}^{\xi^{a-1}\!\!}\,{1\over u-t^{a-1}_i\}}\;+\>
\sum_{i=1}^{\xi^a}\,{1\over u-t^a_i}\biggr)\,.\kern-2em
\Tag{MxiL}
$$
Let \>$\Zxi{\)1}(u\);t\);z\);\La)\lc\Zxi{\)N\}}(u\);t\);z\);\La)$
be its coefficients:
$$
\Mxi(u\),\der_u\);t\);z\);\La)\,=\,
\sum_{k=0}^N\,(-1)^k\,\Zxi{\>k}(u\);t\);z\);\La)\,\der_u^{\>N\]-k}\).
\Tag{Zxiz}
$$
\par
The expression~\(FF) defines a vector\)-valued \raf/
$\FF^{\vox}_{\]\xi}(t\);z)$ of the \var/s $\txn$, $\zn$,
$$
\FF^{\vox}_{\]\xi}(t\);z)\,=\,\FF_{\]\xi}(t)\>\vox\;,
\vv-.1>
\Tag{FFtz}
$$
with the denominator
\vvn-.1>
$\prod_{i=1}^n\,\prod_{a=1}^{N-1}\,\prod_{j=1}^{\xi^a\!}\,(t^a_j-z_i)\,
\prod_{a=1}^{N-2}\,\prod_{i=1}^{\xi^a\!}\;
\prod_{j=1}^{\xi^{a+1}\!\!}\;(t^{a+1}_j\}-t^a_i)$, \,cf.~Lemma~\[FFpol].
The \fn/ $\FF^{\vox}_{\]\xi}(t\);z)$ coincides with
the \em{universal \wt/ \fn/} used in integral formulae for \hgeom/ \sol/s
of the \KZv/ \eq/s, see~\cite{Ma}\), \cite{SV}\), \cite{RSV}\).
\vsk-.3>
\Th{main4}
Let $\Mn$ be \hw/ \gnmod/s with \hw/s $\Lan$ and \hwv/s $v_1\lc v_n$.
Let $K$ be a diagonal matrix and $\tti=(\ttxn)$ a \sol/ of system \(Bae2).
Then for any $\kN\}$, we have
\vvn.2>
$$
\Gc^{\Moxs}_{k,K}(u\);z)\;\FF^{\vox}_{\]\xi}(\tti\);z)\,=\,
\Zxi{\>k}(u\);\tti\);z\);\La)\;\FF^{\vox}_{\]\xi}(\tti\);z)\,.
\Tag{GcFox}
$$
\endpro
\Pf.
The statement follows from Theorem~\[main3] and formulae \(Gcz)\), \(FFtz)\).
\epf
\Rem
Since \coeffs/ $\Gc_{1,K}(u)$ are central in $\Unx$, formulae~\(Gc12) and
\(MxiL) yield
\vvn-.4>
$$
\Gc^{\Moxs}_{1,K}(u\);z)\,=\,
\sum_{a=1}^N\)\biggl(\]K_a+\>\sum_{i=1}^n\){\La^a_i\over u-z_i}\)\biggr)
\,=\,\Zxi{\)1}(u\);t\);z;\)\La)\,,
$$
which implies Theorem~\[main3] for $k=1$.
\enddemo
\Prop{Fsing}
\back\;\cite{RV}
Let $K$ be the zero matrix, that is, ${K_1\]\lsym=K_N\]=0}$,
and $\tti=(\ttxn)$ an \offdiag/ \sol/ of the \Bae/s \(Bae2)\).
Then the \Bv/ $\FF^{\vox}_{\)\xi}(\tti\);z)$ is a singular vector \wrt/
the \,\$\gln$-\)action in $\Mox$ of \wt/ $\bigl(\)
\sum_{i=1}^n\La^1_i\]-\)\xi^1\],\,\sum_{i=1}^n\La^2_i\]+\)\xi^1\}-\)\xi^2\]
\lc\)\sum_{i=1}^n\La^N_i\!+\)\xi^{N\]-1}\)\bigr)$.
\endpro
\Pf.
The claim follows from identites for the universal \wt/ \fn/~\cite{RSV} and
quasiclassical asymptotics of \hgeom/ \sol/s of the \KZv/ \eq/s \cite{RV}\).
\epf
\Rem
Let $\Mn$ be any \gnmod/s, not necessarily \hw/ ones, and let ${v\in\Mox}$ be
any \wt/ singular vector \wrt/ the action of $\glnx$ in the module ${\Moxz}$.
Define \raf/s $r^v_1(u)\lc r^v_N(u)$ by the rule
$L_{aa}(u)\>v\)=\)r^v_a(u)\>v$, $\an\}$, and consider the \difl/ operator
\vvn-.2>
$$
\Mxi^{\>v}(u\),\der_u\);t)\,=
\prodr_{1\le a\le N\!}\biggl(\der_u-K_a-\>r^v_a(u)\,-\>
\sum_{i=1}^{\xi^{a-1}\!\!}\,{1\over u-t^{a-1}_i\}}\;+\>
\sum_{i=1}^{\xi^a}\,{1\over u-t^a_i}\biggr)\>.
$$
Let \>$\Zxi{\)1}^{\>v}(u\);t\);z)\lc\Zxi{\)N\}}^{\>v}(u\);t\);z)$
be its coefficients:
\vvn-.4>
$$
\Mxi^{\>v}(u\),\der_u\);t\);z\);\La)\,=\,
\sum_{k=0}^N\,(-1)^k\,\Zxi{\>k}^{\>v}(u\);t\);z)\,\der_u^{\>N\]-k}\).
\vv-.3>
$$
Similarly to Theorem~\[main4], for any \sol/ $\tti$ of the system
\vvn-.1>
$$
r^v_a(t^a_i)\,+\,\sum_{j=1}^{\xi^{a-1}\!\!}\;{1\over t^a_i\]-t^{a-1}_j}\;-\,
2\)\sum_{\tsize{j=1\atop j\ne i}}^{\xi^a}\,{1\over t^a_i\]-t^a_j}\;+\,
\sum_{j=1}^{\xi^{a+1}\!\!}{1\over t^a_i\]-t^{a+1}_j}\;=\,K_{a+1}\]-K_a\,,
\vv-.4>
$$
$\an-1$, \>$i=1\lc\xi^a\}$, we have
\vv.1>
$$
\Gc^{\Moxs}_{k,K}(u\);z)\;\FF_{\]\xi}(\tti\))\>v\,=\,
\Zxi{\>k}^{\>v}(u\);\tti\);z)\;\FF_{\]\xi}(\tti\))\>v\,.
$$
Moreover, for any $\xi$ and $t$,
\vvn-.5>
$$
\Gc^{\Moxs}_{1,K}(u\);z)\;\FF_{\]\xi}(t)\>v\,=\>
\Zxi{\)1}^{\>v}(u\);t\);z)\;\FF_{\]\xi}(t)\>v\,=\>
\sum_{a=1}^N\,\bigl(K_a+\>r^v_a(u)\bigr)\;\FF_{\]\xi}(t)\,v\,.
\vv-.5>
$$
\enddemo

\Sect[filt]{Filtration on $\Yn$}
In this section we are going to relate the series
$\Te_{1,Q}(u)\lc\Te_{N\},Q}(u)$ \wcoeff/ the Yangian $\Yn$ and the series
$\Gc_{1,K}(u)\lc\Gc_{N\},K}(u)$ \wcoeff/ the algebra $\Unx$, and to prove
Propositions~\[ADLA]\,--\,\[Gcsym].
The results of this section are essentially borrowed from \cite{T}\).
\vsk.2>
Given a filtered algebra $\Ac$ with an ascending filtration
${{}\ldots\sub\Ac_{\)s-1}\}\sub\Ac_{\)s}\}\lsym\sub\Ac}$, possibly double
infinite, we denote by ${\grad^\Ac_s\]:\Ac_{\)s}\]\to\Ac_{\)s}/\]\Ac_{\)s-1}}$
the canonical projection. We will identify the quotient spaces
$\Ac_{\)s}/\]\Ac_{\)s-1}$ with the corresponding homogeneous subspaces
in the algebra
$$
\grad\Ac\,=\,\Plus_{r\in\Z}\,\Ac_{\)r}/\]\Ac_{\)r-1}
\vv-.2>
$$
and will regard $\grad^\Ac_s$ as a map from $\Ac_{\)s}\]$ to $\grad\Ac$.
Abusing notation, we will write just $\grad_s$, dropping the superscript,
if it can cause no confusion.
\vsk.1>
For a filtered algebra $\Ac$, we consider the algebra ${\End(V)\ox\Ac}$
to be filtered by the subspaces
$\bigl(\End(V)\ox\Ac\bigr)\vpp s\}=\)\End(V)\ox\Ac_s\}$.
\par
The Yangian $\Yn$ admits a degree \fn/ \st/ ${\deg\)T_{ab}\#s\}=s-1}$ for
any $a\),\bn$, ${s=1,2,\ldots}$,
see Section 1.20 in~\cite{MNO}\); the \fn/ $\deg$ here coincides with
the \fn/ $\deg_{\)2}\]$ in \cite{MNO}\). Then $\Yn$ is a filtered algebra
with an ascending filtration $\Yn\vpp0\}\sub\Yn\vpp1\}\lsym\sub{}\Yn$, where
the subspace $\Yn\vpp s$ contains elements of degree not greater than $s$.
\Th{gradY}
\back\;\cite{MNO, Theorem~1.26\)}
We have ${\grad\bigl(\Yn\bigr)\>=\>\Unx}$.\uline\ Moreover,
\>$\grad_{s-1}(\)T_{ab}\#s)\)=\)e_{ba}\3s\]$, \>$a\),\bn\}$, \,$s=1,2,\ldots$.
\endpro
The \ainv/ $\opi$ is compatible with the filtration on $\Yn$.
For any ${X\]\in\Yn}$ we have $\deg\)\opi(X)\)=\)\deg\)X$,
and if $\deg\)X=r$, then
\vvn-.2>
$$
\grad_{\)r}\bigl(\opi(X)\bigr)\,=\,\tau(\grad_{\)r}\}X)\,.
\vv-.2>
\Tag{gradopi}
$$
\par
We extend the filtration from the Yangian $\Yn$ to the algebra $\Yn\uder$ by
the rule $\deg\)u\1\}=\deg\)\der_u\}=-\)1$. Notice that the indeterminates
$u\1\}$ and $\der_u\}$ do not commute obeying the following relation
$\der_u{\cdot}\>u^{-s}\}=u^{-s}\>\der_u\]-s\)u^{-s-1}\}$ instead. Clearly,
$$
\grad\bigl(\Yn\uder\bigr)\,=\,\Unx\uder\,.
$$
It is easy to see that the series ${T_{ab}(u)-\dl_{ab}\]\in\Yn\uder}$ has
degree $-\)1$ and
\vvn.2>
$$
\grad_{-1}\bigl(T_{ab}(u)-\dl_{ab}\bigr)\>=\>L_{ab}(u)\,,\qqq\Rlap{a\),\bn\).}
\Tag{TL}
$$
\vsk.1>
Further on we assume that in the definition \(Te) of the \tmcs/ $Q$ is a
series from $\End(V)\lbc\)\zt\)\rbc$ instead of being an element of $\End(V)$.
Then the obtained \tmcs/ $\Te_{1,Q}(u)\lc\Te_{N\},Q}(u)$ are power series in
$u\1\}$ and $\zt$ \wcoeff/ in $\Yn$. We also regard them as elements of
${\Yn\udzt}$. The results of Sections~\[:Yn] and~\[:Ba] naturally generalize to
the described setting. We extend the filtration from $\Yn\uder$ to ${\Yn\udzt}$
by setting $\deg\)\zt=-\)1$. Similarly, in formulae \(detL)\), \(DG) we will
assume $K$ being a series from $\End(V)\lbc\)\zt\)\rbc$.
\vsk.2>
Recall that $\Te_{0,Q}(u)=1$ by convention. For any $\koN\}$, set
\vvn-.4>
$$
\Se_{k,Q}(u)\,=\,{1\over(N\]-k)\)!}\,
\sum_{l=0}^k\,(-1)^{k-l}\,{(N\]-l\))\)!\over(k-l\))\)!}\ \Te_{l,Q}(u)\,.
\Tag{ST}
$$
For example, $\Se_{0,Q}(u)\)=\)1$ and $\Se_{1,Q}(u)\)=\)\Te_{1,Q}(u)-N\}$.
Formulae \(ST) for all $\koN\}$ taken together are \eqv/ to the identity
$$
\sum_{k=0}^N\,(-1)^k\,\Se_{k,Q}(u)\,y^{\)N\]-k}\)=\,
\sum_{l=0}^N\,(-1)^l\,\Te_{l,Q}(u)\,(y+1)^{\)N\]-\)l}\).
\Tag{STx}
$$
The series $\Se_{1,Q}(u)\lc\Se_{N\},Q}(u)$ are called the \em{modified \tmcs/}.
\Th{SkQ}
\back\;\cite{T}
Let ${\deg\)(Q-1)\le\]-\)1}$. Then
\>${\deg\)\bigl(\Se_{k,Q}(u)\bigr)\)=\)-\)k}$ \)for any $\koN\}$. Moreover,
\vvn-.4>
$$
\grad_{-k}\]\bigl(\Se_{k,Q}(u)\bigr)\>=\,\Gc_{k,K}(u)
\vv-.4>
\Tag{SeGc}
$$
where ${K=\grad_{-1}(Q-1)}$.
\endpro
\Pf.
Taking the \)\${(N\]-m)\)}$-th derivative of formula~\(STx) \wrt/ $y$ and
making there the substitution ${y=e^{\der_u}\}-1}$, we get
$$
\align
\sum_{k=0}^m\,(-1)^k\,
{(N\]-k)\)!\over(m-k)\)!}\ \Se_{k,Q}(u)\,(e^{\)\der_u}\}-1)^{m-k}\)=
\,\sum_{l=0}^m\,(-1)^l\,{(N\]-l\))\)!\over(m-l\))\)!}
\ \Te_{l,Q}(u)\,e^{(m-l)\)\der_u}\)={}\kern-2.79em &
\Tagg{DS}
\\
\nn4>
{}=\,(N\]-m)\)!\ \Dg_{m,Q}(u\),\]\der_u)\,e^{\)m\der_u}\>,\kern-3em &
\\
\cnn-.7>
\endalign
$$
\cf. \(DgmT)\). Exploiting the standard identity
\vvn-.2>
$$
\gather
\sum_{r=0}^s\,{(-1)^r\over r\)!\,(s-r)\)!}\;=\,0\,,\quad\Rlap{\qquad s\ge 1\,,}
\\
\nn-1>
\Text{we obtain that}
\nn-2>
\Se_{k,Q}(u)\,=\,{1\over(N\]-k)\)!}\,\sum_{m=0}^k\,
(-1)^m\,{(N\]-m)\)!\over(k-m)\)!}\ \Dg_{m,Q}(u\),\der_u)\,
e^{\)m\der_u}\)(e^{\)\der_u}\}-1)^{k-m}\).\kern-2em
\Tag{SD}
\\
\cnn-.2>
\endgather
$$
Since \,${\deg\)\bigl(T(u)-1\bigr)=-\)1}$,
\,${\deg\)(e^{\)\der_u}\}-1)=-\)1}$, \>and \,${\deg\)(Q-1)\le\]-\)1}$,
formula \(Dgm) implies that ${\deg\)\bigl(\Dg_{m,Q}(u\),\]\der_u)\bigr)=-\)m}$.
Hence, ${\deg\)\bigl(\Se_{k,Q}(u)\bigr)\)=\)-\)k}$ by formula \(SD)\).
\vsk.2>
Observe that ${\grad_{-1}\bigl(T(u)-1\bigr)=L(u)}$, \cf.~\(TL)\),
and ${\grad_{-1}(e^{\)\der_u}\}-1)=\der_u}$. Then computing
$\grad_{-m}\}\Dg_{m,Q}(u\),\der_u)$ in two ways using either
\uugood
formula~\(Dgm) or formula~\(DS), we get
\vvn-.5>
$$
\align
({\tr_{V^{\]\ox m}}\:\}}\ox\id\))\Bigl(\]
\bigl(\der_u\]-K\"1\}-L\"{1,\)m+1}(u)\bigr)\)\ldots\)
\bigl(\der_u\]-K\"m\}-L\"{m,\)m+1}(u)\bigr)\;
\AA\"{1\ldots m}\6m\)\Bigr)\,={} &
\\
\nn5>
{}=\,{1\over(N\]-m)\)!}\;\sum_{k=0}^m\,(-1)^k\,{(N\]-k)\)!\over(m-k)\)!}\;
\grad_{-k}\bigl(\Se_{k,Q}(u)\bigr)\,\der_u^{\>m-k}\> & .
\\
\cnn-.2>
\endalign
$$
The last step is to compare the obtained formula for $m=N$ with
formulae~\(detL) and \(DG)\), which completes the proof.
\epf
\Pf of Propositions~\[ADLA]\,--\,\[Gcsym].
Let $Q=1+\zt K$. Propositions~\[Tekl], \[Tesym] and \[ADTA] together with
formulae~\(ST)\), \(SeGc), \(gradopi) and the equality ${\grad_{-N}
\bigl(\Dg_{\]N\],Q}\](u\),\der_u)\bigr)\)=\)\De_{\zt\]K}\](u\),\der_u)}$
imply \resp/ Propositions~\[Gckl], \[Gcsym] and~\[ADLA] with the matrix $K$
replaced by $\zt K$. Since all the expressions
$\Gc_{1,\)\zt\]K}(u)\lc\Gc_{N\},\)\zt\]K}(u)$ and $\De_{\zt\]K}\](u\),\der_u)$
are \pol/s in $\zt$, we can evaluate them at $\zt=1$, which completes
the proof.
\epf

\Sect[proof]{Proof of Theorem~\[main]}
We will prove Theorem~\[main] by induction \wrt/ $N\]$. The key points are
Proposition~\[BB1] and \[ADB], and formula \(TeX)\).
\vsk.2>
Set $W\}=\Cnn\!$. Let $\wbn$ be the standard basis of $W\}$ and $\vbn\}$
the standard basis of $V\}=\Cn\!$. We identify $W\}$ with its image in $V\}$
under the embedding $\wb_a\map\)\vb_{a+1}$, $\an-1$.
\vsk.2>
Let ${P\+{N-1}\]\in\End(W^{\ox2})}$ be the flip map and
${R\+{N-1}(u)=\)u+P\+{N-1}}\}$ the \rat/ \Rm/ used in the definition
of the Yangian $\Ynn$. The \Rm/ $R(u)$ preserves the subspace
${W^{\ox2}\}\sub V^{\ox2}}\!$ and the restriction of $R(u)$ on $W^{\ox2}\!$
coincides with $R\+{N-1}(u)$.
\vsk.2>
We will use the embedding ${\psi\):\)\Ynn\>\hto\>\Yn}$:
\,${\psi\bigl(T_{ab}\+{N-1}(u)\bigr)\)=\>T_{a+1,\)b+1}(u)}$,\uline
\;$a\),\bn-1$, where $T_{ab}\+{N-1}(u)$ are series \(Tab) for the algebra
$\Ynn$. It is clear that $\psi\bigl(\Ynnp\bigr)\)\sub\)\Ynp$.
\vsk.3>
Let $W(x)$ be the \eval/ \Ynnmod/ and ${\pi(x):\Ynn\to\End(W)}$ be
the corresponding \hom/. In other words
\vvn.3>
$$
\pi(x)\):\)T\+{N-1}(u)\,\map\,(u-x)\1R\+{N-1}(u-x)\,.
\vv.1>
\Tag{pix}
$$
Introduce a map \,$\psi(\xr)\):\)\Ynn\>\to\>\End\bigl(W^{\ox r}\bigr)\ox\Yn$,
$$
\psi(\xr)\,=\,\bigl(\pi(x_1)\lox\pi(x_r)\ox\psi\bigr)\o
\bigl(\Dl\*{N-1}\bigr)\2r\),
\Tag{psih}
$$
where \>${\bigl(\Dl\*{N-1}\bigr)\2r\!:
\)\Ynn\>\to\>\bigl(\Ynn\]\bigr)^{\]\ox(r+1)}}$
is the multiple coproduct. For any ${X\}\in\End(W)}$ set ${\nu(X)=X\)\wb_1}$.
Define a map \,${\psit(\xr)\):\)\Ynn\>\to\>W^{\ox r}\!\ox\Yn}$,
$$
\psit(\xr)\,=\,(\nu^{\>\ox r}\!\]\ox\id\))\o\psi(\xr)\,.
$$
\Lm{psiY}
${\psit(\xr)\)\bigl(\)\Ynnp\bigr)\)\sub\)W^{\ox r}\!\ox\Ynp}$.
\endpro
\Pf.
The claim follows from the facts that $\Ynnp$ is a coideal in $\Ynn$, and
the vector $\wb_1\}$ is a singular vector \wrt/ the action of $\Ynn$.
\epf
Consider the embedding ${\pho\):\)Y(\gl_{N-2})\>\hto\>\Ynn}$:
\vvn.1>
\,${\pho\bigl(T_{ab}\+{N-2}(u)\bigr)\)=\>T_{a+1,\)b+1}\+{N-1}(u)}$,\uline
\;$a\),\bn-2$, similar to the embedding $\psi$. Recall that
$\tht:\Yn\to\Yn/\Ynp$ is the canonical projection.
\Lm{phipsi}
We have \,$(\id^{\ox r}\!\]\ox\)\tht)\o\psit(\xr)\o\pho\,=\,
\wb_1^{\)\ox r}\!\]\ox(\tht\o\psi\o\pho)$.
\endpro
For any element ${g\in W^{\ox r}\!\ox\Yn}$ we define its components
$g^{\)a_1\lc a_r}\}$ by the rule \uline
$g\>=\!\sum_{a_1\lc a_r=1}^{N-1}\!\wb_{a_1}\!\lox\)\wb_{a_r}\!\ox
g^{\)a_1\]\lc a_r}$.
\vsk.3>
Set \>$\bar\xi=\)(\xi^2\]\lc\xi^{\)N-1})$ and \,$\tbar\,=\,
(t^2_1\lc t^2_{\xi^2};\;\ldots\;;t^{\)N-1}_1\}\lc t^{\)N-1}_{\xi^{N-1}})$.
\Prop{BB1}
\back\;\cite{TV1, Theorem~3.4.2\)}
We have
\vvn-.1>
$$
\BB_{\)\xi}(t)\,=\!\sum_{a_1\]\lc a_{\xi^1}=\)1}^{N-1}\!
T_{1,a_1+1}(t^1_1)\ldots\>T_{1,a_{\xi^1}+1}(t^1_{\xi^1}\])\,
\Bigl(\psit(\txine\])\bigl(\)\BB\*{N-1}_{\bar\xi}(\)\tbar)\bigr)\}\Bigr)
^{\}a_1\]\lc a_{\xi^1}}.\kern-1em
\vvn-.4>
\Tag{BBe}
$$
\endpro
\Pf of Lemma \[thtBB].
Lemma~\[phipsi], Proposition~\[BB1] and formulae \(pix)\), \(psih) imply
that the denominator of $\tht\bigl(\BB_{\)\xi}(t)\bigr)$ is at most
${\prod_{a=1}^{N-2}\)\prod_{i=1}^{\xi^a\!}\>\prod_{j=1}^{\xi^{a+1}\!}\,
(t^{a+1}_j\!-t^a_i)}$. Now the statement follows from formula~\(BB)\).
\epf
Since ${V\}=\)\C\)\vb_1\]\oplus W}\}$ we have
${\Vw k\}=\Wtw{(k-1)}\oplus\Ww k}\}$,
where the first summand is spanned by the vectors
${\vb_1\]\wg\vb_{a_1}\!\lwg\vb_{a_{k-1}}}$, \,$2\le a_1\lsym<a_{k-1}\le N$,
and the second summand is spanned by the vectors ${\vb_{b_1}\!\lwg\vb_{b_k}}$,
\,$2\le b_1\lsym<b_k\le N$. Abusing notation, we will often identify
$\Ww{(k-1)}\}$ and $\vb_1\]\wg\Ww{(k-1)}\}$ taking a vector $\xb$ to
$\vb_1\}\wg\xb$.
\vsk.2>
The \Rm/ $\Rwk(u)$ acting in ${\Vw k\}\ox V}$ preserves the subspaces
$$
\Wtw{(k-1)}\}\ox\C\)\vb_1\,,\ \quad
\bigl(\Wtw{(k-1)}\ox W\)\bigr)\)\oplus\)(\Ww k\}\ox\C\)\vb_1)\,,\ \quad
\Ww k\}\ox W\>.
$$
Its restrictions onto these subspaces are
\vvn-.2>
$$
\gather
\Rwk(u)\vst{\Wtws{(k-1)}\]\ox\)\C\)\vb_1}\!=\,u+1\,,\qquad
\Rwk(u)\vst{\Ww k\]\ox\)W}=\>\Rwk\+{N-1}(u)\,,
\Tag{RN-1}
\\
\nn10>
\Rwk(u)\vst{(\Wtws{(k-1)}\]\ox W)\)\oplus\)(\Wws k\]\ox\)\C\vb_1)}\)=\,
\pmatrix\]\Rwe{(k-1)}\+{N-1}(u) & \;\Sti'\,\\\nn5> \Sti & \;u\,\endpmatrix\,,
\\
\cnn-.2>
\endgather
$$
where $\Sti\bigl((\vb_1\}\wg\xb)\ox\wb\bigr)\)=\)(\wb\]\wg\xb)\ox\vb_1$,
and we will not use the explicit form of $\Sti'$. A similar decomposition
can be written for the \Rm/ $\Rwwk(u)$.
\vsk.2>
In what follows we regard $T(u)$ and $\Tw k(u)$ as matrices over the algebra
$\Yn$ and introduce their submatrices induced by the decompositions
${V\}=\)\C\)\vb_1\]\oplus W}\}$ and \uline
${\Vw k\}=\)\Wtw{(k-1)}\oplus\Ww k}\}$:
\vvn-.2>
$$
T(u)\,=\,\pmatrix A(u) & B(u)\\\nn2> C(u) & D(u)\endpmatrix\,,\qqq
\Tw k(u)\,=\,\pmatrix \Ah(u) & \Bh(u)\\\nn2> \Ch(u) & \Dh(u)\endpmatrix\,.
\vv-.2>
\Tag{ABCD}
$$
For example, $A(u)=T_{11}(u)$, $B(u)=\bigl(T_{12}(u)\lc T_{1N}(u)\bigr)$,
$D(u)\)=\)\bigl(T_{ij}(u)\bigr)_{i,j\)=\)2}^N$.
\vsk.1>
Denote by $\HY(L\),M)$ the space ${\Hom(L\),M)\ox\Yn}$ of matrices
with noncommuting entries. We call $L$ the \em{domain} of those matrices.
\Coeffs/ $B\),\>D\),\>\Ah\>,\>\Bh\),\alb\>\Dh$ \resp/ belong to $\HY(W,\C)$,
$\HY(W,W)$, $\HY(\Ww{(k-1)},\Ww{(k-1)})$,
$\HY(\Ww k,\Ww{(k-1)})$, $\HY(\Ww k,\Ww k)$. For ${k=1}$ we have
${\Ww{(k-1)}\]=\)\C}$, ${\Ah=A}$, ${\Bh=B}$, ${\Dh=D}$.
\vsk.1>
Set
\vvn-.4>
$$
\Rbar(u)\,=\,\){1\over u}\,R\+{N-1}(u)\,,\kern1.6em
\Rt(u)\,=\,\){1\over u+1}\,\Rwe{(k-1)}\+{N-1}(u)\,,\kern1.6em
\Rh(u)\,=\,\){1\over u}\,\Rwk\+{N-1}(u)\,.
\vv-.2>
$$
Define a map ${S:\Ww{(k-1)}\}\ox W\]\to\)\Ww k}\}$ by the rule
\vvn-.2>
$$
S(\xb\ox\wb)\,=\,\wb\]\wg\xb\,.
\vv-.2>
\Tag{Sxw}
$$
We will use the following commutation relations that are obtained from
formulae~\(RTTw)\), \(RN-1)\), Corollary~\[inw] and the equality
$S\>\Rt(u-1)\)=\)(u-k)\>S$:
\vvn-.2>
\<BhB>\<AhB>
$$
\gather
\Bh\41(u)\>B\42(v)\,=\,\){u-v\over u-v+1}\;
B\42(v)\>\Bh\41(u)\,\Rh\"{12}(u-v)\,,\kern-1.6em
\Tag{BhB}
\\
\ald
\nn7>
\Ah\"1(u)\>B\42(v)\,=\,B\42(v)\>\Ah\"1(u)\>\Rt\"{12}(u-v-1)\,+\,
{1\over u-v}\;\Bh(u)\>S\4{12}A(v)\,,\kern-2em
\Tag{AhB}
\\
\nn7>
\Dh\"1(u)\>B\42(v)\,=\,B\42(v)\>\Dh\"1(u)\>\Rh\"{12}(u-v)\,-\,
{1\over u-v}\;S\>\Bh\41(u)\>D\"2(v)\,.\kern-2em
\Tag{DhB}
\endgather
$$
The equalities hold respectively in ${\HY(\Ww k\}\ox W,\Ww{(k-1)})}$,
${\HY(\Ww{(k-1)}\}\ox W,\Ww{(k-1)})}$ and ${\HY(\Ww k\}\ox W,\Ww k)}$,
and the superscripts in brackets indicate what tensor factors are domains
of the respective matrices. If $k=1$, relations \(BhB)\,--\,\(DhB) become
\<BBR>\<AB>
\vvn-.2>
$$
\gather
B\41(u)\>B\42(v)\,=\,\){u-v\over u-v+1}\;
B\42(v)\>B\41(u)\,\Rbar\"{12}(u-v)\,,\kern-1.6em
\Tag{BBR}
\\
\nn7>
A(u)\>B(v)\,=\,\){u-v-1\over u-v}\;B(v)\>A(u)\,+\,
{1\over u-v}\;B(u)\>A(v)\,,\kern-1.4em
\Tag{AB}
\\
\nn7>
D\"1(u)\>B\42(v)\,=\,B\42(v)\>D\"1(u)\>\Rbar\"{12}(u-v)\,-\,
{1\over u-v}\;B\41(u)\>D\"2(v)\,.\kern-2em
\Tag{DB}
\\
\ungood
\endgather
$$
\par
Let \)${\Rc(u)\)=\)(u+1)\1\)P\+{N-1}\>R\+{N-1}(u)}$.
For an expression $f(\ur)$ \wcoeff/ matrices with the domain
$W^{\ox r}\}$ and a simple transposition $(i\),i+1)$, ${\ik-1}$, set
\vvn-.1>
$$
\0{(i,i+1)}f(\ur)\,=\,f(u_1\lc u_{i-1},u_{i+1},u_i\),u_{i+2}\lc u_r)\,
\Rc\"{i,i+1}(u_i\]-u_{i+1})\,.
\Tag{Ract}
$$
The matrix $\Rc(u)$ has the properties $\Rc(u)\>\Rc(-u)\>=\>1$ and
\vvn.1>
$$
\Rc\"{12}(u-v)\>\Rc\"{23}(u)\>\Rc\"{12}(v)\,=\,
\Rc\"{23}(v)\>\Rc\"{12}(u)\>\Rc\"{23}(u-v)\,,
\vvn-.1>
$$
\cf. \(YB). This yields the following lemma.
\Lm{RSym}
Formula \(Ract) extends to the action of the \symg/ $S_r$ on expressions
$f(\ur)$ \wcoeff/ matrices with the domain $W^{\ox r}\}$:
${f\)\map\0{\si\!}f}$, ${\si\in S_r}$.
\endpro
\nt
By formula \(BBR) the expression $B\41(u_1)\ldots\)B\4r(u_r)$ is invariant
under the action \(Ract) of the \symg/ $S_r$.
\vsk.2>
For an expression $f(\ur)$ \wcoeff/ matrices with the domain
$W^{\ox r}\!$, set
$$
\Symr_{\)\ur}\"\lr f(\ur)\,=\sum_{\si\in S_r}\0{\si\!}f(\ur)\,.
\vv->
\Tag{SymR}
$$
\Prop{ADB}
$$
\align
& \Ah\"0(u)\>B\41(u_1)\ldots\)B\4r(u_r)\,={}
\Tagg{ABB}
\\
\nn8>
& \kern-2.78em {}=\,B\41(u_1)\ldots\)B\4r(u_r)\,
\Ah\"0(u)\>\Rt\"{0\)r}(u-u_r\]-1)\ldots\Rt\"{0\)1}(u-u_1\]-1)\,+{}
\\
\nn4>
{}+\;{1\over(r-1)\)!}\, & \;\Bh(u)
\Symr_{\)\ur}\"{\lr}\biggl({1\over u-u_1}\;
\prod_{i=2}^r\,{u_1\]-u_i\]-1\over u_1-u_i}\;
S\4{01}\>B\42(u_2)\ldots\)B\4r(u_r)\>A(u_1)\biggr)\,,
\endalign
$$
\vvn-1.3>
$$
\alignat2
& \,\Dh\"0(u)\>B\41(u_1)\ldots\)B\4r(u_r)\,={}
\Tagg{DBB}
\\
\nn8>
& \Rlap{{}=\,B\41(u_1)\ldots\)B\4r(u_r)\,
\Dh\"0(u)\>\Rh\"{0\)r}(u-u_r)\ldots\)\Rh\"{01}(u-u_1)\,-{}}
\\
\nn5>
& {}\)-\;{1\over(r-1)\)!}\,\;S\>\Bh\40(u)
\Symr_{\)\ur}\"{\lr}\biggl({1\over u-u_1}\;
&& B\42(u_2)\ldots\)B\4r(u_r)\>D\"1(u_1)\,\x{}\kern-3.8em
\\
\nn->
&&& \Llap{{}\x\,{}}
\Rbar\"{1\)r}(u_1\]-u_r)\ldots\)\Rbar\"{12}(u_1\]-u_2)\}\biggr)\,,\kern-3.8em
\\
\cnn-.5>
\endalignat
$$
where the tensor factors are counted by \)$0\),\)1\)\lc r$.
\endpro
\Pf.
The statement follows from relations \(BhB)\,--\,\(DB) by induction \wrt/ $r$.
We apply formula \(AhB) or \(DhB) to the product of the first factors in \lhs/
and then use the induction assumption.
\epf
\Rem
Formulae \(ABB) and \(DBB) have the following structure. The first term in
\rhs/ comes from repeated usage of the first term in \rhs/ of relation \(AhB)
or \(DhB), \resp/. The second term, involving symmetrization, is effectively
determined by the fact that the whole expression in \rhs/ is regular at $u=u_i$
for any $i=1\lc r$, and is invariant \wrt/ action \(Ract) of the \symg/ $S_r$.
The symmetrized expression is obtained by applying once the second term in
\rhs/ of the relevant relation \(AhB) or \(DhB) followed by repeated usage
of the first term of the respective relation.
\enddemo
\vsk-.3>
Let $Q\)=\]\sum_{a=1}^N Q_aE\+N_{aa}\in\End(V)$ be the diagonal matrix
from Theorem~\[main]. Set
\vvn-.6>
$$
\Qbar\)=\)\sum_{a=1}^{N-1} Q_{a+1}\>E\+{N-1}_{aa}\>\in\,\End(W)\,.
\vv-.2>
\Tag{Qbar}
$$
Denote $\Qt\)=\Qbw{(k-1)}$ and \>$\Qh\)=\Qbw k\}$.
Then formula~\(Te) implies that
\vvn.1>
$$
\Te_{k,Q}(u)\,=\,Q_1\>\tr_{\Ww{(k-1)}}\bigl(\Qt\>\Ah(u)\bigr)\,+\,
\tr_{\Wws k}\bigl(\Qh\>\Dh(u)\bigr)\,.
\vv.1>
\Tag{TAD}
$$
\vsk.1>
To prove formula \(TeB) we will employ Propositions~\[BB1] and \[ADB].
In our present notation, relation~\(BBe) reads as follows:
\vvn.1>
$$
\BB_{\)\xi}(t)\,=\,B\41(t^1_1)\ldots\)B\4{\xi_1}(t^1_{\xi^1}\])
\ \psit(\txine\])\bigl(\)\BB\*{N-1}_{\)\bar\xi}(\)\tbar)\bigr)\,.
\vv.1>
\Tag{BBB}
$$
We substitute \(TAD) and \(BBB) into \lhs/ of formula \(TeB) and use
relations \(ABB)\), \(DBB) to move $\Ah(u)$ and $\Dh(u)$ through the product
$B\41(t^1_1)\ldots\>B\4{\xi_1}(t^1_{\xi^1}\])$.
\vsk.1>
To simplify writing, let $r=\)\xi^1\}$ and $u_i=\)t^1_i$, $i=1\lc r$.
We have
\vvn.2>
$$
\align
& \hp{B\41(u_1)\ldots\){}}\Llap{\Te_{k,Q}(u)\>}B\41(u_1)\ldots\)B\4r(u_r)\,={}
\Tagg{TeBB}
\\
\nn8>
{}=\,{} &B\41(u_1)\ldots\)B\4r(u_r)\,\x{}
\\
\nn5>
& \aligned
{}\!\!\]\x\)\Bigl(\>Q_1
({\tr_{\Wws{(k-1)}}}\}\ox\id^{\ox r})\bigl(\Qt\"0\>\Ah\"0(u)\>
\Rt\"{0\)r}(u-u_r\]-1)\ldots\Rt\"{0\)1}(u-u_1\]-1)\bigr)\,+{} &
\\
\nn5>
{}+\,({\tr_{\Wws k}}\}\ox\id^{\ox r})\bigl(\Qh\"0\>\Dh\"0(u)\>
\Rh\"{0\)r}(u-u_r)\ldots\)\Rh\"{01}(u-u_1)\bigr)\}\Bigr)\>& {}\)+{}\kern-.6em
\endaligned
\\
\nn6>
{}+\,{} & {1\over(r-1)\)!}\,\Symr_{\)\ur}\"{\lr}
\biggl[\;{1\over u-u_1}\,\Be\41_Q(u)\>B\42(u_2)\ldots\)B\4r(u_r)\,\x{}
\\
\nn3>
& \ \){}\x\,\biggl(Q_1\,
\prod_{i=2}^r\,{u_1\]-u_i\]-1\over u_1-u_i}\;\>A(u_1)\,-\,
\Qbar\"1\)D\"1(u_1)\>\Rbar\"{1\)r}(u_1\]-u_r)\ldots\)
\Rbar\"{12}(u_1\]-u_2)\]\biggr)\biggr]\,,\kern-.6em
\endalign
$$
where the tensor factors for the products under the traces are counted by
\)$0\),\)1\)\lc r$. Here
$$
\Be_Q(u)\,=\,
({\tr_{\Wws{(k-1)}}}\}\ox\id\))\bigl(\)\Qbw{(k-1)}\Bh(u)\>S\)\bigr)
\vvu.5>
\uugood
\vvu-.5>
\Tag{BeQ}
$$
is a series \wcoeff/ $\HY(W,\C)$. Besides this, we use the equality
$$
({\tr_{\Wws k}}\}\ox\id\))\Bigl(\Qbw kS\>\bigl(\Bh(u)\ox\id\)\bigr)\]
\Bigr)\,=\,\Be_Q(u)\,\Qbar
$$
that follows from the cyclic property of the trace and
the formula ${\Qbw kS\)=\)S\,(\Qbw{(k-1)}\}\ox\Qbar)}$.
Notice that $\Be_Q(u)=B(u)$ if $k=1$.
\vsk.2>
To proceed further, we use the next lemmas.
\Lm{DRR}
\vvn-.2>
$$
\gather
D\"0(u)\>\Rbar\"{0\)r}(u-u_r)\ldots\)\Rbar\"{0\)1}(u-u_1)\,=\,
\psi(\ur)\bigl(\)T\+{N-1}(u)\bigr)\,,
\\
\nn5>
\Dh\"0(u)\>\Rh\"{0\)r}(u-u_r)\ldots\)\Rh\"{0\)1}(u-u_1)\,=\,
\psi(\ur)\bigll(\]\bigl(\)T\+{N-1}(u)\bigr)^{\]\wg k}\bigrr)\,.
\\
\cnn-.3>
\endgather
$$
\endpro
\Pf.
The equalities follow from formula \(DlT)\), Lemma~\[DlTw] and
formula \(psih)\).
\epf
\Lm{ARR}
For any $X\]\in\Ynn$,
\<Apsih>
$$
\gather
A(u)\,\)\psi(\ur)(X)\,\simeq\,\psi(\ur)(X)\,A(u)\,,
\Tag{Apsih}
\\
\nn8>
{\align
\Ah\"0(u)\>\Rt\"{0\)r}(u-u_r\]-1)\ldots\)\Rt\"{0\)1}(u-u_1\]-1)\;
\psi(\ur)(X)\,\simeq{} &
\Tagg{AhRpsih}
\\
\nn4>
{}\simeq\,\prod_{i=1}^r {u-u_i\]-1\over u-u_i}\;\>
A(u)\;\psi(\ur)\bigll(\]\bigl(\)T\+{N-1}(u-1)\bigr)^{\]\wg(k-1)}X\;
& \}\!{}\bigrr)\,.
\endalign}
\\
\cnn-.2>
\endgather
$$
\endpro
\Pf.
Denote by $\Ynx$ the left ideal in $\Yn$ generated by \coeffs/ $T_{a1}(u)$
for $2\le a\le N\}$. It is a subideal of $\Ynp$. Relations \(Tabcd) imply
that for any ${Z\]\in\Ynn}$ and ${C\]\in\Ynx}$ \coeffs/
$\bigl[\)T\+N_{11}(u)\>,\psi(Z)\)\bigr]$ and the product $C\>\psi(Z)$ belong
to $\Ynx$. Then formula \(Apsih) holds because ${A(u)=T\+N_{11}(u)}$ and
the coefficients of components of $\psi(\ur)(X)$ belong to
$\psi\bigl(\Ynn\]\bigr)$.
\vsk.1>
To prove formula~\(AhRpsih) we first observe that the coefficiens of entries
of the matrix ${\Ah(u)-A(u)\>D^{\wg(k-1)}(u-1)}$ belong to $\Ynx$; this fact
follows from the formula for $\Ah(u)$ similar to the first formula of~\(qdet)
for the identity \perm/ $\si$. Then we proceed like in Lemma~\[DRR] and use
the fact that $\psi(\ur)$ is a homomorphism of algebras.
The calculation goes as follows:
$$
\alignat2
& \Ah\"0(u)\>\Rt\"{0\)r}(u-u_r\]-1)\ldots\)\Rt\"{0\)1}(u-u_1\]-1)\;
\psi(\ur)(X)\,\simeq{}
\\
\nn6>
& \Rlap{{}\simeq\,A(u)\>\bigl(D^{\wg(k-1)}(u-1)\bigr)^{\][\)0\)]}\>
\Rt\"{0\)r}(u-u_r\]-1)\ldots\)\Rt\"{0\)1}(u-u_1\]-1)\;
\psi(\ur)(X)\,\simeq{}}
\\
\nn6>
& {}\simeq\,A(u)\;\psi(\ur)\bigll(\]\bigl(\)T\+{N-1}(u-1)\bigr)^{\]\wg(k-1)}\bigrr)
\,\psi(\ur)(X)\,\prod_{i=1}^r {u-u_i\]-1\over u-u_i}\;\simeq{} &&
\\
\nn2>
&& \Llap{{}\simeq\,A(u)\;\psi(\ur)\bigll(\]\bigl(\)T\+{N-1}(u-1)
\bigr)^{\]\wg(k-1)}X\bigrr)\,\prod_{i=1}^r {u-u_i\]-1\over u-u_i}\,\)} &
\).\kern1.6em
\\
\cnn->
\Text{\qed}
\cnn-.5>
\ungood
\endalignat
$$
\edemo
\vsk.5>
For an expression $f(v)$ set
\>$\res_{v\)=\)w}f(u)\)=\)\bigl((v-w)\)f(v)\bigr)\vst{v\)=\)w}\]$,
\vvn-.2>
if the substitution makes sense. Observe that
$$
\alignat2
& \Qbar\"1\)D\"1(u_1)\>
\Rbar\"{1\)r}(u_1\]-u_r)\ldots\)\Rbar\"{12}(u_1\]-u_2)\,={}
\Tag{Dres}
\\
\nn5>
& \!\}{}=\res_{v\)=\)u_1\!}\>\Bigl(({\tr_W}\}\ox\id^{\ox r})
\bigl(\)\Qbar\"0\)D\"0(v)\>
\Rbar\"{0\)r}(v-u_r)\ldots\)\Rbar\"{0\)1}(v-u_1)\bigr)\}\Bigr)\,={}\!
\kern-2em &&
\\
\nn3>
&&\Llap{{}=\,\res_{v\)=\)u_1\!}\Bigl(\psi(\ur)
\bigl(\)\Te\+{N-1}_{1,\Qbar}(v)\bigr)\}\Bigr)\,,\!\!}\kern-2em &
\\
\cnn-.8>
\endalignat
$$
because $\res_{v\)=\)u_1\!}\)\bigl(\Rbar(v-u_1)\bigr)\)=P\+{N-1}\}$.
\vsk.2>
Altogether, for any $X\]\in\Ynn$ we have
$$
\align
& \hp{{}\!\!\simeq\,\biggl\lb B\41}\,
\Te_{k,Q}(u)\,B\41(u_1)\ldots\)B\4r(u_r)\;\psi(\ur)(X)\,\simeq{}
\Tagg{TeX}
\\
\nn6>
& \aligned
{}\!\!\!\!\simeq\,\biggl\lb B\41(u_1)\ldots\)B\4r(u_r)\,
\Bigl(\psi(\ur)\)\bigl(\)\Te\+{N-1}_{k-1,\Qbar}(u-1)\,X\bigr)\;
Q_1\>A(u)\,\prod_{i=1}^r {u-u_i\]-1\over u-u_i}\>+{} &
\\
\nn6>
{}+\,\psi(\ur)\)\bigl(\)\Te\+{N-1}_{k,\Qbar}(u)\,X\bigr)\}\Bigr)\> &{}+{}
\endaligned
\\
\nn8>
& {}\!\!+\,{1\over(r-1)\)!}\,\Symr_{\)\ur}\"{\lr}
\biggl[\;{1\over u-u_1}\;\Be\41_Q(u)\>B\42(u_2)\ldots\)B\4r(u_r)\>\x{}
\\
\nn6>
& {}\!\!\x\>\biggl(\psi(\ur)(X)\ Q_1\,A(u_1)\,
\prod_{i=2}^r\,{u_1\]-u_i\]-1\over u_1\]-u_i}\,-\)
\res_{v\)=\)u_1\!}\Bigl(\psi(\ur)\bigl(\)\Te\+{N-1}_{1,\Qbar}(v)\,X\bigr)
\Bigr)\}\biggr)\biggr]\biggr\rb\,.\kern-.8em
\endalign
$$
\vsk-.1>
The next step is to apply both sides of formula \(TeX) to the vector
${\wb_1^{\)\ox r}}\!$. That amounts to replacing there the map $\psi(\ur)$
by $\psit(\ur)$, and changing the symmetrization $\Symr_{\)\ur}\"{\lr}\}$,
\)cf.~\(SymR)\), to the ordinary symmetrization $\Sym_{\)\ur}\}$, because
${\Rc(u)\>\wb_1\}\ox\)\wb_1=\wb_1\}\ox\)\wb_1}$.
\vsk.2>
Now we replace $X$ in formula~\(TeX) by \)$\BB\*{N-1}_{\)\bar\xi}(\)\tbar)$,
cf.~\(BBB)\), and employ the induction assumption. We recall that
${A(u)=T_{11}(u)}$, ${r=\xi^1\}}$, and ${u_i\]=t^1_i}$, $i=1\lc r$.
We also use Lemma~\[DlTaa] and equalities
\vvn-.8>
$$
\gather
\pi(x)\bigl(T\+{N-1}_{11}(u)\bigr)\>\wb_1\>=\,{u-x+1\over u-x}\;\wb_1\,,
\\
\nn5>
\pi(x)\bigl(T\+{N-1}_{aa}(u)\bigr)\>\wb_1\>=\,\wb_1\,,\qquad a=2\lc N-1\,,
\endgather
$$
cf.~\(pix)\), which imply that
\vvn->
$$
\gather
\psit(\txine\])\bigl(T\+{N-1}_{11}(u)\bigr)\,\simeq\,\wb_1^{\)\ox r}\!\]\ox
T_{22}(u)\;\prod_{i=1}^{\xi^1\!}\,{u-t^1_i\]+1\over u-t^1_i\]}\;,
\\
\nn5>
\psit(\txine\])\bigl(T\+{N-1}_{aa}(u)\bigr)\,\simeq\,\wb_1^{\)\ox r}\!\]\ox
T_{a+1,a+1}(u)\,,\qquad a=2\lc N-1\,.
\endgather
$$
Besides this, we use Lemmas~\[Taa0] and \[psiY].
As a result, we transform formula \(TeX) to
\vvn-.2>
$$
\alignat2
& \Rlap{\Te_{k,Q}(u)\>\BB_{\)\xi}(t)\,\simeq{}}
\Tagg{TeB2}
\\
\nn3>
{}\!\simeq\,\BB_{\)\xi} &{} (t)\,\biggl(\}Q_1\)T_{11}(u)\,&
\prod_{i=1}^{\xi^1\!}\,{u-t^1_i\]-1\over u-t^1_i}\,
\sum_{\ab}\,\prod_{j=2}^k\,\Xxi{\>a_j}(u-j+1\);t)\,+{} &
\\
\nn2>
&& {}+\,\sum_{\bb}\,\prod_{j=1}^k\,\Xxi{\)b_j}(u-j+1\);t)\}\biggr)\> &
{}+\,\Uxi(u\);t)\,,
\endalignat
$$
cf.~\(Xxi)\), where the first sum is taken over all \$(k-1)\)$-tuples
$\ab=(a_2\lc a_k)$ \st/ $2\le a_2\lsym<a_k\le N\]$, the second sum is taken
over all \$k\)$-tuples $\bb=(b_1\lc b_k)$ \st/ $2\le b_1\lsym<b_k\le N\]$, and
\vvn-.2>
$$
\alignat2
& \hp{{}\x\>\Bigl(\psit(t^1_1}
\Uxi(u\);t)\,\simeq{}
\Tagg{UxikQ}
\\
\nn5>
& \aligned
\!{}\simeq\,B\41(t^1_1)\ldots\)B\4r(t^1_{\xi^1}\])\,
\Bigl(\psit(\txine\])\bigl(\)
\Ue\+{N-1}_{\)\bar\xi\],\)k-1,\Qbar}(u-1\);\tbar\))\bigr)\;Q_1\>A(u)\,
\prod_{i=1}^{\xi^1\!}\,{u-t^1_i\]-1\over u-t^1_i}\>+{}
\\
\nn6>
{}+\,\psit(\txine\])\bigl(\)
\Ue\+{N-1}_{\)\bar\xi\],\)k,\Qbar}(u\);\tbar\))\bigr)\}\Bigr)\> & \){}+{}
\endaligned
\\
\nn7>
& {}\)+\,{1\over(\xi^1\}-1)\)!}\,\Sym_{\,\txine}\biggl[\;
{1\over u-t^1_1}\,\Be\41_Q(u)\>B\42(t^1_2)\ldots\)B\4r(t^1_{\xi^1}\])\,\x{}
\\
\nn8>
& \hp{{}\)+{}\}}\x\>\biggl(
\psit(\txine\])\bigl(\)\BB\*{N-1}_{\)\bar\xi}(\)\tbar)\bigr)\,\x{}
\\
\nn3>
& \hp{{}\)+{}\x{}\!}\x\>\biggl(Q_1\,T_{11}(t^1_1)\,
\prod_{i=2}^{\xi^1\!}\,{t^1_1\]-t^1_i\]-1\over t^1_1\]-t^1_i}\;-
\,Q_2\,T_{22}(t^1_1)\,
\prod_{i=2}^{\xi^1\!}\,{t^1_1\]-t^1_i\]+1\over t^1_1\]-t^1_i}\;
\prod_{j=1}^{\xi^2\!}\,{t^1_1\]-t^2_j\]-1\over t^1_1\]-t^2_j}\;
\biggr)\>-{}\] &&
\\
\nn5>
&& \Llap{{}-\>\res_{v\)=\)t^1_1\!}\Bigl(\psit(\txine\])
\bigl(\)\Ue\+{N-1}_{\)\bar\xi\],\)1,\Qbar}(v\);\tbar\))\bigr)\}\Bigr)\!
\biggr)\biggr]} & \Rlap{\,.}
\endalignat
$$
\par
The last formula determines $\Uxi(u\);t)$ modulo
\vvn-.1>
terms with coefficients belonging to $\Ynp$. By the induction assumption
the expressions $\Ue\+{N-1}_{\)\bar\xi\],\)j,\Qbar}(v\);\tbar\))$ belong
to $I\+{N-1}_{\)\bar\xi,\)j,\Qbar}\bra\)v\);\tbar\)\ket$. Then together with
Lemmas~\[Rwk], \[Taa0], \[DlTaa], \[psiY], \[phipsi], and formulae~\(Ie)\),
\(pix)\), \(psih), the induction assumption implies that there exists
an expression $\Uxi(u\);t)$ which belongs to $\Ixik\]$ and obeys formula
\(UxikQ)\). Therefore, formula~\(TeB2) coincides with formula~\(TeB).
Theorem~\[main] is proved.
\qed
\Pf of Proposition~\[Bsing].
We prove the statement by induction \wrt/ $N\}$ using Proposition~\[BB1].
It is straightforward to check that the induction assumption yields
\vvn.2>
$$
e_{ab}\,\BB^{\vox}_{\)\xi}(\tti\);z)\,=\,0\,,\qqq\Rlap{2\le a<b\le N\).}
$$
Let $C(u)$ be the left bottom block in the decomposition \(ABCD) of $T(u)$.
Consider the coefficients of the series $C(u)$,
${C(u)=\]\sum_{s=1}^\8C_s\)u^{-s}}\}$. Then $C_1\}$ is a matrix from
$\HY(\C,W)$ whose entries are $e_{12}\lc e_{1N}$. Therefore, it remains
to check that $C_1\>\BB^{\vox}_{\)\xi}(\tti\);z)\)=\)0$.
\vsk.1>
Relations \(Tabcd) implies that
$$
C_1\)B(u)\>-\>B(u)\>C_1\)=\,D(u)\)-\)A(u)\,.
$$
Then, similarly to the proof of Proposition~\[ADB], we get
\vvu0>
$$
\align
C_1\>B\41(u_1)\ldots\)B\4r(u_r)\,=\,B\41(u_1)\ldots\)B\4r(u_r)\,C_1\)+{} &
\\
\nn4>
{}+\,{1\over(r-1)\)!}\,
\Symr_{\)\ur}\"{\lr}\biggl(B\42(u_2)\ldots\)B\4r(u_r)\,\x{} &
\\
\nn3>
{}\x\,\Bigl(D\41(u_1)\>\Rbar\"{1\)r}(u_1\]-u_r)\ldots\)
\Rbar\"{12}(u_1\]-u_2)\>-{} &\>
\prod_{i=2}^r\,{u_1\]-u_i\]-1\over u_1-u_i}\;A(u_1)\Bigr)\}\biggr)\,.
\\
\cnn-.2>
\endalign
$$
The rest of the proof of Proposition~\[Bsing] is similar to considerations
in the proof of Theorem~\[main] after Lemma~\[DRR] and in the proof of
Theorem~\[main2].
\epf

\Sect[proof2]{Proof of Theorem~\[main3]}
We will prove Theorem~\[main3] by induction \wrt/ $N\}$, suitably adapting
the proof of Theorem~\[main].
\vsk.1>
We will continue to use notation introduced in the previous section.
The maps $\psi$, $\psi(\xr)$ and $\psit(\xr)$ defined there for the Yangian
$\Ynn$ have their counterparts for the current Lie algebra $\glnnx$.
Abusing notation we will denote them by the same letters.
\vsk.1>
The embedding ${\psi\):\)\Unnx\>\to\>\Unx}$ is defined by the rule
$\psi\bigl(L_{ab}\*{N-1}(u)\bigr)\)=\)L_{a+1,\)b+1}(u)$, \>${a\),\bn-1}$,
where $L_{ab}\*{N-1}(u)$ is series \(Lab) for $\glnnx$. The \hom/
${\psi(\xr)\):\)\Unnx\>\to\>\End\bigl(W^{\ox r}\bigr)\ox\Unx}$ is given by
$$
\psi(\xr)\bigl(L_{ab}\*{N-1}(u)\bigr)\>=\>\one^{\ox r}\!\ox L_{a+1,\)b+1}(u)
\>+\>\sum_{i=1}^r\,{\one^{\ox(i\)-1)}\]\ox E_{ba}\+{N-1}\]\ox\one^{\ox(r-i)}
\over u-x_i}\>\ox\one\,,
\vv-.2>
$$
${a\),\bn-1}$. The map $\psit(\xr)\):\)\Unnx\>\to\>W^{\ox r}\!\ox\Unx$ is \st/
$$
\psit(\xr)\bigl(L_{ab}\*{N-1}(u)\bigr)\>=\>
\wb_1^{\)\ox r}\!\ox L_{a+1,\)b+1}(u)\>+\>\dl_{a1}\>\sum_{i=1}^r\,
{\wb_1^{\)\ox(i-1)}\!\ox\wb_b\}\ox\wb_1^{\)\ox(r-i)}\}\over u-x_i}\>\ox\one\,,
\vv-.3>
$$
${a\),\bn-1}$, and ${\psit(\xr)(X_1X_2)\)=\)\psi(\xr)(X_1)\,\psit(\xr)(X_2)}$
for any ${X_1,X_2\}\in\Unnx}$.
\Lm{psiU}
We have
\vvn-.1>
$$
\psit(\xr)\)\bigl(\Unnx\>\ngpxo\bigr)\)\sub\)W^{\ox r}\!\ox\Unnx\>\ngpxo\,.
\vv-.3>
$$
\endpro
\nt
The proof is straightforward.
\vsk.5>
Define the embedding ${\pho\):\)U(\gl_{N-2}\lsb\)x\)\rsb)\>\hto\>\Unnx}$,
\,${\pho\bigl(L_{ab}\+{N-2}(u)\bigr)\)=\>L_{a+1,\)b+1}\+{N-1}(u)}$,
\;$a\),\bn-2$, similar to the embedding $\psi$.
\Lm{phipsi2}
We have \,$\psit(\xr)\o\pho\,=\,\wb_1^{\)\ox r}\!\]\ox(\psi\o\pho)$.
\endpro
\nt
The proof is straightforward.
\vsk.5>
For any element ${g\in W^{\ox r}\!\ox\Unnx}$ we consider its components
$g^{\)a_1\lc a_r}\}$ given by \uline
$g\>=\!\sum_{a_1\lc a_r=1}^{N-1}\!\wb_{a_1}\!\lox\)\wb_{a_r}\!\ox
g^{\)a_1\]\lc a_r}$.
\vsk.3>
Recall that \>$\bar\xi=\)(\xi^2\]\lc\xi^{\)N-1})$ and \,$\tbar\,=\,
(t^2_1\lc t^2_{\xi^2};\;\ldots\;;t^{\)N-1}_1\}\lc t^{\)N-1}_{\xi^{N-1}})$.
\Prop{FF1}
We have
\vvn-.3>
$$
\FF_{\]\xi}(t)\,=\!\sum_{a_1\]\lc a_{\xi^1}=\)1}^{N-1}\!
L_{1,a_1+1}(t^1_1)\ldots\>L_{1,a_{\xi^1}+1}(t^1_{\xi^1}\])\,
\Bigl(\psit(\txine\])\bigl(\)\FF\*{N-1}_{\]\bar\xi}(\)\tbar)\bigr)\}\Bigr)
^{\}a_1\]\lc a_{\xi^1}}.\kern-1em
\vv->
\Tag{FFe}
$$
\endpro
\Pf.
The statement follows from formula \(FF)\).
\epf
In Section~\[:filt] we introduced the degree \fn/ $\deg$ on the algebra
$\Yn\udzt$ by the rule ${\deg\)u\1\}=\)\deg\)\der_u\}=\)\deg\)\zt=-\)1}$.
We define the degree to \rat/ expressions in $\txn\}$ \wcoeff/ $\Yn\udzt$,
assuming in addition that \>${\deg\)t^a_i\]=\)1}$ and
\>${\deg\)(t^a_i\}-t^b_j)\1\]=\)-\)1}$ for all possible $a\),b\),i\),j$.
We use similar definition of the degree for the Yangian $\Ynn$ and related
\vvn-.16>
\rat/ expressions in ${t^2_2\lc t^{\)N-1}_{\xi^{N-1}}}\}$. Notice that the maps
$\psi$, $\psi(\txine)$ and $\psit(\txine)$ agree with the introduced degree
\vvn.1>
\fn/s and the projections to the graded algebras. For example,
if ${X\]\in\Ynn}$ and ${\deg\)X=\)k}$, then ${\deg\)\psi(X)=\)k}$ and
${\psi\bigl(\grad_{\)k}\]X\bigr)=\grad_{\)k}\bigl(\psi(X)\bigr)}$.
\Prop{gradB}
We have \>$\deg\)\bigl(\BB_{\)\xi}(t)\bigr)\)=\)-\>|\)\xi\)|$, and
\>$\grad_{-|\xi|}\bigl(\BB_{\)\xi}(t)\bigr)\)=\)\FF_{\]\xi}(t)$.
\endpro
\Pf.
The statement follows from formula~\(TL) and Propositions~\[BB1], \[FF1]
by induction \wrt/ $N\}$.
\epf
\Pf of Theorem~\[main3].
As in Section~\[:filt], set ${Q=1+\zt K}$. We will employ Theorems~\[gradY],
\[SkQ], and Proposition~\[gradB], and adapt arguments used in the proof of
Theorem~\[main] to prove the statement of Theorem~\[main3] by induction
\wrt/ $N\}$.
\vsk.1>
Consider formula~\(TeX) and write there $\Be_{k,Q}(u)$ instead of $\Be_Q(u)$
to indicate the dependence on $k$. Recall that we have $r=\)\xi^1\}$,
\vvn.1>
${u_i=\)t^1_i}$, $i=1\lc\xi^1\}$, and ${A(u)=T_{11}(u)}$. Taking the sum
\vv-.2>
of equalities \(TeX) for \)${k=1\lc l}$, \)with coefficients
${\dsize{(-1)^{l-k}(N\]-k)\)!\over(N\]-l\))\)!\,(l-k\))\)!}}$
and using formula~\(ST)\), for any ${X\]\in\Ynn}$, we obtain
\vvu0>
$$
\align
& \hp{{}\!\!\simeq\,\biggl\lb B\41}\,\Se_{\)l,Q}(u)\,
B\41(t^1_1)\ldots\)B\4{\xi^1\}}(t^1_{\xi^1}\])\;\psi(\txine\])(X)\,\simeq{}
\Tagg{SeX}
\\
\nn6>
& \alignedat2
{}\simeq\,\biggl\lb {} & B\41(t^1_1)\ldots\)B\4{\xi^1\}}(t^1_{\xi^1}\])\,\x{}
\\
\nn4>
&{}\!\x\biggl(\psi(\txine\])\)\bigl(\)\Se\+{N-1}_{\)l-1,\Qbar}(u-1)\,X\bigr)\;
\Bigl(Q_1\>T_{11}(u)\,\prod_{i=1}^{\xi^1\!}\,{u-t^1_i\]-1\over u-t^1_i}\,-\)1
\Bigr)\,+{} &&
\\
\nn6>
&& \Llap{{}+\,\psi(\txine\])\)\bigll(\]\bigl(\)\Se\+{N-1}_{\)l-1,\Qbar}(u-1)-
\Se\+{N-1}_{\)l-1,\Qbar}(u)+\Se\+{N-1}_{\)l,\Qbar}(u)\bigr)\>X\bigrr)\}\biggr)
\)}& \){}+{}
\endalignedat
\\
\nn9>
& \hp{{}\simeq\}{}}
\aligned
{}+\,{1\over(\xi^1\}-1)\)!}\,\Symr_{\)\txine}\"{1\ldots\)\xi^1\]}{} &
\biggl[\;{1\over u-t^1_1}\;\Pe\41_{l,Q}(u)\>B\42(t^1_2)\ldots\)B\4{\xi^1\}}
(t^1_{\xi^1}\])\>\x{}
\\
\nn4>
{}\x\){} & \)\biggl(\psi(\txine\])(X)\>\Bigl(Q_1\,T_{11}(t^1_1)\,
\prod_{i=2}^{\xi^1\!}\,{t^1_1\]-t^1_i\]-1\over t^1_1\]-t^1_i}\,-\)1\Bigr)\>+{}
\\
\nn4>
{}+\){} & \>\Bigl(1-\]\res_{v\)=\)t^1_1\!}
\bigll(\psi(\txine\])\bigl(\)\Te\+{N-1}_{1,\Qbar}(v)\bigr)\]\bigrr)\}\Bigr)\,
\psi(\txine\])(X)\}\biggr)\biggr]\biggr\rb\,,
\endaligned
\\
\cnn-.6>
\endalign
$$
where
\vvn-.8>
$$
\Pe_{l,Q}(u)\>=\>{1\over(N\]-l\))\)!}\,
\sum_{k=1}^l\,(-1)^{l-k}\,{(N\]-k)\)!\over(l-k)\)!}\ \Be_{k,Q}(u)\,.
\vv-.4>
\Tag{Pel}
$$
Notice that in formula \(TeX) for ${k=N}\}$ we have
${\Te\+{N-1}_{\]N\},\Qbar}(u)\)=\)0}$ and ${\Be_{N\},Q}(u)\)=\)0}$, and in
formula \(SeX) for ${l=N}\}$ we assume that ${\Se\+{N-1}_{N,\Qbar}(u)\)=\)0}$.
We remind that ${\Te_{0,Q}(u)=1}$ and
${\Te\+{N-1}_{0,\Qbar}(u)\)=\)\Se\+{N-1}_{0,\Qbar}(u)\)=\)1}$
by convention.
\vsk.2>
Recall that $\Be_{k,Q}(u)$ is an ${(N-1)\]\x\]1}$ matrix whose entries are
series in $u\1\!$ \wcoeff/ $\Yn$, and $\Pe_{l,Q}(u)$ is an ${(N-1)\]\x\]1}$
matrix of the same kind,
\vvn.1>
$$
\Be_{k,Q}(u)\)=\)\bigl(\>\Be_{k,Q\);\)1}(u)\)\lc\)\Be_{k,Q\);\)N-1}(u)\bigr)\,,
\kern1.6em
\Pe_{l,Q}(u)\)=\)\bigl(\>\Pe_{l,Q\);\)1}(u)\)\lc\)\Pe_{l,Q\);\)N-1}(u)\bigr)\,.
\vv.1>
$$
\Lm{Bew}
For any ${m=k\lc N}\}$, \)distinct ${i_1\lc i_k\]\in\lb\)1\lc m\)\rb}$,
and any $s=1\lc k$, we have
\vvn-.5>
$$
\align
& \Be_{k,Q\);\)a}(u)\,=\,{m\)!\,(N\]-m)\)!\over(k-1)\)!\,(N\]-k\))\)!}\;\x{}
\\
\nn6>
& \quad{}\x\,{}({\tr_{V^{\]\ox m}}\:\}}\ox\id\))\bigl(E_{a+1,1}\"{i_s}\,
Q\"{i_1}\!\]\ldots\>Q\"{i_k}\>\>T\"{i_1,\)m+1}(u)\>\ldots\>
T\"{i_k,\)m+1}(u-k+1)\;\AA\"{1\ldots m}\6m\)\bigr)\,.\!
\\
\cnn-.4>
\endalign
$$
\endpro
\Pf.
Formulae \(BeQ)\), \(Qbar)\), \(ABCD)\), and \(Sxw) give that
$$
\Be_{k,Q\);\)a}(u)\,=\,
k\,({\tr_{V^{\]\ox k}}\:\}}\ox\id\))\bigl(E_{a+1,1}\"1\,
Q\"1\!\]\ldots\>Q\"k\>\>T\"{1,\)k+1}(u)\>\ldots\>
T\"{k,\)k+1}(u-k+1)\;\AA\"{1\ldots k}\6k\)\bigr)\,.
$$
The remaining consideration is similar to the proof of Lemma~\[Tetr].
\uugood
\epf
Lemma~\[Bew] and formula~\(Pel) allow us to obtain $\Be_{k,Q\);\)a}(u)$
and $\Pe_{l,Q\);\)a}(u)$ as coefficients of a formal \dif/ operator
similarly to formulae~\(Dgm)\), \(DS)\),
$$
\align
& \sum_{k=1}^m\,(-1)^{(k-1)}\,{(N\]-k)\)!\over(m-k)\)!}
\ \Be_{k,Q\);\)a}(u)\,e^{(m-k)\)\der_u}\)={}
\Tagg{BP}
\\
\nn5>
& \,{}=\,\sum_{l=1}^m\,(-1)^{(l-1)}\,{(N\]-l\))\)!\over(m-l\))\)!}
\ \Pe_{l,Q\);\)a}(u)\,(e^{\)\der_u}\}-1)^{m-l}\)={}
\\
\nn6>
& \,\,{}=\, m\,(N\]-m)\)!\;({\tr_{V^{\]\ox m}}\:\}}\ox\id\))
\Bigl(E_{a+1,1}\"1\>Q\"1\)T\"{1,\)m+1}(u)\,e^{-\der_u}\>\x{}
\\
\nn6>
& \hp{{}={}}\,\)\x\,\bigl(1-Q\"2\)T\"{2,\)m+1}(u)\,e^{-\der_u}\bigr)\)\ldots\)
\bigl(1-Q\"mT\"{m,\)m+1}(u)\,e^{-\der_u}\bigr)\;
\AA\"{1\ldots m}\6m\)\Bigr)\,e^{\)m\der_u}\).\kern-3.8em
\\
\cnn-.2>
\endalign
$$
Introduce the following ${(N-1)\]\x\]1}$ matrices
$$
\Fe_{\)l,K}(u)\>=\>\bigl(\>\Fe_{\)l,K\];\)1}(u)\)\lc\)
\Fe_{\)l,K\];\)N-1}(u)\bigr)\,,\qqq\Rlap{l=1\lc N\),}
$$
whose entries are defined by the rule
\vvn-.2>
$$
\align
& \){1\over N}\,
\sum_{l=1}^N\,(-1)^{(l-1)}\,\Fe_{\)l,K\];\)a}(u)\,\der_u^{\>N\]-l}\)={}
\Tagg{Fe}
\\
\nn3>
{}=\,(\trVN & {}\ox\id\))\Bigl(E_{a+1,1}\"1\>
\bigl(\der_u\]-K\"2\}-L\"{2,\)N+1}(u)\bigr)\)\ldots\)
\bigl(\der_u\]-K\"N\}-L\"{N\],\)N+1}(u)\bigr)\;\AA\"{1\ldots N}\6N\)\Bigr)\,.\!
\endalign
$$
It is easy to see that $\Fe_{\>1,K}(u)=0$. Like in the proof of Theorem~\[SkQ],
formulae~\(BP) and~\(Fe) imply that
${\deg\)\bigl(\Pe_{l,Q}(u)\bigr)\)=\)1-l}\)$ for any $l=1\lc N\}$, and
$$
\grad_{1-\)l}\]\bigl(\Pe_{l,Q}(u)\bigr)\,=\,\Fe_{\)l,\)\zt\]K}(u)\,.
$$
\vsk-.5>
\par
Set \>$\Kbar=\]\sum_{a=1}^{N-1}\]K_{a+1}\>E\+{N-1}_{aa}\]$,
\,so that $\Qbar=1+\zt\Kbar$.
\vsk-.2>
\vsk0>
\Lm{T1S2}
We have \,$\deg\)\Bigl(1-\]\res_{v\)=\)t^1_1\!}\bigll(\psi(\txine\])
\bigl(\)\Te\+{N-1}_{1,\Qbar}(v)\bigr)\]\bigrr)\}\Bigr)=-\)1$ \>and
\vvn-.2>
$$
\align
& \grad_{-1}\]\Bigl(1-\]\res_{v\)=\)t^1_1\!}\bigll(
\psi(\txine\])\bigl(\)\Te\+{N-1}_{1,\Qbar}(v)\bigr)\]\bigrr)\}\Bigr)\>={}
\\
\nn5>
& \hp{\roman{g}\]}=\,\Bigl(\){1\over v-t^1_1}\>-
\psi(\txine\])\bigl(\Gc\+{N-1}_{1,\)\zt\]\Kbar}(v)\bigr)\]\Bigr)
\Big|_{v\)=\)t^1_1}\]-\)\res_{v\)=\)t^1_1\!}\bigll(\psi(\txine\])
\bigl(\Gc\+{N-1}_{2,\)\zt\]\Kbar}(v)\bigr)\]\bigrr)\,.
\\
\cnn-.5>
\endalign
$$
\endpro
\Pf.
The statement follows from formulae~\(Dres) and~\(Gc12) by straightforward
calculations. Notice that the denominator of the expression
${(v-t^1_1)\1\}-\psi(\txine\])\bigl(\Gc\+{N-1}_{1,\Kbar}(v)\bigr)}$ does not
actually contain the factor ${v-t^1_1}$, which makes the substitution
${v=t^1_1}$ in the expression legal.
\epf
The Poincar\'e\>-Birkhoff\>-Witt theorem for the Yangian $\Yn$,
see~\cite{MNO, Corollary~1.23}\), \cite{Mo, Theorem~2.6}\), yields that
if ${X\]\in\Ynp}$, ${\deg\)X\]=d}$, then ${\grad_d\]X\]\in\Unx\>\ngpx}$.
\vsk.3>
Let ${F(u)=\bigl(L_{12}(u)\)\lc L_{1N}(u)\bigr)=\)\grad_{-1}\bigl(B(u)\bigr)}$,
\vvn.1>
and ${d=\deg X}$ in formula~\(SeX)\). We apply the map $\grad_{d-l}\}$
\vvn.06>
to both sides of \(SeX)\), and then evaluate the result at ${\zt=1}$.
Recall that we have
$$
\grad_{-l}\bigl(\Se_{\)l,Q}(u)\bigr)=\)\Gc_{\)l,\)\zt\]K}(u)\,,\qqq
\grad_{-s}\bigl(\Se\+{N-1}_{\)s,\Qbar}(u)\bigr)=\)
\Gc\+{N-1}_{\)s,\)\zt\]K}(u)\,.
$$
Thus we obtain that for any ${Z\]\in\Unnx}$,
\vvn-.2>
$$
\alignat2
& \hp{{}\!\!\simeq\,\biggl\lb F\41}\,\Gc_{\)l,K}(u)\,
F\41(t^1_1)\ldots\)F\4{\xi^1\}}(t^1_{\xi^1}\])\;\psi(\txine\])(Z)\,\simeq{}
\Tagg{GcZ}
\\
\nn6>
& \aligned
{}\simeq\,\biggl\lb F\41(t^1_1)\ldots\)F\4{\xi^1\}}(t^1_{\xi^1}\])\,
\biggl(\psi(\txine\])\)\bigl(\)\Gc\+{N-1}_{\)l-1,\Kbar}(u)\,Z\bigr)\,
\Bigl(K_1\]+L_{11}(u)\)-\]\sum_{i=1}^{\xi^1\!}\,{1\over u-t^1_i}\)\Bigr)\>+{} &
\\
\nn6>
{}+\,\psi(\txine\])\)\bigll(\]\bigl(\Gc\+{N-1}_{\)l,\Kbar}(u)
-\)\dot\Gc\+{N-1}_{\)l-1,\Kbar}(u)\bigr)\>Z\bigrr)\}\biggr)\)& \)+{}
\endaligned &&
\\
\nn7>
& \hp{{}\simeq\}{}}
\aligned
{}+\,{1\over(\xi^1\}-1)\)!}\,\Sym_{\)\txine}{} & \biggl[\;
{1\over u-t^1_1}\;\Fe\41_{l,K}(u)\>F\42(t^1_2)\ldots\)F\4{\xi^1\}}
(t^1_{\xi^1}\])\>\x{}
\\
\nn4>
{}\x\){} & \)\biggl(\psi(\txine\])(Z)\>\Bigl(K_1\]+L_{11}(t^1_1)\)-\]
\sum_{i=2}^{\xi^1\!}\,{1\over t^1_1\]-t^1_i}\)\Bigr)\,+{}
\endaligned
\\
\nn6>
&& \Llap{{}+\,\biggl(\!\Bigl(\){1\over v-t^1_1}\>-
\psi(\txine\])\bigl(\Gc\+{N-1}_{1,\Kbar}(v)\bigr)\]\Bigr)
\Big|_{v\)=\)t^1_1}\]-\)\res_{v\)=\)t^1_1\!}\bigll(\psi(\txine\])
\bigl(\Gc\+{N-1}_{2,\Kbar}(v)\bigr)\]\bigrr)\!\biggr)\>\x{}} &
\\
\nn1>
&& \Llap{{}\x\,\psi(\txine\])(Z)\}\biggr)\biggr]\biggr\rb\)} & \>.\!
\\
\cnn-.9>
\endalignat
$$
Here \;$\dsize
\dot\Gc\+{N-1}_{\)l-1,\Kbar}(u)\>=\>{d\over du}\,\Gc\+{N-1}_{\)l-1,\Kbar}(u)$.
\vsk.5>
The rest of the proof of Theorem~\[main3] is similar to the proof of
Theorem~\[main] after formula~\(TeX)\). We apply both sides of formula
\(GcZ) to the vector ${\wb_1^{\)\ox r}}\!\]$, that amounts to replacing
there the map $\psi(\ur)$ by $\psit(\ur)$, and then take
${Z=\FF\*{N-1}_{\]\bar\xi}(\)\tbar)}$. Since formula~\(FFe) is \eqv/ to
$$
\FF_{\]\xi}(t)\,=\,F\41(t^1_1)\ldots\)F\4{\xi_1}(t^1_{\xi^1}\])
\ \psit(\txine\])\bigl(\)\FF\*{N-1}_{\]\bar\xi}(\)\tbar)\bigr)\,,
$$
\lhs/ of \(GcZ) becomes $\Gc_{\)l,K}(u)\,\FF_{\]\xi}(t)$, and to transform
\rhs/ to the required form, we employ the induction assumption.
\epf

\uupage
\Appendix

\Sect[naive]{The Gaudin model as a limit of the \XXX/-\)type model}
Let $\Mn$ be \gnmod/s. In Sections~\[:emods] and~\[:Sform], we describe
the \XXX/-\)type model on the tensor product ${\Mox}$ considering the action
of the Yangian $\Yn$ on the tensor product ${\Moxz}$ of \emod/s. The objects
associated with the \XXX/-\)type model are the operators
\vvn.2>
$$
T^{\Moxs}_{ab}(u\);z)\,=\,T_{ab}(u)\vst{\Moxz}\>,\qquad\Rlap{a\),\bn\),}
$$
\rat/ly depending on $u$ and $\zn$;
the \tmcs/ $\Te^{\Moxs}_{k,Q}(u\);z)$, cf.~\(Tez)\); the modified \tmcs/
\vvn.2>
$$
\align
\Se^{\Moxs}_{k,Q}(u\);z)\, &{}=\,\Se_{k,Q}(u)\vst{\Moxz}\>={}
\\
\nn5>
&{}=\,{1\over(N\]-k)\)!}\,\sum_{l=0}^k\,
(-1)^{k-l}\,{(N\]-l\))\)!\over(k-l\))\)!}\ \Te^{\Moxs}_{l,Q}(u)\,,
\endalign
$$
$k=0\lc N\}$, \;cf.~\(ST)\); the \dif/ operator
$$
\align
\Dg^{\Moxs}_{N,Q}(u\);\der_u\);z)\, &{}=\,
\sum_{k=0}^N\,(-1)^k\,\Te^{\Moxs}_{k,Q}(u\);z)\,e^{-k\der_u}\)={}
\Tagg{DgNTz}
\\
\nn4>
& {}=\,\sum_{k=0}^N\,(-1)^k\,\Se^{\Moxs}_{k,Q}(u\);z)\,
(e^{\)\der_u}\}-1)^{(N\]-k)\der_u}\)e^{-N\der_u}\);\kern-1em
\endalign
$$
the deformed Shapovalov form $S_{\Moxs}^{\>z}$ on ${\Mox}$, see~\(SMz)\);
the universal \wt/ \fn/ $\BB^{\vox}_{\)\xi}(t\);z)$, cf.~\(BBtz)\); the \raf/s
$$
\align
& \Qxe^{\>a,i}(t\);z\);\La)\,={}
\\
\nn4>
& \!{}=\,{Q_a\over Q_{a+1}}\,
\prod_{j=1}^n\>{t^a_i\]-z_j+\La^a_j\over t^a_i\]-z_j+\La^{a+1}_j}\,
\prod_{j=1}^{\xi^{a-1}\!\!}\;{t^a_i\]-t^{a-1}_j\}+1\over t^a_i\]-t^{a-1}_j}\,
\prod_{\tsize{j=1\atop j\ne i}}^{\xi^a}\,
{t^a_i\]-t^a_j\]-1\over t^a_i\]-t^a_j\]+1}\,
\prod_{j=1}^{\xi^{a+1}\!\!}\;{t^a_i\]-t^{a+1}_j\over t^a_i\]-t^{a+1}_j\}-1}\;,
\kern-2em
\\
\cnn-.4>
\endalign
$$
$\an-1$, \>$i=1\lc\xi^a\}$, \>that define the \Bae/s \(Bae)\);
and the \dif/ operator
\vvn.3>
$$
\Mg_{\>\xi\],Q}(u\),\]\der_u\);t\);z\);\La)\,=\,
\bigl(1-\Xxi{\)1}(u\);t\);z\);\La)\,e^{-\der_u}\bigr)\)\ldots\)
\bigl(1-\Xxi N(u\);t\);z\);\La)\,e^{-\der_u}\bigr)\,,
\vv.3>
$$
cf.~\(1-Xz)\), \(Xxiz)\), whose coefficients are \eva/s of the \tmcs/.
\vsk.2>
In Section~\[:cemods], we describe the Gaudin model on the tensor product
${\Mox}$ considering the action of the current algebra $\glnx$ on the tensor
product ${\Moxz}$ of \emod/s over $\glnx$. For the Gaudin model, the associated
objects are the operators
\vvn.2>
$$
L^{\Moxs}_{ab}(u\);z)\,=\,\sum_{i=1}^n\){e_{ba}\"i\over u-z_i}\;,
\qqq\Rlap{\an\),}
$$
\rat/ly depending on $u$ and $\zn$; the \tmcs/ $\Gc^{\)\Moxs}_{k,K}(u\);z)$,
cf.~\(Labz)\), \(Gcz)\); the \difl/ operator
\vvn-.2>
$$
\De^{\Moxs}_{\]K}\](u\),\der_u\);z)\,=\,
\sum_{k=0}^N\,(-1)^k\,\Gc^{\Moxs}_{\)k,K}(u\);z)\,\der_u^{\>N\]-k}
\Tag{detLz}
$$
with \$\End(\Mox)\)$-valued \rat/ coefficients, cf.~\(detL)\);
the Shapovalov form $\SMox\}$ on ${\Mox}$; the universal \wt/ \fn/
$\FF^{\vox}_{\)\xi}(t\);z)$, cf.~\(FFtz)\); the \raf/s
\vvn.3>
$$
\align
& \Kxe^{\>a,i}(t\);z\);\La)\,={}
\\
\nn4>
& \!{}=\,K_a\]-K_{a+1}\]+
\sum_{j=1}^n\,{\La^a_j\]-\La^{a+1}_j\over t^a_i\]-z_j}\;+\,
\sum_{j=1}^{\xi^{a-1}\!\!}\;{1\over t^a_i\]-t^{a-1}_j}\;-\,
2\)\sum_{\tsize{j=1\atop j\ne i}}^{\xi^a}\,{1\over t^a_i\]-t^a_j}\;+\,
\sum_{j=1}^{\xi^{a+1}\!\!}{1\over t^a_i\]-t^{a+1}_j}\;,
\\
\cnn-.3>
\endalign
$$
$\an-1$, \>$i=1\lc\xi^a\}$, \>that define the \Bae/s \(Bae2)\);
and the \difl/ operator
\vvn-.5>
$$
\Mxi(u\),\der_u\);t\);z\);\La)\,=
\prodr_{1\le a\le N\!}\biggl(\der_u-K_a-\>\sum_{i=1}^n\){\La^a_i\over u-z_i}
\;-\>\sum_{i=1}^{\xi^{a-1}\!\!}\,{1\over u-t^{a-1}_i\}}\;+\>
\sum_{i=1}^{\xi^a}\,{1\over u-t^a_i}\biggr)\,.\kern-.4em
$$
cf.~\(MxiL)\), whose coefficients are \eva/s of the \tmcs/.
\vsk.2>
The objects associated with the Gaudin model can be obtained from
the corresponding objects for the \XXX/-\)type model in the following limit.
\vsk.2>
Set $\eps\1\]z=(\eps\1\]z_1\lc\eps\1\]z_n)$ and
\>$\eps\1t=(\eps\1t^1_1\lc\eps\1t^{\)N-1}_{\xi^{N-1}})$.
\Th{limit}
Let $Q=1+\eps K$. As $\eps\to 0$, we have
\vvn.4>
$$
T^{\Moxs}_{ab}(\eps\1u\);\eps\1\]z)\,=\,\dl_{ab}+\)
\eps\)L^{\Moxs}_{ab}(u\);z)\)+\)O(\eps^2)\,,
\vv.2>
\Tag{TLlim}
$$
$a\),\bn\}$,
\vvn.2>
$$
\Se^{\Moxs}_{k,Q}(\eps\1u\);\eps\1\]z)\bigr)\,=\,
\eps^{\)k}\>\Gc^{\Moxs}_{k,K}(u\);z)\)+\)O(\eps^{k+1})\,,
\vv.2>
\Tag{SGlim}
$$
$\koN\}$,
\vvn.1>
$$
\gather
\Dg^{\Moxs}_{N,Q}(\eps\1u\),\eps\)\der_u\);\eps\1\]z)\,=\,
\eps^{\)N}\>\De^{\Moxs}_{\]K}\](u\);\der_u\);z)\)+\)O(\eps^{\)N\]+1})\,,
\Tag{Dlim}
\\
\nn10>
S_{\Moxs}^{\>\eps^{\}-1}\}z}\}\to\,\SMox\),
\Tag{Slim}
\endgather
$$
$$
\align
\BB^{\vox}_{\)\xi}(\eps\1t\);\eps\1\]z)\,
\prod_{i=1}^n\,\prod_{a=1}^{N-1}\,\prod_{j=1}^{\xi^a\!}\,{\eps\over t^a_j-z_i}
\,\prod_{a=1}^{N-2}\,\prod_{i=1}^{\xi^a\!}\;
\prod_{j=1}^{\xi^{a+1}\!\!}\;{\eps\over t^{a+1}_j\}-t^a_i}\;={}\kern-2em &
\Tag{BFlim}
\\
\nn7>
{}=\,\eps^{\)|\xi|}\>\FF^{\vox}_{\)\xi}(t\);z)\)+\)O(\eps^{\)|\xi|+1})\,,\!\]
\kern-2em &
\endalign
$$
$$
\Qxe^{\>a,i}(\eps\1t\);\eps\1z\);\La)\,=\,
1\)+\)\eps\>\Kxe^{\>a,i}(t\);z\);\La)\)+\)O(\eps^2)\,,
\Tag{QKlim}
$$
$$
\Mg_{\>\xi\],Q}(\eps\1u\),\eps\)\der_u\);\eps\1t\);\eps\1\]z\);\La)\,=\,
\eps^N\)\Mxi(u\),\der_u\);t\);z\);\La)\,.
\Tag{Mlim}
$$
\endpro
\Pf.
Formulae~\(TLlim) and~\(QKlim) are straightforward. The proof of
formulae~\(SGlim) and~\(Dlim) is similar to the proof of Theorem~\[SkQ].
Formula~\(Slim) follows from formula~\(R8)\). Formula~\(BFlim) is similar
to the claim of Proposition~\[gradB].
\epf
\Rem
Using Theorem~\[limit] one can see that Theorem~\[Shapz] implies
Theorem~\[Shap], and Theorem~\[main2] implies Theorem~\[main4]
for an isolated \sol/ $\tti=(\ttxn)$ of the \Bae/s \(Bae2).
\enddemo

\Sect[dyn]{Dynamical Hamiltonians as parts of \tmcs/}
Let \>$C_1=\sum_{a=1}^Ne_{aa}$, \>and
\>$C_2=\!\!\sum_{1\le a<b\le N}\!\}(e_{aa}e_{bb}\]-e_{ab}e_{ba}+e_{aa})$.
They are central elements of the universal enveloping algebra $\Un$.
\par
Let $\Mn$ be \gnmod/s. Let $K$ be a diagonal matrix,
$K=\]\sum_{a=1}^N\]K_aE_{aa}$. The Gaudin Hamiltonians
$H_{1,K}(z)\lc H_{n,K}(z)$, see~\cite{G}\), the \rat/ dynamical Hamiltonians
$G_{1,K}(z)\lc G_{N\],K}(z)$, see~\cite{FMTV}\), \cite{TV2}\), and the \tri/
dynamical Hamiltonians $X_{1,K}(z)\lc X_{N\],K}(z)$, see~\cite{TV2}\),
\cite{TV3}\), acting in ${\Mox}$, are the following operators:
\vvn->
$$
\gather
H_{i,K}(z)\,=\,\sum_{a=1}^N\,K_a\>e_{aa}\"i\>+\>
\sum_{a,b=1}^N\>\sum_{\tsize{j=1\atop j\ne i}}^n\,
{e_{ab}\"i\>e_{ba}\"j\over z_i\]-z_j}\;,
\\
\nn6>
G_{a,K}(z)\,=\,\sum_{i=1}^n\>z_i\>e\"i_{aa}\,+\>
\sum_{\tsize{b=1\atop b\ne a}}^N\>{e_{ab}\>e_{ba}\]-\)e_{aa}\over K_a\]-K_b}\;.
\\
\nn6>
X_{a,K}(z)\,=\,-\,{(e_{aa})^2\over2}\>+\,\sum_{i=1}^n\>z_i\>e\"i_{aa}\,+
\>\sum_{b=1}^N\,\sum_{1\le i<j\le n}\!e\"i_{ab}\>e\"j_{ba}\>+\>
\sum_{\tsize{b=1\atop b\ne a}}^N\>{K_b\over K_a\]-K_b}\;
(e_{ab}\>e_{ba}\]-\)e_{aa})\,.
\\
\cnn-.5>
\endgather
$$
Recall that $e_{ab}\}$ acts on ${\Mox}$ as $\sum_{i=1}^ne\"i_{ab}$.
\vsk.2>
It is well known that the Gaudin Hamiltonians come from the residues of
the \tmx/ $\Gc^{\Moxs}_{2,K}(u\);z)$:
\vvn.3>
$$
\align
& \Gc^{\Moxs}_{2,K}(u\);z)\,={}
\\
\nn5>
& \!{}=\,\tr K^{\wg2}\)+\>\sum_{i=1}^n\,{1\over u-z_i}\,\biggl[\>C\"i_1\>
\biggl(\]\tr K\)+\)\sum_{\tsize{j=1\atop j\ne i}}^n\,{C\"j_1\over z_i\]-z_j}
\biggr)\)-\)H_{i,K}(z)\>\biggr]\>+\>\sum_{i=1}^n\,{C\"i_2\over(u-z_i)^2}\;.
\endalign
$$
The \rat/ and \tri/ dynamical Hamiltonians also can be recovered from
\vvn.1>
the \tmcs/ $\Gc^{\Moxs}_{\)k,K}(u\);z)$ and $\Te^{\Moxs}_{k,K}(u\);z)$, \resp/.
To this end, consider the following \$\End(\Mox)\)$-valued \raf/ of $u,x,\zn$:
\vv-.2>
$$
\align
\Gch_K\:(u\),x\);z)\, &{}=\,
\sum_{k=0}^N\,(-1)^k\,\Gc^{\Moxs}_{\)k,K}(u\);z)\,x^{\)N\]-k}
\>\prod_{a=1}^N\>{1\over x-K_a}\;,
\\
\nn4>
\Teh_K\:(u\),x\);z)\, &{}=\,
\sum_{k=0}^N\,(-1)^k\,\Te^{\Moxs}_{k,K}(u\);z)\,x^{\)N\]-k}
\>\prod_{a=1}^N\>{1\over x-K_a}\;.
\\
\cnn.2>
\endalign
$$
\Prop{dynHam}
We have
\vvn-.3>
$$
\align
\Gch_K\:(u\),x\);z)\,=\,1\,&{}-\,u\1\,\sum_{a=1}^N\>{e_{aa}\over x-K_a}\,-{}
\\
\nn4>
& {}-\,u^{-2}\,\sum_{a=1}^N\>{1\over x-K_a}\>
\biggl(G_{a,K}(z)\>-\sum_{\tsize{b=1\atop b\ne a}}^N\>
{e_{aa}\>e_{bb}\over K_a\]-K_b}\,\biggr)\)+\,O(u^{-3})
\\
\cnn-1.1>
\endalign
$$
and
\vvn-.8>
$$
\align
\Teh_K\:(u\),x\);z)\, &{}=\,1\,-\,
u\1\,\sum_{a=1}^N\>{K_a\over x-K_a}\,e_{aa}\>-{}
\\
\nn4>
& \>{}-\,u^{-2}\,\sum_{a=1}^N\>{K_a\over x-K_a}\>\biggl(X_{a,K}(z)\>+\>
{(e_{aa})^2\over2}\>-\)\sum_{\tsize{b=1\atop b\ne a}}^N\>
{K_b \over K_a\]-K_b}\,e_{aa}\>e_{bb}\biggr)\)+\,O(u^{-3})\,.
\\
\cnn-.5>
\endalign
$$
\endpro
\Pf.
To get the first formula we compute the expansion of
\vvn.06>
$\Gc^{\Moxs}_{\)k,K}(u\);z)$ as $u\to\8$ using formulae~\(detL)\), \(DG)\),
\vvn.08>
\(Labz)\), \(Gcz)\), and Proposition~\[ADLA]. For the second formula we
compute the expansion of $\Te^{\Moxs}_{k,K}(u\);z)$ using formulae~\(Te)\),
\(Tez)\), and the definition of the module $\Moxz$ over the Yangian $\Yn$.
\epf

\Sect[real]{Eigenvalues of \tmcs/ for real $\zn$}
Denote by $\glnr\}$ the real part of the Lie algebra $\gln\}$ generated
by $e_{ab}$, $a\),\bn\}$. Let $\Mn$ be \fd/ \irr/ \gnmod/s with \hw/s $\Lan$,
and \hwv/s $v_1\lc v_n$. Denote by $\Mnr$ be their real parts generated over
$\glnr\}$ by the vectors $v_1\lc v_n$, \resp/. It is clear that for real
\vvn.06>
$u\),\)\zn$, real matrices $Q\),K$, and any ${\kN}\}$, the operators
\vvn-.16>
$\Te^{\Moxs}_{k,Q}(u\);z)$ and $\Gc^{\)\Moxs}_{k,K}(u\);z)$ preserve
the tensor product $\Moxr$.
\vsk.2>
Since $\Mn$ are \fd/ \irr/ \gnmod/s, the tensor product $\SMox\!$ of Shapovalov
\vvn.06>
forms restricted to $\Moxr\]$ is a positive definite bilinear form.
Hence, an operator, which is \sym/ \wrt/ this form, has real \eva/s.
Therefore, Theorems~\[Shap] and~\[main4] yield the following statement.
\Th{realZ}
\vvn.06>
Let $K$ be a real diagonal matrix, $\zn$ real numbers and $\tti=(\ttxn)$
a \sol/ of the \Bae/s \(Bae2)\). If \,${\FF^{\vox}_{\]\xi}(\tti\);z)\ne 0}$,
then for any ${\kN}\}$, the \fn/ $\Zxi{\>k}(u\);\tti\);z\);\La)$,
cf.~\(Zxiz)\), takes real values for real values of $u$.
\endpro
For a dominant \wt/ $\La=(\LaN)$, set
${\La'\}=\)\La^N\!\]+1-\max\>\lb\)a\vert \La^a\}>\]\La^N\rb}$.
\vsk.2>
For \fd/ \irr/ \gnmod/s $L\>,\]M\}$, the \Rm/ $R_{LM}\:(u)$ is a \raf/ of $u$.
Its poles and degeneracy points are known from the \rep/ theory of the Yangian
$\Yn$, see~\cite{Mo}\). In particular, the following proposition holds.
\Prop{singR}
The \Rm/ $R_{M_iM_j}\:\](u)$ is well\)-\)defined and \ndeg/ for real values
of $u$ \st/ \,${u>\]\La_i^1\}-\La_j'}$ or \,${u<\La_i'\}-\La_j^1\]}$.
Moreover, $R_{M_iM_j}\:\](u)$ has a simple pole at \,${u=\La_i'\}-\La_j^1}$
with the degenerate residue, and the value ${R_{M_iM_j}\:\](\La_i^1\}-\La_j')}$
is well\)-\)defined, non\)-\)zero, and degenerate.
\endpro
It is clear that for real $\zn$, the operator $R_{\Mnd}\:\}(z)$, cf.~\(RMox)\),
preserves the tensor product $\Moxr$.
\Prop{SMz+}
Let $\zn$ be real numbers \st/ \>${z_i\]-z_j>\]\La_i^1\}-\La_j'\]}$ for
any ${i\ne j}$, \>$i\),\jn$. Then the bilinear form $S_{\Moxs}^{\>z}\!$,
cf.~\(SMz)\), restricted to $\Moxr\]$ is positive definite.
\endpro
\Pf.
The statement follows from Proposition~\[singR] and formula~\(R8)\).
\epf
\Th{realX}
\vvn.06>
Let $Q$ be a real diagonal matrix, $\zn$ real numbers \st/ \uline
\,${z_i\]-z_j>\]\La_i^1\}-\La_j'}$ for any ${i\ne j}$, \>$i\),\jn$,
and \,${\tti\)=\)(\ttxn)}$ a \sol/ of the \Bae/s \(Bae)\).
If \,${\BB^{\vox}_{\)\xi}(\tti\);z)\ne 0}$, then for any ${\kN}\}$, the \fn/
\vvn.2>
$$
\sum_{1\le a_1\lsym<a_k\le N\!\!\!}\;\;
\prod_{r=1}^k\,\Xxi{\>a_r}(u-r+1\);\tti\);z\);\La)\,,
\vv.2>
$$
cf.~\(Xxiz)\), takes real values for real values of $u$.
\endpro
\Pf.
The statement follows from Theorems~\[Shapz] and~\[main2],
and Proposition~\[SMz+].
\epf

\myRefs

\ref\Key D
\by \Dri/
\paper Quantum groups
\inbook in Proceedings of the ICM, Berkeley, 1986
\yr 1987 \pages 798\>\~\>820 \publ \AMSa/
\endref

\ref\Key FFR
\by B\&Feigin, E\&Frenkel, and \Reshet/
\paper Gaudin model, \Ba/ and critical level
\jour \CMP/ \vol 166 \yr 1994 \issue 1 \pages 27\)\~\>62
\endref

\ref\Key FMTV
\by \Feld/, Ya\&Markov, \VT/ and \Varch/
\paper Differential \eq/s compatible with \KZ/ \eq/s
\jour Math\. Phys., Analysis and Geometry \vol 3 \yr 2000 \pages 139\>\~177
\endref

\ref\Key G
\by M\&Gaudin
\book La fonction d'onde de Bethe \yr 1983 \publ Masson \publaddr Paris
\endref

\ref\Key Ju
\by B\&Jur\v co
\paper Classical \YB/s and quantum integrable systems
\jour \JMP/ \vol 30 \yr 1989 \pages 1289\>\~1293
\endref

\ref\Key KBI
\by \Kor/, N\&M\&Bogoliubov and A\&G\&Izergin
\book Quantum inverse scattering method and correlation \fn/s
\yr 1993 \publ \CUP/ 
\endref

\ref\Key KR1
\by \Kulish/ and \Reshy/
\paper Diagonalization of $GL(n)$ invariant transfer-mat\-rices and quantum
\$N\}$-wave system (Lee model)
\jour \JPA/ \vol 15 \yr 1983 \pages L\)591\~\>L\)596
\endref

\ref\Key KR2
\by \Kulish/ and \Reshy/
\paper On \$GL_3$-invariant \sol/s of the \YB/ and associated quantum systems
\jour \JSM/ \vol 34 \yr 1986 \pages 1948\>\~1971
\endref

\ref\Key KS
\by \Kulish/ and \Skl/
\paper Quantum spectral transform method. Recent developments
\jour Lect. Notes in Phys. \vol 151 \yr 1982 \pages 61\~119
\endref

\ref\Key Ma
\by A\&Matsuo
\paper An application of Aomoto-Gelfand \hgeom/ \fn/s to the $SU(n)$
\KZv/ \eq/ \jour \CMP/ \vol 134 \yr 1990 \pages 65\>\~\)77
\endref

\ref\Key Mo
\by A\&Molev
\paper Yangians and their applications \inbook Handbook of Algebra, \volume 3
\yr 2003 \pages 907\)\~\>959 \publ North\)-Holland
\ifUS\adjustnext{\kern-4pt}\fi \publaddr Amsterdam
\endref

\ref\Key MNO
\by A\&Molev, \MN/, and G\&Olshanski
\paper Yangians and classical Lie algebras
\jour Russian Math\. Surveys \vol 51 \yr 1996 \pages 205\>\~\)282
\endref

\ref\Key MTV
\by E\&Mukhin, \VT/, and \Varch/
\paper The B.\,\,and~M\&Shapiro conjecture in real algebraic geometry and
the \Ba/ \jour Preprint \yr 2005 \pages 1\~17 \info\tt math.AG/0512299
\endref

\ref\Key MV1
\by E\&Mukhin and \Varch/
\paper Critical points of master functions and flag varieties
\jour Commun. Contemp. Math. \vol 6 \yr 2004 \issue 1 \pages 111\~163
\endref

\ref\Key MV2
\by E\&Mukhin and \Varch/
\paper Solutions to the \XXX/ type \Bae/s and flag varieties
\jour Central Europ. J. Math. \vol 1 \yr 2003 \pages 238\>\~\)271
\endref

\ref\Key N
\by \MN/
\paper Yangians and Capelli identities \jour \AMST/ \vol 181 \yr 1998
\pages 139\>\~163
\endref

\ref\Key NO
\by \MN/ and G\&Olshanski
\paper Bethe subalgebras in twisted Yangians
\jour \CMP/ \vol 178 \yr 1996 \issue 2 \pages 483\>\~\)506
\endref

\ref\Key RSV
\by R\&Rim\' anyi, L\&Stevens, and \Varch/
\paper Combinatorics of \raf/s and Poincar\'e\)-Birchoff\>-\]Witt expansions
of the canonical \>\$U(\ngm)$-\)valued differential form
\jour Annals Combin. \vol 9 \yr 2005 \issue 1 \pages 57\)\~\)74
\endref

\ref\Key RV
\by \Reshet/ and \Varch/
\paper Quasiclassical asymptotics of \sol/s to the \KZ/ \eq/s
\inbook Geometry, Topology, \ampersand~Physics \yr 1995 \pages 293\>\~\)322
\publ Internat. Press \publaddr Cambridge, MA
\endref

\ref\Key SV
\by \SchV/ and \Varn/
\paper Arrangements of hyperplanes and Lie algebras homology
\jour Invent. Math. \vol 106 \yr 1991 \pages 139\>\~194
\endref

\ref\Key T
\by D\&\]Talalaev
\paper Quantization of the Gaudin system \jour Preprint \yr 2004 \pages 1\~19
\adjustnext\uline \info\tt hep-th/0404153
\endref

\ref\Key TV1
\by \VT/ and \Varch/
\paper Jackson integral \rep/s for \sol/s to the quantized \KZv/ \eq/
\jour \LpMJ/ \vol 6 \yr 1994 \issue 2 \pages 90\>\~137
\endref

\ref\Key TV2
\by \VT/ and \Varch/
\paper Duality for \KZv/ and dynamical \eq/s
\jour Acta Appl. Math. \vol 73 \yr 2002 \issue 1\~-2 \pages 141\~154
\endref

\ref\Key TV3
\by \VT/ and \Varch/
\paper Dynamical \difl/ \eq/s compatible with \rat/ \qKZ/ \eq/s
\jour \LMP/ \vol 71 \yr 2005 \issue 2 \pages 101\~108
\endref

\endRefs

\bye